\documentclass[preprint,11pt]{elsarticle}
\usepackage[utf8]{inputenc}
\usepackage[T1]{fontenc}
\usepackage{lmodern}
\usepackage{amsmath}
\usepackage{amssymb}
\usepackage{amsfonts}
\usepackage{amsthm}
\usepackage{bbm}
\usepackage[
centering,
left = 2.3cm,
right = 2.3cm,
top = 2.5 cm,
bottom = 2.2 cm
]{geometry}
\usepackage{graphicx}
\usepackage{xcolor}
\usepackage{mathtools}
\usepackage{bm}
\usepackage{mathrsfs}
\usepackage{stmaryrd}
\usepackage{tikz}
\usepackage{pgfplots}
\pgfplotsset{compat=1.18}
\usetikzlibrary{pgfplots.groupplots}
\usepackage[
pdfencoding=auto,
colorlinks = true,
citecolor = blue,
linkcolor = black,
anchorcolor = blue,
urlcolor = blue
]{hyperref}

\makeatletter
\let\ps@pprintTitle\ps@empty
\makeatother

\allowdisplaybreaks[4]

\newcommand{\be}{\begin{equation}}
\newcommand{\bea}{\begin{equation}\begin{aligned}}
\newcommand{\beas}{\begin{equation*}\begin{aligned}}
\newcommand{\eeas}{\end{aligned}\end{equation*}}
\newcommand{\eea}{\end{aligned}\end{equation}}
\newcommand{\ee}{\end{equation}}

\renewcommand{\div}{{\rm div}}

\begin{document}
\begin{frontmatter}
\title{
On the Abels--Garcke--Gr\"un Model for Incompressible Two-Phase Flows with Unmatched Densities: 
Phase Separation and Rayleigh--Taylor Instability}

\author[FJ]{Fei Jiang}
\ead{jiangfei0591@163.com}
\author[FDU]{Maoyin Lv}
\ead{mylv22@m.fudan.edu.cn}
\author[FDU]{Hao Wu}
\ead{haowufd@fudan.edu.cn}

\address[FJ]{College of Mathematics and
Computer Science, Fuzhou University, Fuzhou,   China}
\address[FDU]{School of Mathematical Sciences, Fudan University, Shanghai, China}

\begin{abstract}
Instability phenomena are ubiquitous in multi-fluid flows. 
For binary fluid mixtures, phase separation may occur under specific conditions, 
in which a homogeneous state becomes unstable and spontaneously segregates into different phases, 
leading to the formation of distinct patterns. 
Furthermore, when a heavier fluid is superposed over a lighter one in a gravitational field, 
the fluid interface can be sensitive to perturbations that trigger the Rayleigh--Taylor (RT) instability, 
i.e., the lighter fluid penetrates upward into the heavier fluid against gravity. 
We investigate such instability phenomena and the inhibiting conditions for incompressible, 
immiscible two-phase flows with unmatched densities, 
within the framework of a diffuse-interface model, 
which is called the Abels--Garcke--Gr\"un (AGG) model. 
More precisely,
we adopt energy methods to establish several stability/instability conditions for the perturbed system of the AGG model, 
which is obtained by slightly perturbing the equilibrium profile $(\mathbf{0},\overline{\varphi})$, 
where $\overline{\varphi}=\overline{\varphi}(x_3)$ depends solely on the vertical coordinate $x_3$. 
Notably, for a constant equilibrium profile $\overline{\varphi}=\alpha\in(-1,1)$ satisfying $\Psi''(\alpha)<0$,
we establish a threshold for stability and instability of the phase-separation problem.
Similar stability results are also established for the perturbed system 
with $\overline{\varphi}$ satisfying the RT instability condition in the presence of gravity. 
In particular, if the phase-separation condition fails and $\partial_{x_3}\overline{\varphi}$ is a  constant, 
a larger value of $\partial_{x_3}\overline{\varphi}$, a greater interfacial thickness 
and a higher surface energy density contribute to inhibiting RT instability.
\end{abstract}
\end{frontmatter}


\newtheorem{thm}{Theorem}[section]
\newtheorem{lem}{Lemma}[section]
\newtheorem{pro}{Proposition}[section]
\newtheorem{cor}{Corollary}[section]
\newproof{pf}{Proof}
\newdefinition{rem}{Remark}[section]
\newtheorem{definition}{Definition}[section]

\section{Introduction}
\numberwithin{equation}{section}

Instability phenomena are widely observed in multi-fluid flows. 
For example, in binary fluid mixtures, phase separation may occur under certain conditions 
where a homogeneous state becomes unstable and spontaneously separates into distinct phases, 
leading to distinct patterns \cite{Anderson}. 
When a heavier fluid is superposed over a lighter one in a gravitational field, 
the fluid interface can be sensitive to perturbations that trigger the Rayleigh--Taylor (RT) instability, 
i.e., the lighter fluid penetrates upward into the heavier fluid against gravity \cite{Rayleigh,Taylor50}. 
These instabilities are closely related to morphological changes in the free interface between macroscopically immiscible fluid components.

Two main approaches are commonly employed to describe the interfacial evolution: 
the sharp-interface model and the diffuse-interface model. 
The former treats the interface as a lower-dimensional hypersurface, 
and the corresponding free-boundary problem consists of equations of motion in each fluid complemented by boundary conditions at the free surface. 
The sharp-interface model has been studied extensively and is successful in applications. 
However, it encounters difficulties in the analytic and numerical studies when topological changes such as pinch-off or coalescence occur. 
In contrast, the diffuse-interface model assumes that the interface has a finite thickness and allows for partial miscibility within a thin transition layer. 
This framework avoids explicit tracking of free interfaces and naturally accommodates large deformations and topological transitions, 
facilitating both mathematical analysis and numerical simulations \cite{Anderson,Gur96,LT1998}.

The Cahn--Hilliard (CH) equation provides a diffuse-interface description for phase separation in isothermal, isotropic binary mixtures, 
for instance, spinodal decomposition in alloys \cite{Cahn1965,Cahn1968,CH}. 
Its extension to two-phase flows is naturally realized by coupling the CH equation with the Navier--Stokes (NS) equations. 
The resulting Navier--Stokes--Cahn--Hilliard (NSCH) system turns out to be a promising and efficient tool 
for studying the complex dynamics of two incompressible Newtonian viscous fluids \cite{Anderson,LCSJPD179}. 
The classical ``Model H'' \cite{Gur96,HH77} addresses the specific case where the fluid constituents have equal constant densities. 
Subsequently, several extended diffuse-interface models have been developed for the general scenario of unmatched densities, 
incorporating different choices of the mean velocity of the mixture and constitutive assumptions; 
see, e.g., \cite{AHGHGGMMASTCDIMI,Aki14,Boyer02,LT1998,SYW13}. 
We refer to \cite{ten23} for a unified framework based on continuum mixture theory and thermodynamic principles, as well as a comparative study of existing models.

In this work, our aim is to investigate the stability conditions, under which  phase separation and RT instability do not occur, 
for incompressible immiscible two-phase flows with unmatched densities, 
based on the following thermodynamically consistent diffuse-interface model proposed by Abels, Garcke and Gr\"un in \cite{AHGHGGMMASTCDIMI}:
\begin{align}
	\begin{cases}
		\partial_t(\rho(\varphi)\mathbf{v})+\text{div}(\mathbf{v}\otimes(\rho(\varphi)\mathbf{v}+\widetilde{\mathbf{J}}))
        -\text{div}(2\nu(\varphi)\mathbb{D}\mathbf{v})+\nabla P=-\sigma l \div(\nabla\varphi\otimes\nabla\varphi)&\text{in }\Omega_T,\\
		\text{div}\, \mathbf{v}=0&\text{in }\Omega_T,\\
		\partial_t\varphi+\mathbf{v}\cdot\nabla\varphi=\text{div}(m(\varphi)\nabla\mu)&\text{in }\Omega_T,\\
		\mu=\sigma(l^{-1}\Psi'(\varphi)-l\Delta\varphi)&\text{in }\Omega_T,
	\end{cases}\label{AGG}
\end{align}
where $\Omega_T:=\Omega\times(0,T)$,  $\Omega\subset\mathbb{R}^3$ is a bounded domain 
with sufficiently smooth boundary $\Gamma:=\partial\Omega$ and $T\in(0,+\infty)$ is a given final time.
The state variables are the volume-averaged velocity $\mathbf{v}=\mathbf{v}(x,t):\Omega_{T}\to\mathbb{R}^3$, 
the pressure $P=P(x,t):\Omega_{T}\to\mathbb{R}$, and the phase-field variable $\varphi=\varphi(x,t):\Omega_{T}\to[-1,1]$. 
The volume‑averaged velocity of a mixture of incompressible fluids remains divergence‑free. In the momentum equation, 
the mean density $\rho$ and the mean viscosity coefficient $\nu$ are given by
\[
\rho(\varphi)=\rho_1\frac{1+\varphi}{2}+\rho_2\frac{1-\varphi}{2},\quad
\nu(\varphi)=\nu_1\frac{1+\varphi}{2}+\nu_2\frac{1-\varphi}{2},
\]
where $\rho_1$, $\rho_2>0$ and $\nu_1$, $\nu_2>0$, are the homogeneous densities and viscosity {coefficients}  of the two fluid components, resp. 
Without loss of generality, we assume that  $\rho_1>\rho_2>0$ in the subsequent analysis.
The operator $\mathbb{D}$ denotes the symmetric gradient, that is, $\mathbb{D}{\mathbf{v}} =(\nabla {\mathbf{v}}+\nabla {\mathbf{v}}^{\top})/2$,
and $\widetilde{\mathbf{J}}$ is the relative flux related to the diffusion of the fluid components such that
\[\widetilde{\mathbf{J}}= -\rho'(\varphi)m(\varphi)\nabla\mu =\frac{\rho_2-\rho_1}{2}m(\varphi)\nabla\mu.\]
The function $m(\varphi)\geq 0$ denotes the mobility coefficient for diffusion that may depend on $\varphi$.
The scalar function $\mu:\Omega_{T}\to\mathbb{R}$ is the chemical potential, which can be obtained by taking the variational derivative of the free energy
\[E_{\text{free}}(\varphi)=\int_\Omega \sigma\Big(\frac{l}{2}|\nabla\varphi|^2+\frac{1}{l}\Psi(\varphi)\Big)\,\mathrm{d}x,\]
that is, $\mu=\delta_\varphi E(\varphi)$.
Here, the constant $\sigma>0$ is related to the surface energy density and $l>0$ denotes the interfacial thickness.
The nonlinear function $\Psi$ denotes the bulk free energy density, which usually has a double-well structure, {see Appendix \ref{2608071516} for further details.}

In order to solve the hydrodynamic system \eqref{AGG} in a bounded domain, suitable boundary conditions are necessary. 
The classical boundary conditions are the no-slip boundary condition for the velocity field $\mathbf{v}$ 
and the homogeneous Neumann boundary conditions for the phase-field variable $\varphi$ and the chemical potential $\mu$, i.e.,
\begin{align}
    \mathbf{v}= \boldsymbol{0},\quad \partial_\mathbf{n}\varphi=\partial_\mathbf{n}\mu=0\quad \text{on }\Gamma_T,\label{classical-BC}
\end{align}
where $\Gamma_T:=\Gamma\times(0,T)$, $\mathbf{n}=\mathbf{n}(x)$ denotes the outward unit normal vector on $\Gamma$ 
and $\partial_\mathbf{n}$ represents the outward normal derivative on the boundary.
In addition, Dirichlet type boundary conditions for $(\varphi, \mu)$ are admissible:
\begin{align}
    {\varphi}={\varphi}_\mathrm{b},\quad {\mu}={\mu}_\mathrm{b} \quad\text{on }\Gamma_T,\label{classical-BC2}
\end{align}
where $\mu_\mathrm{b}$ is a harmonic function in $\Omega$; see, e.g., \cite{Elliott}.
We note that the Dirichlet boundary conditions for the CH equation have been widely used in the modeling of tumor growth; see, e.g., \cite{Chen2014,Frieboes,GarckeLam,Wise}.

For the AGG system \eqref{AGG} subject to suitable initial data $(\mathbf{v}_0,\varphi_0)$ and the boundary conditions \eqref{classical-BC}, 
extensive studies have been conducted in the literature. 
We refer to \cite{AHDDGHEJMFM,ADG-HP} for the existence of global weak solutions under non-degenerate or degenerate mobilities. 
The local well-posedness of strong solutions in three dimensions for polynomial-like potentials $\Psi$ was established in \cite{AW-JEE}.
In \cite{Gior21}, local well-posedness of strong solutions to the AGG system with logarithmic potential in two dimensions was proved; 
moreover, in the case of periodic boundary conditions, the strong solution was shown to exist globally.
The results in \cite{Giorgini2022IFB} further extended the local well-posedness from \cite{Gior21} to a three-dimensional bounded domain. 
We refer to \cite{AGG2024} for global regularity and asymptotic stabilization of global weak solutions, 
and to \cite{GLW26} for regularity criteria of energy equality and Lyapunov stability of energy-minimizing equilibrium. 
Recently, variants and generalizations of the AGG model have been extensively explored. 
For example, a multi-component version was considered in \cite{AGP-JMFM}, while an evolving surface formulation was investigated in \cite{AGP-CVPDE,AGP2024}. 
Nonlocal variants have been addressed in \cite{AT} for fractional free energies 
and in \cite{Frigeri,Frigeri2021,GalGraGiorPoia} for free energy with regular convolution kernels. 
The relationship between local and nonlocal AGG models was examined in \cite{AT2022}. 
Recent progress on the AGG model subject to a generalized Navier boundary condition for the fluid velocity 
and a dynamic boundary condition for the phase-field variable can be found in \cite{GGW19,GLW24,GK23}.

\subsection{Problem setting}
To better understand how the magnitude of the physical parameters influences {the motion} of the binary fluids, for simplicity, we assume that
$$
m(\varphi)\equiv m>0,\ \nu(\varphi)=\nu_1=\nu_2\equiv\nu>0,\text{ i.e., } m \text{ and } \nu \text{ are positive constants.}
$$
Then using the identity (cf. \cite{AHDDGHEJMFM})
\[ \text{div}(\nabla\varphi\otimes\nabla\varphi)=  \Delta\varphi\nabla\varphi+ \nabla\Big(\frac{|\nabla\varphi|^2}{2}\Big),\]
we can rewrite the system \eqref{AGG} in the presence of a uniform gravitational field as
\begin{equation}\label{01011x}
\begin{cases}
 \rho(\varphi) \partial_t \mathbf{v}+   \rho(\varphi)  ( \mathbf{v}\cdot\nabla)\mathbf{v}-
 m\rho'(\varphi)(\nabla  \mu\cdot\nabla) \mathbf{v}-\nu\Delta \mathbf{v}+\nabla p \\
\quad   =-\sigma {l}
  \Delta \varphi\nabla \varphi -\rho(\varphi) g \mathbf{e}_3&\text{in }\Omega_{T},\\
\mathrm{div}\,\mathbf{v} =0&\text{in }\Omega_{T},  \\
\partial_t \varphi+\mathbf{v}\cdot\nabla \varphi=
m\Delta  \mu  &\text{in }\Omega_{T},\\
\mu =\sigma(l^{-1}\Psi'(\varphi)-l\Delta\varphi)&\text{in }\Omega_{T},
\end{cases}
\end{equation}
where $p$ is {called} the modified hydrodynamic pressure, $g\geq0$ is the gravitational constant, 
$\mathbf{e}_3=(0,0,1)^{\top}$ denotes the vertical unit vector, 
and $-\rho(\varphi) g \mathbf{e}_3$ represents the gravitational force in the opposite direction of the $x_3$-axis.

To present a unified treatment of phase separation and the RT instability, we focus, in the subsequent analysis, 
on perturbations of a prescribed $x_3$-dependent equilibrium profile $\overline{\varphi}=\overline{\varphi}(x_3)$ for the phase‑field variable, 
which is independent of the horizontal coordinates $(x_1,x_2)$. 
To this end, let us consider an equilibrium state $(\boldsymbol{0}, \overline{\varphi})$, 
where $\overline{\varphi}=\overline{\varphi}(x_3)$ satisfies $-1<\overline{\varphi}<1$ and
\begin{equation}
\label{20171805}
 \partial_{3}^{2}(l^{-1}\Psi'(\overline{\varphi})- l\partial_3^2\overline{\varphi}) =0 \text{ in}\ \Omega,
\end{equation}
with nonhomogeneous Dirichlet boundary conditions
\begin{align}
    \overline{\varphi}= \varphi_\mathrm{b},\quad \overline{\mu}=\mu_\mathrm{b}\quad\text{on }\Gamma,\label{Dirichlet'}
\end{align}
or  homogeneous Neumann boundary conditions
\begin{align}
    \partial_\mathbf{n} \overline{\varphi}=\partial_\mathbf{n}\overline{\mu}=0\quad\text{on }\Gamma,\label{Neumann'}
\end{align}
where we have defined the equilibrium chemical potential by
\begin{align}
\overline{\mu}:=\sigma (l^{-1}\Psi'(\overline{\varphi})-l\partial_3^2\overline{\varphi}).\notag
\end{align}

\emph{From now on, we only consider the smooth solutions $\overline{\varphi}$, which solve \eqref{20171805} with  \eqref{Dirichlet'} or \eqref{Neumann'}}. 
We note that the existence of an admissible equilibrium profile $(\overline{\varphi}, \overline{\mu})$ depends on the boundary conditions, 
the nonlinear potential function $\Psi$, as well as the geometry of the spatial domain (see Appendix \ref{2608071516} for further discussions). 
This existence problem is of independent interest and is not pursued here. 
The equilibrium pressure $\overline{p}$ corresponding to $\overline{\varphi}$ is also independent of $(x_1,x_2)$ and can be defined by the relation
$$
\nabla \overline{p}=  -(\sigma {l}
 \partial_{3}^{2} \overline{\varphi} \,\partial_{3}\overline{\varphi}
 +\overline{\rho} g)\mathbf{e}_3,
$$
where $\overline{\rho}$
denotes the equilibrium density profile such that
$$
\overline{\rho}:= \varrho_1 \overline{\varphi}+ \varrho_2 \quad  \text{with}\quad \varrho_1=\frac{\rho_1-\rho_2}{2},\ \ \varrho_2=\frac{\rho_1+\rho_2}{2}.
$$

For convenience, we introduce new variables for perturbations around the equilibrium profile:
$$\mathbf{v}= \mathbf{v}-\boldsymbol{0},\quad \phi= \varphi-\overline{\varphi},
\quad  q =p-\overline{p}.$$
Then $(\mathbf{v},\phi,q)$ satisfies the following perturbed AGG system:
\begin{equation}\label{purterbation}
\begin{cases}
 (\varrho_1 \phi+\overline{\rho}  )\partial_t \mathbf{v}+  (\varrho_1\phi+\overline{\rho})( \mathbf{v}\cdot\nabla)\mathbf{v}
-\varrho_1 m (\nabla \mu(\phi+\overline{\varphi})\cdot \nabla) \mathbf{v}
  -  \nu \Delta  \mathbf{v} +\nabla q
 \\
\quad=-\sigma {l}
( \Delta \phi\nabla(\phi+ \overline{\varphi})+ \Delta\overline{\varphi}  \nabla \phi )-g \varrho_1 \phi \mathbf{e}_3 &\text{in }\Omega_{T},\\
\mathrm{div}\,\mathbf{v} =0&\text{in }\Omega_{T}, \\
 \partial_t\phi+  \mathbf{v} \cdot\nabla (\phi+\overline{\varphi})=m \Delta \widehat{\mu}&\text{in }\Omega_{T},\\
  \widehat{\mu}=\sigma(l^{-1}(\Psi'(\phi+\overline{\varphi})-\Psi'(\overline{\varphi}))-l\Delta\phi)&\text{in }\Omega_{T},
\end{cases}
\end{equation}
where $\mu(\phi+\overline{\varphi})=\sigma(l^{-1}\Psi'(\phi+\overline{\varphi})-l\Delta(\phi+\overline{\varphi}))$.
The system \eqref{purterbation} is subject to the initial conditions
\begin{align}
\mathbf{v}(\cdot,0)=\mathbf{v}_0,\quad \phi(\cdot,0)=\phi_0\quad \text{in } \Omega.\label{initial}
\end{align}
For the velocity field $\mathbf{v}$, we impose the classical no-slip boundary condition
\begin{align}
    \mathbf{v}=\boldsymbol{0}\quad\text{on }\Gamma_{T},\label{no-slip}
\end{align}
while for the phase-field variable $\phi$, we take either the homogeneous Dirichlet boundary conditions
\begin{align}
    \phi=\Delta\phi=0\quad\text{on } \Gamma_{T},\label{Dirichlet}
\end{align}
or 
 the homogeneous Neumann boundary conditions
\begin{align}
    \partial_{\mathbf{n}}\phi =\partial_{\mathbf{n}}\Delta\phi=0 \quad \text{on }\Gamma_{T}.\label{Neumann}
\end{align}
The boundary conditions \eqref{no-slip}--\eqref{Neumann} can be derived from \eqref{classical-BC}, \eqref{classical-BC2}, \eqref{Dirichlet'} and \eqref{Neumann'}, 
because the solution $(\varphi, \mu)$ and the equilibrium profile $(\overline{\varphi},\overline{\mu})$ satisfy the same boundary data. 
For convenience, we call the set of initial-boundary conditions \eqref{initial}--\eqref{Dirichlet} the ``\emph{initial-boundary condition $(\mathbf{D})$}'', 
and \eqref{initial}, \eqref{no-slip} and \eqref{Neumann} the ``\emph{initial-boundary condition $(\mathbf{N})$}''.

\subsection{Phase separation} 
We first study the phase separation process in a binary fluid mixture without the influence of the gravitational field, i.e., $g=0$. 
Phase separation is ubiquitous in materials science, engineering, and biology, manifesting across diverse systems from alloys, polymer blends to small organic molecules and proteins, 
where it governs the formation of micro-structures and underpins complex biochemical reactions. 
Under appropriate environmental conditions and compositional properties, a homogeneous mixture spontaneously segregates into distinct phases, 
each characterized by different concentrations of components. 
A typical example of interest is the so-called spinodal decomposition \cite{Cahn1965,Cahn1968}. 
Consider a constant equilibrium profile $\overline{\varphi}=\alpha\in(-1,1)$ within the spinodal region, 
that is, $\Psi''$ satisfies the phase-separation condition
\begin{equation}
\label{constant-ins}
\Psi''(\alpha)<0.
\end{equation}
In the spinodal region, a homogeneous mixture may become unstable to small compositional perturbations, potentially leading to phase separation.
We also call such a phenomenon the \textit{phase-separation instability} in this paper. 
In \cite{Cahn1961}, by performing a linear analysis of the CH equation, 
the author showed that solutions starting with initial data near a constant state in the spinodal region exhibit fine-grained decomposition. 
By contrast, the constant profile $\alpha$ is stable, provided that it lies in a metastable region of $\Psi$, 
i.e., the potential function $\Psi$ is convex near $\alpha$. 
These principles agree with the outcome of physical experiments and the rigorous mathematical analysis can be found in \cite{Bates-Fife,Grant,MW-1,MW-2}. 
A similar phenomenon is expected for fluid mixtures. 
Namely, the condition \eqref{constant-ins} can induce phase separation: when an equilibrium profile $(\mathbf{0},\alpha)$ is perturbed, 
the constituent fluid species tend to assemble, leading to the demixing of the two components. 
For incompressible two-phase flows under shear with matched densities $\rho_1=\rho_2$, 
the stability of the constant equilibrium profile has been investigated in \cite{Boyer}.

In this work, we consider the general case of binary fluid mixtures with unmatched densities.
For simplicity, we call the perturbed system \eqref{purterbation} with $g=0$ subject to initial-boundary condition $(\mathbf{D})$ the \textit{phase-separation problem}, 
since phase separation may occur in \eqref{purterbation} with $g=0$.
For a generic $x_3$-dependent equilibrium profile $(\mathbf{0},\overline{\varphi}(x_3))$, 
we provide stability conditions for inhibiting phase-separation instability, see Theorems \ref{thm-1} and \ref{thm-2}.
In addition, if $\overline{\varphi}=\alpha$ and $\Psi''(\alpha)<0$, where the constant $\alpha$ belongs to $(-1,1)$, 
we also provide an instability condition, under which there exist unstable solutions in the sense of Hadamard, see Theorem \ref{thm-3}. 
In this case, the sign of $l^2\lambda_1+\Psi''(\alpha)$ distinguishes exponential stability (positive sign) from nonlinear escape in $L^2$ (negative sign), where $\lambda_1$ is the first eigenvalue of the negative Laplacian subject to homogeneous Dirichlet boundary condition; see Remark \ref{constant-instability} for further discussions.

\subsection{Rayleigh--Taylor instability}
Next, we investigate the RT instability of a binary fluid mixture using the AGG model \eqref{purterbation} 
in the presence of a uniform gravitational field, i.e., $g>0$. 
In addition, we assume
\begin{equation}\label{0102}
\partial_{3}\overline{\varphi}|_{x_3=x^0_{3}}>0
\end{equation}
for some $x^0_{3}\in \{x_3~|~(x_1,x_2,x_3)^{\top}\in \Omega\}$.
For simplicity, we call condition \eqref{0102} the \textit{RT condition},
and the perturbed system \eqref{purterbation} subject to initial-boundary condition $(\mathbf{D})$ under the RT condition \eqref{0102} the \textit{CHRT problem}.

From the definition of $\overline{\rho}$ and the assumption $\rho_1>\rho_2>0$, 
we see that the equilibrium density profile also satisfies
\begin{align}
\partial_{3} \overline{\rho}|_{x_3=x^0_{3}}>0.\notag
\end{align}
The above condition implies that there is at least a subregion of $\Omega$ in which the density increases with $x_3$, 
i.e., the heavier fluid overlies the lighter fluid within the mixture. 
Due to the presence of gravity, this density distribution is unstable under small perturbations in the classical NS equations \cite{JFJSAM,JJN}, 
see also \cite{JFJSO2014,JFJSWWWOA,JJW,JFJWGCOSdd,FJWGCZXOE} for related results on viscoelastic fluids and magnetohydrodynamics.

The RT instability of two-phase flows has been extensively studied within the sharp-interface framework, 
including the stabilizing effects of surface tension; see, e.g., \cite{JJZ,PS,WT,WYJTIKCT,Wilke}. 
The diffuse-interface framework, which provides a convenient approach to describing complex interfacial dynamics, 
has also been employed to investigate the RT instability through linear analysis and numerical simulations; 
see, e.g., \cite{CMMV,LKK,ZTLH}. In particular, the studies in \cite{CMMV,ZTLH} employ the Boussinesq approximation and focus on small Atwood numbers, 
whereas the AGG model considered here allows for general unmatched densities without invoking this approximation
and incorporates an additional diffusive mass flux in the momentum equation.

Beyond the aforementioned studies, we mention some further developments in the mathematical analysis of RT instability for fluids with variable density; see, e.g., \cite{Aba2010,BM2017,GKS2021,HwangGuo2003,Zhou2021}. 
In particular, we recall some recent results for the case of vanishing mobility, i.e., $m=0$. The AGG system \eqref{01011x} then formally reduces to a gravity-driven Navier--Stokes--Korteweg (NSK) system:
\begin{align}
    \begin{cases}
        \partial_t\rho+\mathbf{v}\cdot\nabla\rho=0&\text{in}\ \Omega_{T},\\[1mm]
        \rho\partial_t\mathbf{v}+\rho(\mathbf{v}\cdot\nabla)\mathbf{v}-\nu\Delta\mathbf{v}+\nabla q
        =-\kappa\Delta\rho\nabla\rho-\rho g\mathbf{e}_3&\text{in}\ \Omega_{T},\\[1mm]
        \text{div}\,\mathbf{v}=0&\text{in}\ \Omega_{T},
    \end{cases}\label{NSK}
\end{align}
where $\kappa=\sigma l \varrho_1^{-2}$ can be considered the capillary coefficient. 
The NSK system \eqref{NSK} describes the motion of a single viscous, incompressible, capillary fluid in the presence of a uniform gravitational field. 
There are various works that study the local and global existence of solutions to the system \eqref{NSK}, 
see \cite{Burtea,SBGLR,TanWang,Wang2017,YYZ} and the references cited therein.
For the RT problem of \eqref{NSK}, 
the impact of capillarity on the RT instability in a horizontally periodic domain $(2\pi L\mathbb{T})^2\times\mathbb{R}$ was analyzed in \cite{ZhangZP}. 
The author showed that for any $\kappa>0$, there exists an unstable solution to the linearized RT problem, 
but the nonlinear instability can only be proved when $\kappa$ is sufficiently small.
Later, \cite{ZTW2023} studied the system \eqref{NSK} in a horizontal periodic domain with infinite height $2\pi\mathbb{T}\times\mathbb{R}$ and improved the result in \cite{ZhangZP}. 
The authors showed that RT instability can occur in fluids with capillary effects for any capillary coefficient $\kappa>0$ if the critical capillary coefficient is infinite.
In \cite{Nguyen}, the author applied Guo--Strauss’s bootstrap instability method (cf. \cite{GuoStrauss,GYSWIC}) to establish the nonlinear RT instability result 
with no-slip boundary condition under some sharp instability condition in a horizontally periodic domain $(2\pi L\mathbb{T})^2\times(0,1)$.
Concerning a smooth bounded domain $\Omega\subset\mathbb{R}^3$, the RT instability for the linear NSK system \eqref{NSK} has been studied in \cite{LZ2023SIMA}.
Moreover, we mention \cite{JLZ2023,JZZ2024} on the RT instability in the nonhomogeneous incompressible NSK equations 
in a horizontally periodic slab domain with Navier (slip) boundary conditions.

In this work, we provide stability conditions, under which the RT instability  does not occur in the \textit{CHRT problem}, see Theorem \ref{thm:0104n}. 
In addition, if $\overline{\varphi}$ is a linear function in $x_3$ and $\Psi''(\overline{\varphi})\geqslant 0$ (i.e., the phase-separation condition fails), 
we have a stability condition, which presents that a higher surface energy density, a greater interfacial thickness 
and a bigger changing rate of the equilibrium profile $\overline{\varphi}$ can directly contribute to inhibiting the RT instability, see Theorem \ref{thm:0104nbx}.  
We recall that surface tension (the tendency of the surface of a liquid to minimize its area due to cohesive forces between molecules) 
has inhibiting effect in the development of RT instability in stratified fluids \cite{GYTI2,WYJTIKCT}. 
Since surface tension is closely related to the surface energy density (cf. \cite{LCSJPD179} for diffuse-interface models), 
we are not surprised that the surface energy density should also present inhibiting effect in the development of RT instability 
within the current framework of the AGG model \eqref{purterbation}.

Finally, we point out that although the RT instability has been extensively studied within the sharp-interface framework, 
the corresponding methods cannot be directly applied to the AGG model. 
A major difficulty arises from the Cahn--Hilliard diffusion, 
which prevents the phase-field perturbation from being directly expressed in terms of the velocity perturbation. 
This complicates the construction of exponentially growing modes. 
Consequently, establishing nonlinear RT instability for the AGG model remains an important open problem.

\medskip
\noindent\textit{Outline of this paper.} 
In Section \ref{2608141458}, we first introduce some notation and basic analytical tools used in this paper, 
and then we state our main results, followed by a sketch of the proofs. 
Section \ref{26081414581} is devoted to the proofs of stability results for the phase-separation problem without the influence of gravity. 
In Section \ref{26081414582}, we further construct unstable solutions to the phase-separation problem and characterize the largest growth rate. 
In Section \ref{26081414583}, we prove the stability results for the {CHRT problem}. In Appendix \ref{2608071516}, we discuss the admissible equilibrium profiles.

\section{Main Results}\label{2608141458}

In this section, we first introduce the notation, basic assumptions, and some preliminary results that are useful in the subsequent analysis. 
Then, we state the main results in detail and sketch their proofs.

\subsection{Preliminaries}

The notation $a\lesssim b$ means that $a\leq cb$ for some constant
$c>0$ depending only on the domain $\Omega$ and the known physical
parameters and functions appearing in \eqref{purterbation}, such as $\rho_1$,
$\rho_2$, $g$, $\overline{\varphi}$, $l$, $m$, $\sigma$, $\nu$, and $\Psi$. We also use $C$,
$c_i$ ($i=1,2$), and $C_{j}^{\mathrm{d}}$ $(1\leq j\leq 4)$ to denote positive constants with the same
dependence. Generic constants may vary from line to line, whereas
indexed constants are fixed once specified.
Moreover, for $i\in\{1,2,3\}$ and $j\in\{0,1,2\}$, the symbols $\partial_i^j$ and $\partial_t^j$ 
denote the $j$-th partial derivatives with respect to $x_i$ and $t$, resp.

For any $p\in [1,+\infty]$, $L^p:=L^p(\Omega)$ denotes the Lebesgue space with norm $\|\cdot\|_{L^p}$. 
For $k\in\mathbb{Z}^+$, $p\in[1,+\infty]$, $W^{k,p}:=W^{k,p}(\Omega)$ denotes the Sobolev space with norm $\|\cdot\|_{W^{k,p}}$. 
When $p=2$, we denote $H^k=W^{k,2}$ with the norm $\|\cdot\|_k$. We denote the dual space of Banach space $X$ by $X'$, and the duality pairing is denoted by $\langle\cdot,\cdot\rangle_{X',X}$. 
In addition, let $\mathbf{X}$ denote the generic space of vectors or matrices whose entries belong to $X$.  
Next, let us introduce the Hilbert spaces of solenoidal vector-valued functions
 \begin{align*} 
 	&{L}_\sigma^2=\{\mathbf{u}\in\mathbf{L}^2:\text{div}\,\mathbf{u}=0\text{ in }\Omega,\ \mathbf{u}\cdot\mathbf{n}=0\text{ on }\Gamma\},\\
 	&{H}_\sigma^1=\{\mathbf{u}\in\mathbf{H}^1:\text{div}\,\mathbf{u}=0\text{ in }\Omega,\ \mathbf{u}=\boldsymbol{0}\text{ on }\Gamma\}.
 \end{align*}    For convenience, we define ${H}_\sigma^2={H}_\sigma^1\cap \mathbf{H}^2$.
Finally, we introduce some function spaces that will be frequently used in the subsequent analysis:
 \begin{align}
    &H_0^1:=\big\{u\in H^1: u=0 \ \text{on }\Gamma\big\},\  \underline{H}_{0}^1:=\Big\{u\in H^1: \int u\,\mathrm{d}x=0\Big\},\notag\\
 	&{H}_\mathrm{D}^3:=\big\{u\in H^3:  u=\Delta u =0\ \text{on }\Gamma\big\}, \ 
    \mathcal{A}_\mathrm{D}: =\big\{u\in {H}_\mathrm{D}^3:\|\nabla u\|_{0} =1\big\}, \notag\\
 	&{H}_\mathrm{N}^3:=
 	\Big\{u\in H^3:\partial_\mathbf{n} u=0\  \text{on }\Gamma\ \text{and}\ \int u\,\mathrm{d}x=0 \Big\},\notag
 	\\
 	&{H}_\mathrm{N}^5:=\Big\{ u \in H^5:  \partial_\mathbf{n} u= \partial_\mathbf{n} \Delta u =0\  \text{on }\Gamma\ \text{and}\ \int u\,\mathrm{d}x=0\Big\}.
 	\nonumber
 	\end{align}
Hereafter, we use the simplified notation for the integral over $\Omega$ such that $\int:=\int_\Omega$.
Thanks to Poincar\'e's inequality and the Poincar\'{e}--Wirtinger inequality (see, e.g., \cite[Lemmas 1.42 and 1.43]{NASII04}), 
it holds
 \begin{align}
 	\label{250508-Poincare}
 	\|u\|_{0} \lesssim  \|\nabla u \|_{0}, \  \forall\, u\in H_0^1\ \text{or}\ u\in \underline{H}_{0}^1.
 \end{align}
 We also recall 
 \begin{itemize}
   \item the product estimates in three dimensions:
 	\begin{align}
 		\|uv\|_{i}&\lesssim
 	\begin{cases}
 			\|u\|_{1}\|v\|_{1} & \hbox{for }i=0,  \\  
 			\|u\|_{i}\|v\|_{2}  & \hbox{for }1\leq i\leq 2, \\ 
 		\|u\|_{2}\|v\|_{i}+\|u\|_{i}\|v\|_{2}& \hbox{for }3\leq i\leq 5.
 		\end{cases} 
        \label{product}  
 \end{align}
   \item  the interpolation inequalities:  
 	for any given $0\leq j<i$, $\epsilon>0$,
 	there exist a constant $c_1(i,j)>0$ depending on $\Omega$, $i$ and $j$, and a constant 
 $\widetilde{c}_1( i,j,\epsilon)>0$ further depending on $\epsilon$ such that
 \begin{align}
        \|u\|_{j}\leq c_1(i,j)\|u\|_{0}^{1-\frac{j}{i}}\|u\|_{i}^{\frac{j}{i}}\leq  \epsilon\|u\|_{i} + \widetilde{c}_1( i,j,\epsilon)\|u\|_{0}.\label{201708190532}
 	\end{align} 
 \end{itemize}
 It should be noted that    we have assumed that the norms on the right-hand side
of \eqref{product}  and \eqref{201708190532} are finite.

Throughout this paper, we make the following basic assumptions:
\begin{itemize}
    \item [$(\mathbf{A0})$] $\Omega\subset\mathbb{R}^3$ is a $C^7$-smooth bounded domain with boundary $\Gamma:=\partial\Omega$.
    \item [$(\mathbf{A1})$] The potential function $\Psi\in C^8(-1,1)$.
    \item [$(\mathbf{A2})$] The equilibrium profile $\overline{\varphi}=\overline{\varphi}(x_3)\in C^5(\overline{\Omega})$ with $\|\overline{\varphi}\|_{L^\infty}<1$. 
    The admissible equilibrium profile $\overline{\varphi}$ satisfies the equation \eqref{20171805} together with either the boundary condition \eqref{Dirichlet'} or  \eqref{Neumann'}.
\end{itemize} 

Finally, we recall the following result on the unique solvability of local strong solutions to the perturbed system \eqref{purterbation} 
with the initial-boundary condition $(\mathbf{D})$ or $(\mathbf{N})$ (see  \cite[Theorem 2.1]{Lv2025}).
\begin{pro} \label{pro:0401n20160203}
Suppose that $g\geq 0$ and the assumptions $(\mathbf{A0})$--$(\mathbf{A2})$ are satisfied.

(1) For any given initial data $( \mathbf{v}_0,\phi_0 )\in {H}^2_\sigma\times (H_{\mathrm{D}}^3\cap H^5)$, satisfying
\begin{align}
	\label{250123}
		&  \|\phi_0+\overline{\varphi}\|_{L^\infty}< 1,
\end{align}
and the compatibility condition
\begin{align}
    \Delta\widehat{\mu}(0):=\sigma l^{-1}\Delta(\Psi'(\phi_0+\overline{\varphi})-\Psi'(\overline{\varphi}))-\sigma l\Delta^2\phi_0=0\quad\text{a.e. on }\Gamma,\label{compatibility}
\end{align}
there exists a time $T_0>0$ such that the perturbed system \eqref{purterbation} subject to the initial-boundary condition $(\mathbf{D})$ 
possesses a unique strong solution $(\mathbf{v},\phi,q)$ on $[0,T_0]$, satisfying the following regularity properties
\begin{align*}
	&\mathbf{v}\in C([0,T_0];{H}_\sigma^2)\cap L^2(0,T_0;\mathbf{H}^3),\  \partial_t \mathbf{v}\in C([0,T_0];{L}_\sigma^2)\cap L^2(0,T_0;{H}_\sigma^1),\\
	&\phi\in C([0,T_0];H_{\mathrm{D}}^3\cap H^5)\cap L^2(0,T_0;H^7),\  \partial_t\phi\in C([0,T_0];H_0^1)\cap L^2(0,T_0;H_{\mathrm{D}}^3),\\
	&\partial_{t}(\rho(\phi+\overline{\varphi})\partial_t\mathbf{v})\in L^2(0,T_0;({H}_\sigma^1)'),\  q\in C([0,T_0];H^1)\cap L^2(0,T_0;H^2).
\end{align*}
Moreover, it holds
$$0< \inf_{x\in \overline{\Omega}}\{ \varrho_1( \phi(x,t)+\overline{\varphi}(x))+\varrho_2\} \text{ and } \int {q}(x,t)\,\mathrm{d}x=0,\ \forall\,t\in (0,T_0).
$$

(2) Similar results hold for the perturbed system \eqref{purterbation} subject to the initial-boundary condition $(\mathbf{N})$, 
provided that the initial data satisfy $( \mathbf{v}_0,\phi_0 )\in{H}^2_\sigma\times {H}_\mathrm{N}^5$ and \eqref{250123}. Correspondingly, the regularity of $\phi$ should be modified as
\[
\phi\in C([0,T_0];H_\mathrm{N}^5)\cap L^2(0,T_0;H^7),\ \partial_t\phi\in C([0,T_0];\underline{H}_{0}^1)\cap L^2(0,T_0;H_\mathrm{N}^3).
\]
\end{pro}

\begin{rem}
The compatibility condition \eqref{compatibility} is necessary in the Dirichlet case, 
because $\partial_t\phi|_{t=0}=m\Delta\widehat{\mu}(0)-\mathbf{v}_0\cdot\nabla(\phi_0+\overline{\varphi})$ belongs to $H_0^1$ and $\mathbf{v}_0$ has zero trace.
Although this condition was omitted in \cite[Theorem 2.1]{Lv2025}, the proof therein remains valid once it is imposed. 
Moreover, the class of admissible initial data is nonempty. 
Indeed, one may choose $\phi_0=\epsilon\zeta$ with $\zeta\in C_c^\infty(\Omega)$ and sufficiently small $\epsilon>0$, 
so that both \eqref{250123} and \eqref{compatibility} are satisfied. 
For constant equilibrium profiles, Proposition \ref{modified-initial-pro} 
further provides compatible initial data associated with an unstable linear mode.
\end{rem}

\subsection{Stability and instability for the phase-separation problem}
To address the stability result for the phase-separation problem, we define the following energy and dissipation functionals:
    \begin{align}
		&\mathcal{E}(t):=\|\mathbf{v}(t) \|_{2}^2+\|\partial_t \mathbf{v}(t) \|_{0}^2+\| \phi(t)\|_{5}^2+\|\partial_t\phi(t)  \|_{1}^2,
        \label{E}\\
		&\mathcal{D}(t):=\|\mathbf{v}(t)\|_{3}^2+\|\partial_t\mathbf{v}(t)\|_{1}^2+\|\phi(t)\|_{7}^2+\|\partial_t\phi(t)\|_{3}^2.
        \label{D}
	\end{align}
We also introduce the following \textit{competitive energy} between the destabilizing and stabilizing factors:
\begin{align}
& E_\mathrm{d}(u):=m\sigma\left( l\|\nabla\Delta u\|_{0}^2+ l^{-1}\int \Delta(\Psi''(\overline{\varphi})u)\Delta u\,\mathrm{d}x\right),
\label{E-d}
\end{align}
where $u$ is a scalar function. 
Our first result reads as follows.

\begin{thm}\label{thm-1}
Suppose that $g=0$ and the assumptions $(\mathbf{A0})$--$(\mathbf{A2})$ are satisfied.
Assume, in addition, that one of the following conditions is satisfied:
\begin{itemize}
     \item [(a)] it holds that
    \begin{align}
    \omega_\mathrm{d}:=\sup_{u\in H_{\mathrm{D}}^3,\,\|\nabla\Delta u\|_{0}\neq 0}\frac{-\int \Delta(\Psi''(\overline{\varphi})u)\Delta u\,\mathrm{d}x}{\|\nabla\Delta u\|_{0}^2}<l^2;
    \label{assume1}
\end{align}
    \item [(b)] it holds that
    \begin{align}
        \int \Psi''(\overline{\varphi})|\nabla u|^2\,\mathrm{d}x-\frac{1}{2}\int\partial_{3}^{2}(\Psi''(\overline{\varphi}))u^2\,\mathrm{d}x\geq0, \ \  \text{for all }\,u\in H_0^1.
        \label{de-assume-1}
    \end{align}
\end{itemize}
Then, for sufficiently large viscosity $\nu$ or mobility $m$, there exists a constant $\varepsilon\in (0,1)$ such that, 
for any initial data $(\mathbf{v}_0,\phi_0)\in {H}_\sigma^2\times(H_{\mathrm{D}}^3\cap H^5) $ 
satisfying the compatibility condition \eqref{compatibility} and 
$\sqrt{\|\mathbf{v}_0\|_{2}^2+\|\phi_0 \|_{5}^2}<\varepsilon$,
the phase-separation problem possesses a unique global strong solution
$$
(\mathbf{v},\phi)\in C([0,+\infty);{H}_\sigma^2\times (H_{\mathrm{D}}^3\cap H^5)).
$$
Moreover, the solution satisfies the following exponential decay estimate:
\begin{equation}
\label{201708111550n}
e^{c_1 t} \mathcal{E}(t)+ \int_0^t e^{c_1 \tau} \mathcal{D}( \tau)\,\mathrm{d}\tau\leq c_2( {\|\mathbf{v}_0\|_{2}^2}+\|\phi_0\|_{5}^2),\ \forall\, t\geq 0.
\end{equation}
\end{thm}

Under the stability conditions in the above theorem, all sufficiently small compatible perturbations decay exponentially to the prescribed equilibrium.
Due to the generality of $\Psi$ in Theorem \ref{thm-1},  
we need an extra condition, i.e., large viscosity coefficient $\nu$ or large mobility coefficient $m$,
to establish the stability of equilibrium profile $(\mathbf{0},\overline{\varphi}(x_3))$. However,
 if the equilibrium profile $\overline{\varphi}=\overline{\varphi}(x_3)$ is linear in $x_3$, that is,
 \begin{align}
     \partial_{3}\overline{\varphi}\ \text{ is a  constant},\label{constant-de}
 \end{align}
the additional restriction can be removed. More precisely, we have the following stability result.
\begin{thm}\label{thm-2}
Suppose that $g=0$ and the assumptions $(\mathbf{A0})$--$(\mathbf{A2})$ are satisfied.
Assume, in addition, $\overline{\varphi}$ satisfies \eqref{constant-de}.
If the following stability condition is fulfilled
 \begin{equation}
\label{201708101624nnnn-de}
{l^2}> -{\Psi''}(\overline{\varphi})\varpi_\mathrm{d} \mbox{ with }
\varpi_\mathrm{d}:=\sup_{u\in H_{\mathrm{D}}^3,\,\|\nabla\Delta u\|_{0}\neq 0 }\frac{\|\Delta u\|_{0}^2}{\|\nabla\Delta u\|_{0}^2},
\end{equation}
then there exists a constant $\varepsilon>0$ such that, for any initial data $(\mathbf{v}_0,\phi_0)
\in {H}_\sigma^2\times(H_{\mathrm{D}}^3\cap H^5) $ satisfying the compatibility condition \eqref{compatibility} and $\sqrt{\|\mathbf{v}_0\|_{2}^2+\|\phi_0 \|_{5}^2}<\varepsilon$, 
the phase-separation problem admits a unique global strong solution 
$(\mathbf{v},\phi) \in C([0,+\infty);{H}_\sigma^2\times (H_{\mathrm{D}}^3\cap H^5))$ that satisfies the exponential decay estimate \eqref{201708111550n}.
\end{thm}
\begin{rem}
Recalling \eqref{20171805}, we see that 
       \begin{align}
        \sigma (l^{-1}\Psi'(\overline{\varphi})-l\partial_3^2\overline{\varphi})=a x_3+b\mbox{ for some constants }a\mbox{ and }b,\notag
        \end{align}
        which implies that 
         \begin{align}
         \Psi''(\overline{\varphi})\partial_3\overline{\varphi}=
        l^2\partial_3^3\overline{\varphi} +
        la/\sigma.\notag
        \end{align} 
Hence, it is easy to see that
\begin{align}
\Psi''(\overline{\varphi})\partial_3 \overline{\varphi}= la/\sigma \notag,
\end{align}
if \eqref{constant-de} holds. Consequently,
\begin{align}
\Psi''(\overline{\varphi})\mbox{ in \eqref{201708101624nnnn-de} is also a constant}.\label{quadratic}
\end{align}
In the subsequent analysis, whenever the assumption \eqref{constant-de} holds, we write the constant $\Psi''=\Psi''(\overline{\varphi})$ for simplicity.
\end{rem}

\begin{rem}
\label{sharp-poincare}
It is well known that
 \begin{align}
 	\|u\|_{0} \leq 1/\sqrt{\lambda_1} \|\nabla u \|_{0}, \ \forall\, u\in H_0^1,\notag
 \end{align}
 where $\lambda_1$ is the first eigenvalue of the negative Laplacian subject to homogeneous Dirichlet boundary condition.
In this case, we claim that $\varpi_\mathrm{d}=\lambda_1^{-1}$.
On the one hand, Poincar\'e's inequality implies that $\varpi_\mathrm{d}\leq \lambda_1^{-1}$.
On the other hand, $\|\Delta \xi_1\|_0^2/\|\nabla\Delta \xi_1\|_0^2=\lambda_1^{-1}$, where $\xi_1$ is the eigenfunction corresponding to $\lambda_1$.
This implies that $\varpi_{\mathrm{d}}= \lambda_1^{-1}$.
If $\Omega$ is also a convex domain,
we have $\varpi_\mathrm{d}\leq d_\Omega^2/\pi^2$, where $d_\Omega$ is the diameter of the domain $\Omega$.
\end{rem}

\begin{rem}
    \rm
The conclusions of Theorems~2.1--2.2 remain valid for the perturbed
system~\eqref{purterbation} subject to the initial-boundary condition
{(N)}, with the Dirichlet spaces $H_0^1$, $H_{\mathrm{D}}^3$ and $H_{\mathrm{D}}^3\cap H^5$ replaced by the corresponding
zero-mean Neumann spaces $\underline{H}_0^1$, $H_\mathrm{N}^3$ and $H_\mathrm{N}^5$, resp. In particular, the constants $\omega_{\mathrm{d}}$ and
$\varpi_{\mathrm{d}}$ should be replaced by their Neumann counterparts $\omega_{\mathrm{n}}$ and $\varpi_{\mathrm{n}}$:
\[
\omega_{\mathrm{n}}
    :=
    \sup_{u\in H_{\mathrm{N}}^3,\ 
    \|\nabla\Delta u\|_0\neq0}
    \frac{
    -\int
    \Delta\bigl(\Psi''(\overline{\varphi})u\bigr)
    \Delta u\,\mathrm{d}x}
    {\|\nabla\Delta u\|_0^2},\qquad
    \varpi_{\mathrm{n}}
    :=
    \sup_{u\in H_{\mathrm{N}}^3,\ 
    \|\nabla\Delta u\|_0\neq0}
    \frac{\|\Delta u\|_0^2}
         {\|\nabla\Delta u\|_0^2}.
\]
An argument similar to that in Remark \ref{sharp-poincare} yields $\varpi_{\mathrm{n}}=1/\lambda_{1,\mathrm{N}}^+$, 
where $\lambda_{1,\mathrm{N}}^+$ is the first positive Neumann eigenvalue.
In this case, the proofs of Theorems~2.1 and~2.2 can be adapted to the Neumann case by working in the corresponding mean-zero spaces and using the Poincar\'e--Wirtinger inequality.
Notice that in a general bounded domain, conditions \eqref{Neumann'} and \eqref{constant-de} typically force the equilibrium profile to be constant.
\end{rem}

We next establish a nonlinear instability result for the phase-separation problem.

\begin{thm} \label{thm-3}
Suppose that $g=0$, the assumptions $(\mathbf{A0})$ and $(\mathbf{A1})$ are satisfied, 
$\overline{\varphi}\equiv\alpha\in(-1,1)$ with $\Psi''(\alpha)<0$, and the following instability condition holds:
\begin{align}
\varpi_\mathrm{d} > \varpi:=-{l^2}/{\Psi''(\alpha)},
    \label{unstable}
\end{align}
where $\varpi_\mathrm{d}$ is defined in \eqref{201708101624nnnn-de}.
Then the equilibrium profile $(\boldsymbol{0},\alpha)$
 is unstable to the phase-separation problem. More precisely, there exist positive constants $\Lambda_{\mathrm{M}}$,
  $\epsilon_0$, $\delta_0$ and functions ${\psi}$, $\phi^{\mathrm{r}}\in H_{\mathrm{D}}^3\cap H^5$, 
  depending only on the known physical parameters and  functions in the perturbed system \eqref{purterbation} with $g=0$  and 
  the domain $\Omega$,
such that for any given $\delta\in (0,\delta_0)$ and the initial data
$(\mathbf{v}_0, \phi_0):=(\boldsymbol{0},\delta{\psi}+\delta^2\phi^{\mathrm{r}} )$, 
 the phase-separation problem admits a local-in-time strong solution $(\mathbf{v},\phi)\in C([0,T_{\mathrm{max}} ); {H}_\sigma^2\times (H_{\mathrm{D}}^3\cap H^5))$, but
\begin{equation}
\|  \phi(T^\delta )\|_{0} \geq \frac{ \epsilon_0}{2}\|\psi\|_0 \ \text{ for some escape time } T^\delta:=\Lambda_{\mathrm{M}}^{-1}\mathrm{ln}({\epsilon_0}{ \delta^{-1}})\in  (0,T_{\mathrm{max}}).
\notag
\end{equation}
Moreover, the function $\phi^{\mathrm{r}}\in H_{\mathrm{D}}^3\cap H^5$ is 
such that the initial datum $\phi_0=\delta\psi+\delta^2\phi^{\mathrm{r}}$ satisfies the compatibility condition \eqref{compatibility}.
\end{thm}

\begin{rem}
\label{constant-instability}
In view of Theorems \ref{thm-2} and \ref{thm-3}, 
if the equilibrium profile is constant (i.e., $\overline{\varphi}=\alpha\in(-1,1)$) with $\Psi''(\alpha)<0$, 
 the sign of $l^2\lambda_1+\Psi''(\alpha)$ distinguishes exponential stability (positive sign) from nonlinear escape in $L^2$ (negative sign).
 At the critical threshold $l^2\lambda_1+\Psi''(\alpha)=0$, the stability and instability analysis becomes more delicate.
 On the one hand, the maximal growth rate $\Lambda_{\mathrm{M}}:=\sup_{u\in \mathcal{A}_{\mathrm{D}}}(-E_\mathrm{d}(u))$ of the linearized Cahn--Hilliard equation vanishes. 
 Indeed, for any $u\in\mathcal{A}_{\mathrm{D}}$, write $u=\sum_{j}a_j \xi_j$, where \(\{\xi_j\}_{j\geq1}\) is an orthonormal Dirichlet
eigenbasis with eigenvalues \(\{\lambda_j\}_{j\geq1}\). 
 Then $E_{\mathrm{d}}(u)=m\sigma \sum_j (\lambda_j a_j^2)\lambda_j(\Psi''(\alpha) l^{-1}+l \lambda_j)$.
Since $\Psi''(\alpha)=-l^2\lambda_1$ and $\lambda_1\leq \lambda_j$ for all $j\geq2$, it follows that
\[
-E_{\mathrm d}(u)
=m\sigma l\sum_{j\geq1}\lambda_j^2a_j^2
(\lambda_1-\lambda_j)\leq0
\qquad\forall\,u\in\mathcal A_{\mathrm D}.
\]
Thus, $\Lambda_{\mathrm M}\leq0$. Moreover,
$\xi_1/\sqrt{\lambda_1}\in\mathcal A_{\mathrm D}$ and
$E_{\mathrm d}(\xi_1/\sqrt{\lambda_1})=0$.
Hence, $\Lambda_{\mathrm M}=0$, and the construction of exponentially growing solutions in Proposition \ref{exist-unstable-linear}, 
which is essential to our proof of nonlinear instability, is no longer applicable.
On the other hand, the energy functional $E_{\mathrm d}$ loses its coercivity, since $E_{\mathrm d}(\xi_1)=0$. 
Thus, the estimate \eqref{201708130943n}, which plays a crucial role in our exponential stability analysis, fails. 
In fact, even asymptotic stability cannot be guaranteed in general. To see this, suppose that
$\Psi(r)=-l^2\lambda_1(r-\alpha)^2/2$
in a neighborhood of $\alpha$. 
Then, for any sufficiently small nonzero constant $c$, set
$\mathbf v=\mathbf{0}$, $\phi=c\xi_1$, $q=\sigma l\lambda_1\phi^2/2$.
These are nontrivial stationary solutions of the perturbed system. 
Consequently, the equilibrium $(\mathbf 0,\alpha)$ is not asymptotically stable for this choice of potential.
\end{rem}

\begin{rem}\rm
For the Neumann counterpart of Theorem~2.3, the instability condition
becomes $\varpi_{\mathrm{n}}>\varpi$. 
The main difference is the construction of the unstable mode $\widetilde{\psi}$, which can be found in Remark \ref{Neumann-unstable-mode}.
Unlike the Dirichlet case, no correction of initial datum of the linear solution 
is required, since the local well-posedness result in Proposition \ref{pro:0401n20160203}-(2) does not
require compatibility conditions. Thus, one may
directly choose $(\mathbf{v}_0^\delta,\phi_0^\delta)=(\boldsymbol{0},\delta\widetilde{\psi})$ as the initial data of the nonlinear phase-separation problem. 
The remaining nonlinear instability argument follows as in the proof of Theorem~2.3, with the initial error equal to zero.
\end{rem}

\subsection{Stability results of inhibiting the RT instability}
Now we introduce our first stability result for the CHRT problem in the presence of a uniform gravitational field.  
\begin{thm}\label{thm:0104n}
Suppose that $g>0$ and the assumptions $(\mathbf{A0})$--$(\mathbf{A2})$ are satisfied.
Assume, in addition, that the stability condition \eqref{assume1} or \eqref{de-assume-1} is satisfied.
Then, for sufficiently large viscosity $\nu$ or mobility $m$, there exists a constant $\varepsilon>0$ such that, 
for any initial data $(\mathbf{v}_0,\phi_0)\in {H}_\sigma^2\times(H_{\mathrm{D}}^3\cap H^5) $ satisfying the compatibility condition \eqref{compatibility} and
$\sqrt{\|\mathbf{v}_0\|_{2}^2+\|\phi_0 \|_{5}^2}<\varepsilon$, the CHRT problem admits a unique global strong solution
$$
(\mathbf{v},\phi)\in C([0,+\infty);{H}_\sigma^2\times (H_{\mathrm{D}}^3\cap H^5))
$$
that satisfies the exponential decay estimate
as in \eqref{201708111550n}.
\end{thm}

Similarly to Theorem \ref{thm-2}, if $\overline{\varphi}$ satisfies the constant RT condition, that is,
\begin{equation}
 \label{constant-RT}
 \partial_{3}\overline{\varphi}>0 \mbox{ is   a constant,}
\end{equation}
and $\Psi''(\overline{\varphi})\geqslant 0$ (i.e., the phase-separation condition fails),
we derive a more concise stability condition, under which the RT instability also does not occur. To this purpose, 
we introduce the following \textit{competitive energy} between the destabilizing and stabilizing factors:
\begin{align}
		&  E_{\mathrm{R}} (u_1,u_2) :=\sigma {l} \partial_{3}\overline{\varphi} \|
		u_1\|_{0}^2-g\varrho_1\|u_2\|^2_{0},\label{E-R}
\end{align}
where $u_1$ and $u_2$ can be scalar or vectorial functions. 
The non-negative part in the competitive energy above is called the \textit{inhibition term}. 
Next, we present the second stability result of the CHRT problem.
\begin{thm}\label{thm:0104nbx}
Suppose that $g>0$ and the assumptions $(\mathbf{A0})$--$(\mathbf{A2})$ are satisfied.
Assume, in addition, $\overline{\varphi}$ satisfies
\eqref{constant-RT} and
\begin{align}
\label{201708101624}
&\omega_{ {\mathrm{R}}}:=\frac{2\sigma{l}\partial_{3}\overline{\varphi}}{g(\rho_1-\rho_2)}>\varpi_{\mathrm{R}},
\end{align}
where
\[\varpi_{\mathrm{R}}:=\max\left\{\sup_{u\in  {H}^1_0,\,  \|\nabla u\|_0\neq0}
 \frac{  \|u\|_{0}^2}{\|\nabla u\|_{0}^2},\ 
  \sup_{u\in H^2\cap H^1_0,\,  \|\Delta u\|_0\neq0} \frac{  \|\nabla u\|_{0}^2}{\|\Delta u\|_{0}^2}, \ 
  \sup_{u\in {H}_\mathrm{D}^3,\,  \|\nabla\Delta u\|_0\neq0} \frac{  \|\Delta u\|_{0}^2}{\|\nabla \Delta u\|_{0}^2}  \right\}.\]
Furthermore, we require
 \begin{equation}
\label{201708101624nnnn}
l^2>-\Psi''\sup_{u\in {H}_\mathrm{D}^3,\,  \|\nabla\Delta u\|_0\neq0} \frac{E_{\mathrm{R}} (\Delta
u ,\nabla u )}{ E_{\mathrm{R}} (\nabla \Delta
u ,\Delta u )}, \text{ provided that } \Psi''<0.
\end{equation}
 Then, there exists a constant $\varepsilon>0$ such that, for any initial data $(\mathbf{v}_0,\phi_0)
\in {H}_\sigma^2\times(H_{\mathrm{D}}^3\cap H^5) $ satisfying the compatibility condition \eqref{compatibility} and $\sqrt{\|\mathbf{v}_0\|_{2}^2+\|\phi_0 \|_{5}^2}<\varepsilon$, 
the CHRT problem  possesses  a unique global strong solution 
$(\mathbf{v},\phi) \in C([0,+\infty),{H}_\sigma^2\times (H_{\mathrm{D}}^3\cap H^5))$, 
which satisfies the  exponential decay estimate as in \eqref{201708111550n}.
\end{thm}
\begin{rem}
The result in Theorem \ref{thm:0104n} subject to the initial-boundary condition $(\mathbf{D})$ can be extended to the initial-boundary condition $(\mathbf{N})$ with minor modifications. 
However, the constant RT condition \eqref{constant-RT} used in Theorem \ref{thm:0104nbx} 
is not compatible with the homogeneous Neumann boundary condition \eqref{Neumann'}  for the common bulk free energy density $\Psi$. 
Further discussions of the admissible equilibrium profile can be found in Appendix \ref{2608071516}.
\end{rem}
\begin{rem}
In analogy to Remark \ref{sharp-poincare}, we see that all three quotients defining $\varpi_{\mathrm{R}}$ attain $1/\lambda_1$, which implies that $\varpi_{\mathrm{R}}=1/\lambda_1$.
Moreover, write $u=\sum_{j}a_j\xi_j$ in an orthonormal Dirichlet eigenbasis with eigenvalues $\lambda_j$, and set $\widetilde{a}=\sigma l \partial_3 \overline{\varphi}$, $\widetilde{b}=g \varrho_1$.
Under \eqref{201708101624}, $\widetilde{a}\lambda_j-\widetilde{b}>0$. The numerator and denominator of the quotient in \eqref{201708101624nnnn} 
are respectively $\sum_j\lambda_j(\widetilde{a}\lambda_j-\widetilde{b})|a_j|^2$ and $\sum_j\lambda_j^2(\widetilde{a}\lambda_j-\widetilde{b})|a_j|^2$. 
Since \(\widetilde{a}\lambda_j-\widetilde{b}>0\) for every \(j\), the quotient is bounded above by \(1/\lambda_1\), with equality attained at the first eigenfunction $\xi_1$. 
For $\Psi''<0$, the additional condition \eqref{201708101624nnnn} reads $l^2\lambda_1> -\Psi''$.
\end{rem}
\begin{rem}
Intuitively, a larger magnitude of the density gradient
$\partial_3 \overline{\rho}=\partial_3\overline{\varphi}(\rho_1-\rho_2)/2$ renders the RT instability more likely to occur in the CHRT problem.
However, from the stability condition \eqref{201708101624}, we observe the inhibition effect of $\partial_{3}\overline{\varphi}$. 
Besides, a larger density contrast between the fluid components facilitates the possible onset of the RT instability. 
\end{rem}
\begin{rem}
We can check that the conclusion of Theorem \ref{thm:0104nbx} also holds by the following modifications: the conditions  \eqref{constant-RT} and
\eqref{201708101624} are replaced by \begin{equation}
\label{20170820}
\partial_{3}\overline{\varphi}<0 \mbox{ is   a constant,}
\end{equation} and ``the CHRT  problem'' is replaced by ``the perturbed system \eqref{purterbation}''.
Under the condition \eqref{20170820}, the equilibrium density profile satisfies $\partial_{3} \overline{\rho} <0$ being a constant. 
This density profile is stable for inhomogeneous fluids within the framework of NS equations in the sense of Hadamard; see \cite[Theorem 1.3]{JFJSAM}. 
Now, it is also exponentially stable under the AGG framework, owing to the dissipation of the phase-field variable (i.e., diffusion).
\end{rem}
%
\subsection{Strategy of proofs}
\label{sec:2.4}
The proofs of the stability results are based on the observation that dissipation terms generally inhibit flow instability, 
particularly when a momentum equation with viscosity is coupled with a parabolic equation that also contains a dissipative term. 
This idea can be implemented using the energy method with suitably constructed energy and dissipation functionals. 
By contrast, the proof of the instability result is based on the so-called bootstrap instability method, which originated in \cite{GuoStrauss,GYSWIC}. 
Subsequently, various versions of this approach have been developed \cite{FSSWVMNA,GYHCSDDC}. 
Here, we adapt the abstract formulation presented in \cite[Lemma 1.1]{GYHCSDDC} to establish Theorem \ref{thm-3}.
In the following, we outline the key steps in the proofs of the main results.

\subsubsection{Stability results for the phase-separation problem}

We first sketch the proof of Theorem \ref{thm-1}. 
The key step is to show the \emph{a priori} exponential decay estimate of small solutions. 
For assumptions $(a)$ and $(b)$, we need to derive energy identities in different ways.
\medskip

\emph{Under assumption $(a)$}, we exploit
the energy method to establish
the following energy identity:
 \begin{align}\frac{1}{2} \frac{\mathrm{d}}{\mathrm{d}t}{\mathcal{E}}_{\mathrm{d}}
  +{  \mathcal{D}}_{\mathrm{d}}  =
I_{\mathrm{d},1}+I_{\mathrm{d},2}, \label{Id1}
   \end{align}
where the definitions  of ${\mathcal{E}}_{\mathrm{d}}$, ${\mathcal{D}}_{\mathrm{d}}$ and $I_{\mathrm{d},i}$ with $i=1,2$, can be found in Lemma \ref{identity'}.
By transforming the temporal derivative estimates  of $(\mathbf{v},\phi )$ into the spatial derivative estimates of $(\mathbf{v},\phi)$,
we can use the regularity theory of elliptic and Stokes equations to show that ${\mathcal{E}}_{\mathrm{d}}$ is equivalent to $\mathcal{E}$ (cf. \eqref{E}).
Noting that ${\mathcal{D}}_{\mathrm{d}}$ includes the  \textit{competitive energy}
$$
E_\mathrm{d}(\phi):=m\sigma\left( l\|\nabla\Delta \phi\|_{0}^2+ l^{-1}\int \Delta(\Psi''(\overline{\varphi})\phi)\Delta \phi\,\mathrm{d}x\right),
$$
thus we need the stability condition \eqref{assume1}
to make sure that the inhibition term $ m\sigma  {l}\| \nabla \Delta \phi\|_{0}^2 $
overcomes the destabilizing terms $m\sigma  {l}^{-1} \int \Delta(\Psi''(\overline{\varphi})\phi)\Delta\phi\,\mathrm{d}x$, 
that is, $\| \nabla \Delta \phi\|_{0}^2 \lesssim  E_{\mathrm{d}}(\phi )$.
With the following observation (cf. \eqref{20170811319561}):
\[
\|(\mathbf{v},\partial_t\mathbf{v})\|_{1}^2+\|(\phi,\partial_t\phi)\|_{3}^2\lesssim \mathcal{D}_\mathrm{d}
 \text{ for sufficiently large }\nu \ \text{or}\ m,
\]
we can verify that the modified dissipation functional ${\mathcal{D}}_{\mathrm{d}} $ is equivalent to $\mathcal{D}$ (cf. \eqref{D}).
Moreover, under the assumption \eqref{aprpiosesnew}, 
the two terms on the right-hand side of  \eqref{Id1} can be controlled by the dissipation term ${\mathcal{D}}_{\mathrm{d}}$.
This enables us to deduce from \eqref{Id1} the exponential decay estimate \eqref{201708111550n}, 
which, together with the existence of local strong solutions of the phase-separation problem (cf. Proposition \ref{pro:0401n20160203}), 
yields the conclusion of Theorem \ref{thm-1} under assumption $(a)$.
\medskip

 \emph{Under assumption $(b)$}, we exploit different testing functions to derive the following energy identity:
\begin{align}
   \frac{1}{2}\frac{\mathrm{d}}{\mathrm{d}t} \widetilde{\mathcal{E}}_{\mathrm{d}}
+\widetilde{\mathcal{D}}_{\mathrm{d}} = \widetilde{I}_{\mathrm{d},1}+\widetilde{I}_{\mathrm{d},2},\label{Id2}
\end{align}
where the definitions  of $\widetilde{\mathcal{E}}_{\mathrm{d}}$, $\widetilde{  \mathcal{D}}_{\mathrm{d}}$ 
and $\widetilde{I}_{\mathrm{d},i}$ with $i=1,2$, can be found in Lemma \ref{identity}.
Since we consider the small perturbation solutions, there exists a constant $C_\mathrm{b}$ such that
 \begin{align}
 \Psi''(\overline{\varphi}+\phi(t))\geq -C_\mathrm{b}\notag
 \end{align}  
 for any sufficiently small $\phi(t)$.
Thus we can verify that the energy functional  $\widetilde{\mathcal{E}}_\mathrm{d}$ is equivalent to $\mathcal{E}$.
The condition \eqref{de-assume-1} together with the following estimate (cf. \eqref{201708251307})
\[
\| \partial_t^i\phi \|_{3}^2\lesssim     \|   \Delta \partial_t^i \phi\|^2_{0}+  \|\nabla ({l}^{-1}  \Psi''(\overline{\varphi})
 \partial_t^i\phi-{l} \Delta \partial_t^i \phi)\|_{0}^2,\quad i=0,1,
\]
implies that (cf. \eqref{201708251845})
\[
\|(\mathbf{v},\partial_t\mathbf{v})\|_{1}^2+\|(\phi,\partial_t\phi)\|_{3}^2\lesssim \widetilde{\mathcal{D}}_\mathrm{d}
 \text{ for sufficiently large }\nu \ \text{or}\ m.
\]
Consequently, we find that $\widetilde{\mathcal{D}}_{\mathrm{d}} $ is equivalent to $ \mathcal{D}$.
In addition, under the assumption \eqref{aprpiosesnew}, 
we can show that $\widetilde{\mathcal{E}}_{\mathrm{d}}$ is equivalent to $\mathcal{E}$ 
and the two terms on the right-hand side of \eqref{Id2} can be controlled by $\widetilde{\mathcal{D}}_\mathrm{d}$.
Thus, from \eqref{Id2} we can again deduce the exponential decay estimate \eqref{201708111550n}, 
which together with the existence of local strong solutions to the phase-separation problem (cf. Proposition \ref{pro:0401n20160203}), 
yields the conclusion of Theorem \ref{thm-1} under assumption $(b)$.

We mention that the lower bound of $\Psi''(r)$ is crucial in the aforementioned argument, 
as the constant $C_\mathrm{b}$ reflects the relative strength of the energy functional $\widetilde{\mathcal{E}}_\mathrm{d}$ 
and the dissipation functional $\widetilde{\mathcal{D}}_\mathrm{d}$.
More precisely, from the definitions of $\widetilde{\mathcal{E}}_\mathrm{d}$ and $\widetilde{\mathcal{D}}_\mathrm{d}$, 
we observe that
$\widetilde{\mathcal{E}}_\mathrm{d}$ contains the terms of the form ${l^2}(\|\nabla\phi\|_{0}^2+\|\nabla\partial_t\phi\|_{0}^2)/{2C_\mathrm{b}}$,
while $\widetilde{\mathcal{D}}_\mathrm{d}$ provides higher-order dissipative terms 
${l^2}{C_\mathrm{b}^{-1}}(\|\nabla\Delta\phi\|_{0}^2+\|\nabla\Delta\partial_t\phi\|_{0}^2)$ (see \eqref{Cb} for details).
Formally, as $C_\mathrm{b}$ tends to infinity, the coefficients of these higher-order norms vanish, 
thereby diminishing the regularizing effect of the dissipation. 

\emph{The proof of Theorem \ref{thm-2} is similar to that of Theorem \ref{thm-1} under assumption $(a)$}. 
Due to \eqref{constant-de} and \eqref{quadratic}, we see that the last term in \eqref{rem-1'} and \eqref{v-lower'} cancel each other by multiplying \eqref{rem-1'} by $\sigma l$. 
This yields a refined energy identity
\begin{align}
    \frac{1}{2}\frac{\mathrm{d}}{\mathrm{d}t}\mathcal{E}^\sharp_\mathrm{d}+\mathcal{D}^\sharp_\mathrm{d}
    =I^\sharp_{\mathrm{d},1}+I^\sharp_{\mathrm{d},2},\notag
\end{align}
where the definitions of ${\mathcal{E}}_{\mathrm{d}}^\sharp$, ${  \mathcal{D}}_{\mathrm{d}}^\sharp$ 
and ${I}^\sharp_{\mathrm{d},i}$ with $i=1,2$, can be found in Lemma \ref{lemma 3.6}.
Since $\mathcal{E}_\mathrm{d}^\sharp$ is similar to $\mathcal{E}_\mathrm{d}$, differing only in some coefficients, 
we find that $\mathcal{E}_\mathrm{d}^\sharp$ is equivalent to $\mathcal{E}$.
Moreover, the bad terms in $\mathcal{D}_\mathrm{d}$ caused by the non-constant $\partial_{3}\overline{\varphi}$ now vanish in $\mathcal{D}_\mathrm{d}^\sharp$.
Therefore, it holds
$$
\|(\mathbf{v},\partial_t\mathbf{v})\|_{1}^2+\|(\phi,\partial_t\phi)\|_{3}^2\lesssim  {\mathcal{D}}^\sharp_\mathrm{d}
$$
without additional assumptions on the size of $\nu$ or $m$.
This enables us to complete the proof of Theorem \ref{thm-2} following the argument for Theorem \ref{thm-1}.

\subsubsection{Instability result for the phase-separation problem}
\medskip

In the spirit of the abstract result \cite[Lemma 1.1]{GYHCSDDC}, the proof of Theorem \ref{thm-3} is divided into five steps:
\begin{itemize}

\item Firstly, we establish a Gronwall-type energy inequality for the phase-separation problem (see Lemma \ref{lem-1}), 
which is crucial for determining the escape time.

\item Secondly, we construct the following profile
 $$
 (\mathbf{v}^{\mathrm{L}},\phi^{\mathrm{L}})=(\boldsymbol{0},e^{\Lambda_{\mathrm{M}} t}\psi)\in {H}_\sigma^2\times( H^3_{\mathrm{D}}\cap H^5),
 $$
 using the variational method, where $\phi^{\mathrm{L}}$ is an unstable solution to the linearized perturbed Cahn--Hilliard equation (see Proposition \ref{exist-unstable-linear}). 
 Then, for any $\delta>0$, the function ${\psi}$ yields a family of functions 
 $(\mathbf{v}^\mathrm{a},\phi^\mathrm{a}):=\delta e^{\Lambda_{\mathrm{M}} t}(\boldsymbol{0},{\psi})\in {H}_\sigma^2\times (H_{\mathrm{D}}^3\cap H^5)$,
 where
 $\phi^\mathrm{a}$ solves the linearized perturbed Cahn--Hilliard equation with the initial datum $ \phi^\mathrm{a}|_{t=0}=\delta {\psi}$.

\item Thirdly, the initial datum $(\boldsymbol{0},\delta\psi)$ of the linear solution may fail to satisfy the compatibility condition \eqref{compatibility}.
We therefore construct a correction $\phi^{\mathrm r}$ explicitly using a regularised distance function, so that
$(\mathbf v_0^\delta,\phi_0^\delta):=(\boldsymbol{0},\delta\psi+\delta^2\phi^{\mathrm r})$
satisfies \eqref{compatibility}; see Proposition \ref{modified-initial-pro}.
Then, we obtain solutions $(\mathbf{v}^\delta,\phi^\delta)$ of the (nonlinear) phase-separation problem 
 originating from the initial data $(\boldsymbol{0},\delta\psi+\delta^2\phi^{\mathrm{r}})$.

\item Fourthly, we derive error estimates between solutions of the linearized perturbed problem and the nonlinear phase-separation problem (see Lemma \ref{lem:0401n}).

\item Finally, based on the aforementioned Gronwall-type inequality and error estimates, 
we are able to prove the existence of a finite escape time $T^\delta$ and thus complete the proof of Theorem  \ref{thm-3}.

\end{itemize}

\subsubsection{Stability results for the CHRT problem}
\medskip

The proof of Theorem \ref{thm:0104n} is based on the stability analysis for the phase-separation problem, 
with the sole distinction arising from the gravity term $-g\varrho_1\phi\mathbf{e}_3$.
In fact, under the assumptions in Theorem \ref{thm:0104n}, 
for a sufficiently large viscosity coefficient $\nu$ or mobility coefficient $m$, 
the dissipation terms $m\sigma l \Delta^2 \phi $
and $\nu \Delta \mathbf{v} $ can jointly inhibit the instability terms, 
that is, the corresponding dissipation integral terms can control all other bad integral terms that may result in flow instability.

The proof of Theorem \ref{thm:0104nbx} follows the argument of Theorem \ref{thm-2}, 
with the aid of some further crucial cancellations enabled by the assumption \eqref{constant-RT}.
Below we give some explanations for the stability conditions \eqref{201708101624} and \eqref{201708101624nnnn}. 
Using the energy method, we can deduce an energy identity (see \eqref{201708070905xc}) that includes the competitive energy functionals
 $ E_{\mathrm{R}} (\nabla\partial_t^i \phi ,\partial_t^i\phi)$, $E_{\mathrm{R}} (\Delta\partial_t^i\phi ,\nabla \partial_t^i\phi)$ 
and $E_{\mathrm{R}} (\nabla \Delta\partial_t^i\phi ,\Delta \partial_t^i\phi)$ for $i=0,1$.
Similar to the role of the stability condition \eqref{201708101624nnnn-de} in the phase-separation problem, 
the stability condition \eqref{201708101624} in the CHRT problem together with $\Psi''\geq 0$ 
ensures that the stabilizing terms in the competitive energies overcome the destabilizing terms. 
This property allows us to deduce the exponential decay estimate from \eqref{201708070905xc}. 
Nevertheless, for the more interesting case of $\Psi''<0$,
the destabilizing terms $  m\sigma   \Psi'' l^{-1}E_{\mathrm{R}} (\Delta\partial_t^i\phi ,\nabla \partial_t^i\phi)$ 
and the inhibition terms $  m\sigma   l E_{\mathrm{R}} (\nabla \Delta\partial_t^i\phi ,\Delta \partial_t^i\phi)$, $i=0,1$, 
also form a competitive relation. 
Hence, we need the condition \eqref{201708101624nnnn} to guarantee that the inhibition terms can overcome the destabilizing terms.

\section{Stability Results for the Phase-Separation Problem}\label{26081414581}

This section is devoted to the proof of the
stability results for the phase-separation problem.
The key step is to show the \emph{a priori}
 exponential decay estimate \eqref{201708111550n} of small solutions. 
 To this end, let $(\mathbf{v},\phi)$ be a strong solution to the phase-separation problem and the \emph{a priori} assumption
\begin{equation}\label{aprpiosesnew}
  \sup_{0\leq t <T}\sqrt{\|\mathbf{v}(t)\|_{2}^2 +\|\phi(t)\|_{5}^2}\leq \delta_1\in (0,1)\mbox{ for some  }T>0,
\end{equation}
where $\delta_1$ is sufficiently small such that
\begin{equation}\label{aprpiosesnew12} 
\phi(x,t)+\overline{\varphi}(x)\in\Big[-\frac{1+\|\overline{\varphi}\|_{L^\infty}}{2},
\frac{1+\|\overline{\varphi}\|_{L^\infty}}{2}\Big]\subset(-1,1) \ 
\mbox{ for all  } (x,t)\in\overline{\Omega}\times [0,T).
\end{equation}
Thus there exists a constant $C_\mathrm{b}>0$ such that
 \begin{align}
 \label{2608131333}
 \Psi''(\overline{\varphi}(x)+\phi(x,t))\geq -C_\mathrm{b} \ 
\mbox{ for all  } (x,t)\in\overline{\Omega}\times [0,T).
 \end{align}
\subsection{Proof of Theorem \ref{thm-1} under the assumption $(a)$}
To begin with, we derive the following energy identity.
\begin{lem}
\label{identity'}
Let $(\mathbf{v},\phi)$ be a strong solution to the phase-separation problem. Then the following energy identity holds
\begin{align}
&\frac{1}{2}\frac{\mathrm{d}}{\mathrm{d}t} {\mathcal{E}}_{\mathrm{d}}
+{\mathcal{D}}_{\mathrm{d}} = I_{\mathrm{d},1}+I_{\mathrm{d},2},\label{d-energy'}
\end{align}
with
\begin{align*}
&\mathcal{E}_\mathrm{d}=\|  \sqrt{\overline{\rho}} \mathbf{v}\|_{0}^2
+\|\nabla\phi\|_{0}^2+\|  \sqrt{\overline{\rho}}\partial_t \mathbf{v}\|_{0}^2
+\|\nabla\partial_t\phi\|_{0}^2+\varrho_1\int \phi|\partial_t\mathbf{v}|^2\,\mathrm{d}x,\\
&\mathcal{D}_\mathrm{d}=\nu\|\nabla\mathbf{v}\|_{0}^2+\nu\|\nabla\partial_t\mathbf{v}\|_{0}^2
+E_\mathrm{d}(\phi)+E_\mathrm{d}(\partial_t\phi)\\
&\qquad\quad-\int\partial_{3}\overline{\varphi}\mathbf{v}_3\Delta\phi\,\mathrm{d}x
+  \sigma {l} \int \Big(\partial_{3}^2\overline{\varphi}(\partial_{1} \phi \mathbf{v}_1+ \partial_{2}\phi \mathbf{v}_2)-
\partial_{3}\overline{\varphi}\nabla \phi\cdot \nabla \mathbf{v}_3\Big)  \, \mathrm{d}x\\
&\qquad\quad-\int\partial_{3}\overline{\varphi}\partial_t\mathbf{v}_3\Delta\partial_t\phi\,\mathrm{d}x
+ \sigma {l} \int  \Big(\partial_{3}^2\overline{\varphi}(\partial_{1}\partial_t \phi \partial_t\mathbf{v}_1+ \partial_{2}\partial_t\phi \partial_t\mathbf{v}_2)-
\partial_{3}\overline{\varphi}\nabla \partial_t\phi\cdot \nabla\partial_t \mathbf{v}_3\Big)  \, \mathrm{d}x,\\
&I_{\mathrm{d},1}=\int  (\mathbf{M}_{\mathrm{d}} - \varrho_1 \phi  \partial_t\mathbf{v})\cdot  \mathbf{v}\,\mathrm{d}x
-\int N_\mathrm{d}\Delta\phi\,\mathrm{d}x,\\
&I_{\mathrm{d},2}=\int  (\partial_t\mathbf{M}_{\mathrm{d}}- \frac{\varrho_1}{2}\partial_t \phi  \partial_t\mathbf{v})\cdot\partial_t\mathbf{v}\,\mathrm{d}x
-\int \partial_t N_\mathrm{d}\Delta\partial_t\phi\,\mathrm{d}x.
\end{align*}
Here, the quantities $E_\mathrm{d}$, $\mathbf{M}_{\mathrm{d}}$ and $N_\mathrm{d}$ are defined as in \eqref{E-d}, \eqref{M-d} and \eqref{N-d}, resp.
\end{lem}
\begin{proof}
First, we rewrite the system \eqref{purterbation} with $g=0$ as the following linearized equations with nonhomogeneous terms
\begin{equation}\label{g=0-linear}
\begin{cases}
 \overline{\rho}\partial_t \mathbf{v}- \nu \Delta \mathbf{v}+\nabla  q
-\varrho_1 m\xi  \partial_{3} \mathbf{v}+\sigma {l}
  (\partial_{3}\overline{\varphi} \Delta \phi \mathbf{e}_3 +\partial_{3}^{2}\overline{\varphi}\nabla \phi)
  =\mathbf{M}_{\mathrm{d}}- \varrho_1 \phi \partial_t \mathbf{v},\\
\mathrm{div}\,\mathbf{v} =0, \\
\partial_t \phi+ \partial_{3}\overline{\varphi}\mathbf{v}_3-m\sigma \Delta( {l}^{-1}  \Psi''(\overline{\varphi})
 \phi-{l} \Delta  \phi)=N_{\mathrm{d}},
\end{cases}
 \end{equation}
with $\overline{\mu}
=\sigma(l^{-1}\Psi'(\overline{\varphi})- l\partial_3^2\overline{\varphi})$, 
$\xi:= \partial_{3}\overline{\mu}$, 
and
\begin{align}
\mathbf{M}_{\mathrm{d}}
&:= \varrho_1 m\sigma \Big(\nabla\Big({l}^{-1}\int^{1}_{0} \Psi''( \tau \phi+\overline{\varphi})\phi\,\mathrm{d}\tau
-  {l} \Delta \phi\Big)\cdot \nabla\Big) \mathbf{v}
 \notag\\
 &\quad \ - (\varrho_1\phi+ \overline{\rho} )( \mathbf{v}\cdot\nabla)\mathbf{v}
  -\sigma {l} \Delta \phi\nabla \phi,
  \label{M-d}\\
N_{\mathrm{d}}&:=m\sigma  {l}^{-1} \Delta \int_0^1\int_0^1  \Psi'''(s\tau \phi+\overline{\varphi})\tau \phi^2\,\mathrm{d}s\,\mathrm{d}\tau
-\mathbf{v}\cdot \nabla \phi.
\label{N-d}
\end{align}
We note that
$\xi$ is a constant by \eqref{20171805}.
Taking the inner product of \eqref{g=0-linear}$_3$ and $-\Delta\phi$ in $L^2$, we get
    \begin{align}
        &\frac{1}{2}\frac{\mathrm{d}}{\mathrm{d}t}\|\nabla\phi\|_{0}^2
        +m\sigma l\|\nabla\Delta\phi\|_{0}^2
        +m\sigma l^{-1}\int \Delta(\Psi''(\overline{\varphi})\phi)\Delta\phi\,\mathrm{d}x
        =\int\partial_{3}\overline{\varphi}\mathbf{v}_3\Delta\phi\,\mathrm{d}x
        -\int N_\mathrm{d}\Delta\phi\,\mathrm{d}x.\label{rem-1}
\end{align}
Next, taking the inner product of \eqref{g=0-linear}$_1$ and $\mathbf{v}$ in $L^2$ implies that
\begin{align}
&\frac{1}{2}\frac{\mathrm{d}}{\mathrm{d}t}\|  \sqrt{\overline{\rho}} \mathbf{v}\|_{0}^2
+ \nu \|\nabla \mathbf{v}\|_{0}^2 
+ \sigma {l} \int  \Big(\partial_{3}^{2}\overline{\varphi}(\partial_{1} \phi \mathbf{v}_1+ \partial_{2}\phi \mathbf{v}_2)-
\partial_{3}\overline{\varphi}\nabla \phi\cdot \nabla \mathbf{v}_3\Big)  \, \mathrm{d}x\nonumber
\\
&\quad
 = \frac{\varrho_1 m \xi }{2}\int\partial_{3}|\mathbf{v}|^2\,\mathrm{d}x
 + \int (\mathbf{M}_{\mathrm{d}} - \varrho_1 \phi  \partial_t\mathbf{v})\cdot  \mathbf{v}\,\mathrm{d}x
 = \int  (\mathbf{M}_{\mathrm{d}} - \varrho_1 \phi  \partial_t\mathbf{v})\cdot  \mathbf{v}\,\mathrm{d}x.\label{v-lower}
 \end{align}
Combining \eqref{rem-1} and \eqref{v-lower}, we obtain
\begin{align}
    &\frac{1}{2}\frac{\mathrm{d}}{\mathrm{d}t}\Big(\|  \sqrt{\overline{\rho}} \mathbf{v}\|_{0}^2
+\|\nabla\phi\|_{0}^2\Big)+ \nu \|\nabla \mathbf{v}\|_{0}^2
+m\sigma l\|\nabla\Delta\phi\|_{0}^2+m\sigma l^{-1}\int \Delta(\Psi''(\overline{\varphi})\phi)\Delta\phi\,\mathrm{d}x\notag\\
&\qquad-\int\partial_{3}\overline{\varphi}\mathbf{v}_3\Delta\phi\,\mathrm{d}x
+ \sigma {l} \int  \Big(\partial_{3}^{2}\overline{\varphi}(\partial_{1} \phi \mathbf{v}_1+ \partial_{2}\phi \mathbf{v}_2)-
\partial_{3}\overline{\varphi}\nabla \phi\cdot \nabla \mathbf{v}_3\Big)  \, \mathrm{d}x\nonumber
\\
 &\quad= \int  (\mathbf{M}_{\mathrm{d}} - \varrho_1 \phi  \partial_t\mathbf{v})\cdot  \mathbf{v}\,\mathrm{d}x
 -\int N_\mathrm{d}\Delta\phi\,\mathrm{d}x.\label{re-energy}
\end{align}
Next, differentiating the system \eqref{g=0-linear} with respect to $t$, we have
\begin{equation}\label{g=0-linear-t}
\begin{cases}
 \partial_{t}[(\overline{\rho}+\varrho_1\phi)\partial_t \mathbf{v}]- \nu \Delta\partial_t \mathbf{v}+\nabla  \partial_tq
-\varrho_1 m\xi  \partial_{3}\partial_t \mathbf{v}+\sigma {l}
  (\partial_{3}\overline{\varphi} \Delta\partial_t \phi \mathbf{e}_3 +\partial_{3}^{2}\overline{\varphi}\nabla \partial_t\phi) =\partial_t\mathbf{M}_{\mathrm{d}},\\
\mathrm{div}\,\partial_t\mathbf{v} =0, \\
\partial_{t}^{2} \phi+ \partial_{3}\overline{\varphi}\partial_t\mathbf{v}_3-m\sigma \Delta({l}^{-1}  \Psi''(\overline{\varphi})
 \partial_t\phi-{l} \Delta  \partial_t\phi)=\partial_t N_{\mathrm{d}}.
\end{cases}
 \end{equation}
 Similarly to \eqref{rem-1} and \eqref{v-lower}, formally testing \eqref{g=0-linear-t}$_1$ and \eqref{g=0-linear-t}$_3$ by $\partial_t\mathbf{v}$ and $-\Delta\partial_t\phi$, respectively, and adding the resulting identities, we obtain
 \begin{align}
&\frac{1}{2}\frac{\mathrm{d}}{\mathrm{d}t}\Big(\|  \sqrt{\overline{\rho}} \partial_t\mathbf{v}\|_{0}^2
+\|\nabla\partial_t\phi\|_{0}^2 +\varrho_1\int \phi|\partial_t\mathbf{v}|^2 \,\mathrm{d}x\Big)
+ \nu \|\nabla \partial_t\mathbf{v}\|_{0}^2 
+m\sigma l\|\nabla\Delta\partial_t\phi\|_{0}^2 \notag\\
&\qquad+m\sigma l^{-1}\int \Delta(\Psi''(\overline{\varphi})\partial_t\phi)\Delta\partial_t\phi\,\mathrm{d}x
-\int\partial_{3}\overline{\varphi}\partial_t\mathbf{v}_3\Delta\partial_t\phi\,\mathrm{d}x\notag\\
&\qquad+ \sigma {l} \int  \Big(\partial_{3}^{2}\overline{\varphi}(\partial_{1}\partial_t \phi \partial_t\mathbf{v}_1
+ \partial_{2}\partial_t\phi \partial_t\mathbf{v}_2)-
\partial_{3}\overline{\varphi}\nabla \partial_t\phi\cdot \nabla\partial_t \mathbf{v}_3\Big)  \, \mathrm{d}x\nonumber
\\
 &\quad= \int  \Big(\partial_t\mathbf{M}_{\mathrm{d}} - \frac{\varrho_1}{2}\partial_t \phi  \partial_t\mathbf{v}
 \Big)\cdot \partial_t \mathbf{v}\,\mathrm{d}x
 -\int \partial_tN_\mathrm{d}\Delta\partial_t\phi\,\mathrm{d}x.\label{re-energy-t}
 \end{align}
These formal calculations can be justified by interpreting the first and third equations of \eqref{g=0-linear-t} in $(H_\sigma^1)'$ and $(H_0^1)'$, respectively. Indeed, by \cite[Remark 6]{ChoKim}, we have
\[\langle\partial_t((\overline{\rho}+\varrho_1\phi)\partial_t\mathbf{v}),\partial_t\mathbf{v}\rangle_{(H_\sigma^1)',H_\sigma^1}=\frac{1}{2}\frac{\mathrm{d}}{\mathrm{d}t}\int(\overline{\rho}+\varrho_1\phi)|\partial_t\mathbf{v}|^2\,\mathrm{d}x+\frac{\varrho_1}{2}\int \partial_t\phi|\partial_t\mathbf{v}|^2\,\mathrm{d}x.\]
Moreover, the regularity established in Proposition \ref{pro:0401n20160203}, together with \eqref{g=0-linear-t}$_3$, yields $$\partial_t^2 \phi\in L^2(0,T;(H_0^1)').$$ Applying the chain rule in \cite[Lemma 4.1]{RS2004}, we obtain
\[\langle\partial_t^2\phi,-\Delta\partial_t\phi\rangle_{(H_0^1)',H_0^1}=\frac{1}{2}\frac{\mathrm{d}}{\mathrm{d}t}\|\nabla\partial_t\phi\|_0^2.\]
Thus, \eqref{re-energy-t} holds rigorously.
Combining \eqref{re-energy} and \eqref{re-energy-t}, we arrive at the conclusion  \eqref{d-energy'}.
\end{proof}

Our next aim is to show that ${\mathcal{E}}_{\mathrm{d}}$ and ${\mathcal{D}}_{\mathrm{d}}$ are equivalent to $\mathcal{E}$ and ${\mathcal{D}}$, resp. 
For this purpose, we shall use the following elliptic and inhibition estimates.
\begin{lem}
	\label{2017081921311043}
	\begin{enumerate}[ (1)]
		\item It holds that
		\begin{equation}
			\label{201708131702n}
			\|u\|_{3}\lesssim\|\nabla \Delta  u\|_{0}
	\mbox{ for any }u\in  H_{\mathrm{D}}^3.
		\end{equation}
		\item   It holds that
		\begin{align}
			\label{201708131719n}
			\|u\|_{i+4 }\lesssim \| \Delta^2 u\|_{i}
            \mbox{ for any }u\in H_{\mathrm{D}}^3\cap H^{i+4},\mbox{ where }1\leq i\leq 3.
		\end{align}
		\item Under the stability condition \eqref{assume1}, we have the inhibition estimate
			\begin{equation}\label{201708130943n}
			\|  u\|_{3}^2\lesssim   {E_{\mathrm{d}}(u)}\mbox{ for any }u\in  H_{\mathrm{D}}^3.
		\end{equation}
	\end{enumerate}
\end{lem}
\begin{proof}
	(1) Since $u\in H_{\mathrm{D}}^3$,
    by the classical theory of elliptic equations with the Dirichlet boundary condition (see \cite[Lemma 5.8]{NASII04}), we have
	\begin{equation}
		\label{201708241039}
		\|u\|_{i+2}\lesssim \|  \Delta  u \|_{i}\ \mbox{ for }i=0,1.
	\end{equation}
	Noting that $\Delta u|_\Gamma=0$,
	by Poincar\'e's inequality \eqref{250508-Poincare}, we can conclude
	\eqref{201708131702n}.

	(2)
Noting that
$$
\begin{cases}
	\Delta  (\Delta  u) = \Delta^2 u  & \hbox{in }\Omega, \\
	\Delta   u=0 & \hbox{on }\Gamma,
\end{cases}
$$
similarly, we have
\begin{equation*}
	\|\Delta u\|_{i+2}\lesssim
	\|  \Delta^2 u \|_{i}\ \text{ for }1\leq i\leq 3.
\end{equation*}
In addition, using the elliptic estimate for $u$, we obtain
\begin{equation*}
	\|u\|_{i+4}\lesssim
	\|  \Delta  u \|_{i+2}\ \text{ for }1\leq i\leq 3.
\end{equation*}
Combining the above estimates yields  \eqref{201708131719n}.
(3) If $\|\nabla\Delta u\|_0=0$, then \eqref{201708131702n}
implies that $u=0$, and thus \eqref{201708130943n} holds trivially.
Otherwise, by the stability condition  \eqref{assume1}, we find
		\begin{align}
			E_{\mathrm{d}}(u)&=m\sigma l^{-1}\|\nabla\Delta u\|_{0}^2\Big(l^2
            + {\|\nabla\Delta u\|_{0}^{-2}}{\int\Delta(\Psi''(\overline{\varphi})u)\Delta u\,
            \mathrm{d}x}\Big) \geq m\sigma l^{-1}(l^2-\omega_\mathrm{d})\|\nabla\Delta u\|_{0}^2,\notag
		\end{align}
		which, together with the estimate \eqref{201708131702n}, yields the conclusion \eqref{201708130943n}.
\end{proof}

Using Lemma \ref{2017081921311043} and the
regularity theory of elliptic equations and Stokes equations, 
we can show the equivalence between ${\mathcal{E}}_{\mathrm{d}}$, resp.  ${\mathcal{D}}_{\mathrm{d}}$ and $\mathcal{E}$, resp. ${\mathcal{D}}$.
\begin{lem}
\label{lem:2070811nn}
Under the assumption \eqref{aprpiosesnew} with a sufficiently small $\delta_1$,
the following two assertions hold:
\begin{enumerate}[\quad (1)]
  \item  
  $\mathcal{E}$ is equivalent to $\mathcal{E}_{\mathrm{d}}$.
  If the stability condition  \eqref{assume1} holds and either $\nu$ or $m$ is sufficiently large, then   ${\mathcal{D}}_{\mathrm{d}}$ is  equivalent to $\mathcal{D}$.
  \item ${\mathcal{E}}_{\mathrm{d}}$ is equivalent to $ \|\mathbf{v}\|_{2}^2+\|\phi\|_{5}^2$. In particular,
      we have
\begin{equation}
\label{20170811}
{\mathcal{E}}_{\mathrm{d}}|_{t=0}\lesssim ( {\|\mathbf{v}_0\|_{2}^2}+\|\phi_0\|_{5}^2).
\end{equation}
\end{enumerate}
\end{lem}
\begin{proof}
(1) Using $(\mathbf{A2})$ and H\"{o}lder's inequality, it is straightforward to check that ${\mathcal{D}}_{\mathrm{d}}\lesssim\mathcal{D}$.
In addition, using the Sobolev embedding theorem and \eqref{aprpiosesnew}, we find
$
{\mathcal{E}}_{\mathrm{d}}\lesssim\mathcal{E}+\|\phi\|_{2}\|\partial_t\mathbf{v}\|_{0}^2\lesssim \mathcal{E}.
$
Hence, to get the first assertion of Lemma \ref{lem:2070811nn}, it suffices to show that 
\begin{equation}
 \label{201708122142xx}
  \mathcal{E}\lesssim
 {\mathcal{E}}_{\mathrm{d}}
 \end{equation}
 and for sufficiently large $\nu$ or $m$,
 \begin{equation}
 \label{201708122142xxD}
  \mathcal{D}\lesssim {\mathcal{D}}_{\mathrm{d}}\quad\text{under the stability condition } \eqref{assume1}.
 \end{equation}
Noting that the strong solution $(\mathbf{v},\phi)$ satisfies
$$
\mathbf{v}=\partial_t\mathbf{v}=\boldsymbol{0},\ \phi= \Delta \phi=\partial_t\phi=\Delta\partial_t\phi=0 \mbox{ on }\Gamma, \quad\text{for a.a. } t\in(0,T),
$$
using the assumption \eqref{aprpiosesnew} with a sufficiently small $\delta_1$ and Poincar\'e's inequality \eqref{250508-Poincare}, we can deduce that
\begin{equation}
\label{2017081131956}
 \| (\mathbf{v},\partial_t\mathbf{v})\|_{0}^2 + \| (\phi,\partial_t\phi )\|_{1}^2
\lesssim {\mathcal{E}}_{\mathrm{d}}.
\end{equation}
Under the stability condition \eqref{assume1}, the inhibition estimates \eqref{201708130943n} holds. 
Then, for sufficiently large viscosity $\nu$ or mobility $m$, we can apply H\"older's inequality and Young's inequality to obtain
\begin{equation}
\label{20170811319561}
\|  (\mathbf{v} ,\partial_t \mathbf{v}) \|_{1}^2+\|(\phi,\partial_t\phi)\|_{3}^2\lesssim {\mathcal{D}}_{\mathrm{d}} .
\end{equation}

To complete the proof of \eqref{201708122142xx}, 
we transform the temporal derivative estimates of $(\mathbf{v},\phi)$ into the spatial derivative estimates of $(\mathbf{v},\phi)$.
Recalling $\Psi\in C^8(-1,1)$, we have,  for $1\leq i\leq 8$,
\begin{equation}
\label{201708170716}
\|\Psi^{(i)}(\overline{\varphi}+\phi(t))\|_{L^\infty}\lesssim 1\text{ for all }t\in[0,T),
\end{equation}
and \begin{equation}
\label{201708170716xdfs}
\|\Psi^{(i)}(\overline{\varphi}+\tau\phi(t))\|_{L^\infty}\lesssim1  \mbox{ for all }\tau \in[0,1] \text{ and }t\in [0,T)  .
\end{equation}
Exploiting the product estimate \eqref{product}, \eqref{aprpiosesnew}, \eqref{201708170716}, \eqref{201708170716xdfs} and H\"older's inequality, we obtain
 \begin{align}
\|N_{\mathrm{d}}\|_{i}&\leq
\Big\|m\sigma  {l}^{-1} \Delta \int_0^1\int_0^1  \Psi'''(s\tau \phi+\overline{\varphi})\tau \phi^2\,\mathrm{d}s\,\mathrm{d}\tau
\Big\|_{i} + \|\mathbf{v}\cdot \nabla \phi\|_{i}
\nonumber
\\
& \lesssim
\begin{cases}
 (\| \mathbf{v}\|_{1}+\|\phi\|_{2})\|\phi\|_{3} & \hbox{for }i=1, \\
 (\|\mathbf{v}\|_2+\|\phi\|_4)\|\phi\|_3&\hbox{for }i=2,\\
(\| \mathbf{v}\|_{3}+\| \phi\|_{5})\|\phi\|_{3}+\|\mathbf{v}\|_2\|\phi\|_4
& \hbox{for }i=3.
\end{cases}
\label{201708132007}
\end{align}
Next, using  \eqref{201708131719n} and \eqref{201708132007}, we deduce from \eqref{g=0-linear}$_3$ that
\begin{align}
\|  \phi\|_{i+4}\lesssim \|\Delta^2\phi \|_{i}&
\lesssim  \|(\Delta (\Psi''(\overline{\varphi})\phi),\partial_t\phi,N_{\mathrm{d}},\partial_{3}\overline{\varphi}\mathbf{v}_3)\|_{i}\notag\\
&\lesssim
\begin{cases}
\|\phi\|_{3}+\|\partial_t\phi\|_{1}+\| \mathbf{v}\|_{1}  & \hbox{for }i=1, \\
\|\phi\|_{5}+\|\partial_t\phi\|_{3}+\| \mathbf{v}\|_{3} & \hbox{for }i=3.
\end{cases}                                    \notag
\end{align}
Then using the interpolation inequality \eqref{201708190532},
we further get
\begin{equation}
\label{201708230922}
\|  \phi\|_{i+4}  \lesssim
 \begin{cases}
\|\phi\|_{0}+ \|\partial_t\phi\|_{1}+\| \mathbf{v}\|_{1}  & \hbox{for }i=1, \\
\|\phi\|_{0}+\|\partial_t\phi\|_{3}+ \|\mathbf{v}\|_{3} & \hbox{for }i=3.
\end{cases}
\end{equation}

Consider the Stokes problem
$$
\begin{cases}
- \nu \Delta \mathbf{v} +\nabla  q 
= \mathbf{M}_{\mathrm{d}}+\mathbf{L}_\mathrm{d}-(\overline{\rho}+\varrho_1 \phi  )\partial_t\mathbf{v} &\text{in}\ \Omega_T,\\
\mathrm{div}\,\mathbf{v} =0&\text{in}\ \Omega_T\\
\mathbf{v}=\boldsymbol{0}&\text{on}\ \Gamma_T,
\end{cases}
$$
where
\[
\mathbf{L}_\mathrm{d}:=\varrho_1 m\xi  \partial_{3} \mathbf{v}-\sigma {l}
  (\partial_{3}\overline{\varphi} \Delta \phi \mathbf{e}_3 +\partial_{3}^{2}\overline{\varphi}\nabla \phi).
\]
Using the Stokes estimate (cf. \cite[Theorem IV 6.1]{galdi2011introduction}) and the product estimate \eqref{product}, 
we can deduce that
\begin{align}
    \|\mathbf{v}\|_{i+2}+\|\nabla q\|_{i}\lesssim \|\partial_t\mathbf{v}\|_{i}(1+\|\phi\|_{2})
+ \|\mathbf{M}_{\mathrm{d}}\|_{i}+\|\mathbf{L}_{\mathrm{d}}\|_{i} \quad \text{ for}\ i=0, 1.\label{260902-1}
\end{align}
In addition, using \eqref{product}, \eqref{aprpiosesnew}, \eqref{201708170716} and \eqref{201708170716xdfs}, we have, for $i=0,1$,
 \begin{align}
 \|\mathbf{M}_{\mathrm{d}}\|_{i}
 &\lesssim (1+\|\phi\|_{2}) \|\mathbf{v}\cdot \nabla \mathbf{v}\|_{i}
+\Big\|\Big(\nabla\Big({l}^{-1}\int^{1}_{0} \Psi''( \tau \phi+\overline{\varphi})\phi\,\mathrm{d}\tau 
-  {l} \Delta \phi\Big)\cdot \nabla\Big)  \mathbf{v}\Big\|_{i}
 +\|\Delta \phi\nabla \phi\|_{i} \nonumber \\
 &\lesssim\Big(
\|\mathbf{v}\|_{2}\|\mathbf{v}\|_{i+1}
+\|  \phi\|_{3}(\| \mathbf{v}\|_{i+1} + \| \phi\|_{i+2})+\| \mathbf{v}\|_{i+1}\| \phi\|_{5}\Big),\label{Md}\\
\|\mathbf{L}_\mathrm{d}\|_{i}&\lesssim\|\mathbf{v}\|_{i+1}+\|\phi\|_{i+2}.\label{Ld}
\end{align}
Applying the interpolation inequality \eqref{201708190532} again, 
we infer from \eqref{aprpiosesnew}, \eqref{260902-1}, \eqref{Md} and \eqref{Ld} that
\begin{equation}
\label{201708131955}
\|\mathbf{v}\|_{i+2}+\|\nabla q\|_{i}\lesssim\|\partial_t\mathbf{v}\|_{i}
+\|\mathbf{v}\|_{0}+ \| \phi\|_{i+2}.
\end{equation}
Hence, it follows from \eqref{201708190532}, \eqref{201708230922} and \eqref{201708131955} that
\begin{align*}
& \|\mathbf{v}\|_{2}+\|  \phi\|_{5}  +\|\nabla q\|_{0}\lesssim
\|\partial_t\mathbf{v}\|_{0}+\|\mathbf{v}\|_{0}+\|\phi\|_{0}+\|\partial_t\phi\|_{1},\\
& \|\mathbf{v}\|_{3}+\|  \phi\|_{7}  +\|\nabla q\|_{1}\lesssim
\| \partial_t\mathbf{v} \|_{1}+\|\mathbf{v}\|_{0}+\|\phi\|_{0} +\|\partial_t\phi\|_{3}.
\end{align*}
Combining the above two estimates with \eqref{2017081131956}
and \eqref{20170811319561}, we arrive at \eqref{201708122142xx} and \eqref{201708122142xxD}.

(2) Noting that
$$
{\mathcal{E}} \lesssim {\mathcal{E}}_{\mathrm{d}}\lesssim  
  \| \mathbf{v}\|^2_{0}+\|\partial_t\mathbf{v}\|_{0}^2 +\|\phi\|_{1}^2  + \| \partial_t\phi\|_{1}^2 ,
$$
to prove the desired assertion, it suffices to verify that
 \begin{equation}
\label{201708132040}
\|\partial_t \mathbf{v}\|_{0} +   \| \partial_t \phi\|_{1} \lesssim  {\|\mathbf{v}\|_{2} +\|\phi\|_{5} }.
\end{equation}
Taking the inner product of \eqref{g=0-linear}$_1$ and $\partial_t\mathbf{v}$ in $L^2$, we obtain
\begin{align*}
 \|\sqrt{\overline{\rho}}\partial_t\mathbf{v}\|_{0}^2
 &=\nu \int \Delta \mathbf{v} \cdot \partial_t\mathbf{v}\,\mathrm{d}x
 +\int(\mathbf{M}_{\mathrm{d}}+\mathbf{L}_{\mathrm{d}}-\varrho_1 \phi  \partial_t\mathbf{v})\cdot\partial_t \mathbf{v}\,\mathrm{d}x\\
&\lesssim \| \sqrt{\overline{\rho}}\partial_t\mathbf{v}\|_{0}(\|\mathbf{v}\|_{2}
+ \|\mathbf{M}_{\mathrm{d}}\|_{0}+\|\mathbf{L}_{\mathrm{d}}\|_{0})
+\|\phi\|_2\|\sqrt{\overline{\rho}}\partial_t\mathbf{v}\|_0^2,
\end{align*}
which, together with \eqref{aprpiosesnew}, \eqref{Md} and \eqref{Ld}, yields the following estimate
$$  \|\partial_t\mathbf{v}\|_{0} \lesssim\|\mathbf{v}\|_{2}+\|\phi\|_{5}.$$
Similarly, from \eqref{g=0-linear}$_3$ and \eqref{201708132007} we can deduce that
\begin{equation}
\notag
\|\partial_t\phi\|_{1}\lesssim \|\mathbf{v}\|_{2}+\|\phi\|_{5}.
\end{equation}
Consequently, we see that \eqref{201708132040} holds. 
The conclusion \eqref{20170811} is a natural consequence for the strong solution. 
The proof is complete.
\end{proof}

The subsequent lemma shows that the two integrals on the right-hand side of the energy identity \eqref{d-energy'} can be controlled by $ \mathcal{D}$.
\begin{lem}\label{201708111545}
Under the assumption \eqref{aprpiosesnew} with  $\delta_1\in (0,1)$, we can deduce that
$$
|I_{\mathrm{d},1}|+|I_{\mathrm{d},2}|\lesssim \sqrt{\mathcal{E}}{\mathcal{D}}.$$
\end{lem}
\begin{proof}
It follows from \eqref{product}, \eqref{aprpiosesnew}, \eqref{201708170716}, \eqref{201708170716xdfs} and the Sobolev embedding theorem that
$$
\begin{aligned}
&\| N_{\mathrm{d}}\|_{0}\lesssim (\|\mathbf{v}\|_{1}
 +\| \phi\|_{2} )\|  \phi\|_{2},\\
 &
\|\partial_t \mathbf{M}_{\mathrm{d}}\|_{0} \lesssim
\|\partial_t\mathbf{v}\|_{1}(\|\mathbf{v}\|_{2}+\|\phi\|_{5}) +\|\partial_t\phi\|_{3}( \|\mathbf{v}\|_{3}+\|\phi\|_{2})
 ,\\
&\|\partial_t N_{\mathrm{d}}\|_{0}\lesssim  \|\phi\|_{2} \|\partial_t\phi\|_{2}
+\|\partial_t\mathbf{v}\|_{0}\|\phi\|_{3}+\|\mathbf{v}\|_{2}\|\partial_t\phi\|_{1}.
\end{aligned}
$$
Using  \eqref{product}, \eqref{aprpiosesnew}, \eqref{Md},  the above three estimates and the Sobolev embedding theorem, we obtain
   $$
   \begin{aligned}
  &| I_{\mathrm{d},1}|\lesssim
\|\mathbf{v}\|_{0}(\| \mathbf{M}_{\mathrm{d}}\|_{0}+\|\phi\|_{2}\|\partial_t\mathbf{v}\|_{0}) 
+ \|  N_{\mathrm{d}} \|_{0}\| \phi\|_{2}
    \lesssim \sqrt{\mathcal{E}}{\mathcal{D}},\\
   &  |I_{\mathrm{d},2} | \lesssim \|\partial_t\phi\|_{1} \|\partial_t\mathbf{v}\|^2_{1}
   + \|\partial_t\mathbf{v}\|_{0} \| \partial_t \mathbf{M}_{\mathrm{d}}\|_{0}
   + \|\partial_t\phi\|_{2}
\|\partial_t N_{\mathrm{d}}\|_{0}\lesssim \sqrt{\mathcal{E}}{\mathcal{D}}.
  \end{aligned}$$
Combining the above two estimates, we immediately get the desired conclusion.
\end{proof}

\emph{We are now in a position to prove Theorem \ref{thm-1} under the assumption (a)}. 
By Lemmas \ref{identity'}, \ref{201708111545}, and the equivalence between ${\mathcal{D}}_{\mathrm{d}} $ and $\mathcal{D}$,  
we can deduce from \eqref{aprpiosesnew} and \eqref{d-energy'} that,  for $\delta_1$ being sufficiently small,
 \begin{align}
  \frac{\mathrm{d}}{\mathrm{d}t} {\mathcal{E}}_{\mathrm{d}}
  + C{  \mathcal{D}} \leq 0\mbox{ on }(0,T). \label{201708070905}
  \end{align}
Since $\mathcal{E}_{\mathrm{d}}\lesssim \mathcal{E}$ and ${\mathcal{E}} \leq {\mathcal{D}} $, we have
    ${\mathcal{E}}_{\mathrm{d}} \lesssim {  \mathcal{D}} $.
Thus, it follows from \eqref{201708070905} that
$$
 e^{c_1 t} {\mathcal{E} }_{\mathrm{d}}(t)+C\int_0^t e^{c_1\tau} \mathcal{D} (\tau)\,\mathrm{d}\tau\leq {\mathcal{E}}_{\mathrm{d}}|_{t=0}.
$$
By \eqref{20170811} and the equivalence between  ${\mathcal{E}}_{\mathrm{d}}$ and $\mathcal{E}$,
  we deduce the exponential decay estimate
     \eqref{201708111550n}. To complete the proof, we combine this estimate with Proposition \ref{pro:0401n20160203}.
    Let $\varepsilon$ in Theorem \ref{thm-1} further satisfy $$\varepsilon<\delta_1/(2\max\{1,\sqrt{c_2}\}).$$
    This, together with \eqref{aprpiosesnew12}, implies that $\|\phi_0+\overline{\varphi}\|_{L^\infty}\leq (1+\|\overline{\varphi}\|_{L^\infty})/2<1$. 
     Then, by Proposition \ref{pro:0401n20160203}, there exists a unique solution $(\mathbf{v},\phi)\in C([0,T^{\max});H_\sigma^2\times H^5)$,
     where $T^{\max}$ denotes its maximal existence time.
     Define $$T^\sharp:=\sup\{t\in(0,T^{\max}): \sqrt{\|\mathbf{v}(s)\|_2^2+\|\phi(s)\|_5^2}\leq \delta_1 \text{ for any }s\in[0,t]\}.$$
     Since $\sqrt{\|\mathbf v_0\|_2^2+\|\phi_0\|_5^2}
<\varepsilon<\delta_1/2$,
the continuity of $(\mathbf v,\phi)$ implies that
$T^\sharp>0$.
     Now we prove that $T^\sharp=+\infty$ by contradiction. Suppose that $T^\sharp<+\infty$.
     By the definition of $T^\sharp$
and \eqref{aprpiosesnew12}, we have
\[
\|\mathbf v(t)\|_2+\|\phi(t)\|_5\leq\sqrt{2}\delta_1,
\qquad
\|\overline{\varphi}+\phi(t)\|_{L^\infty}
\leq1-\eta,\quad\text{for all }t\in[0,T^\sharp),
\]
where $\eta:=(1-\|\overline{\varphi}\|_{L^\infty})/2>0$.
The local existence argument in \cite{Lv2025} therefore
provides a uniform existence time when the solution is
restarted at any $t_0\in[0,T^\sharp)$.
If $T^\sharp=T^{\max}<+\infty$, choosing $t_0$ sufficiently
close to $T^{\max}$ would extend the solution beyond
$T^{\max}$, contradicting the maximality of $T^{\max}$.
    Therefore, $T^\sharp<T^{\max}$ whenever $T^\sharp<+\infty$. Then, by the continuity of $(\mathbf{v},\phi),$ we have $\sqrt{\|\mathbf{v}(T^\sharp)\|_2^2+\|\phi(T^\sharp)\|_5^2}=\delta_1$.
    In view of the definition of $T^\sharp$, the solution $(\mathbf{v},\phi)$ also satisfies \eqref{201708111550n} on $[0,T^\sharp)$, which, together with $\varepsilon<\delta_1/(2\sqrt{c_2})$ and the continuity of $(\mathbf{v},\phi)$, implies that
    \[\sqrt{\|\mathbf{v}(T^\sharp)\|_2^2+\|\phi(T^\sharp)\|_5^2}\leq \sqrt{c_2}e^{-c_1 T^{\sharp}/2}\sqrt{\|\mathbf{v}_0\|_2^2+\|\phi_0\|_5^2}\leq \sqrt{c_2}\varepsilon<\frac{\delta_1}{2},\]
    which leads to a contradiction. 
    Hence, $T^\sharp =+ \infty$ and the solution $(\mathbf{v},\phi)$ satisfies \eqref{201708111550n} for any $t>0$.
\hfill$\square$

\subsection{Proof of Theorem \ref{thm-1} under the assumption $(b)$}
We proceed to consider the assumption (b). First of all, we derive the following energy identity:
\begin{lem}
\label{identity}
Let $(\mathbf{v},\phi)$ be a strong solution to the phase-separation problem. Then the following energy identity holds
\begin{align}
&\frac{1}{2}\frac{\mathrm{d}}{\mathrm{d}t} \widetilde{\mathcal{E}}_{\mathrm{d}}
+\widetilde{\mathcal{D}}_{\mathrm{d}} = \widetilde{I}_{\mathrm{d},1}+\widetilde{I}_{\mathrm{d},2},\label{d-energy}
\end{align}
where
\begin{align*}
\widetilde{\mathcal{E}}_{\mathrm{d}}  
& :=\widetilde{\mathcal{E}}_{\mathrm{d},1} (\mathbf{v},\phi) 
+ \widetilde{\mathcal{E}}_{\mathrm{d},1} (\partial_t\mathbf{v},\partial_t\phi)
+ \varrho_1\int \phi |\partial_t\mathbf{v}|^2\,\mathrm{d}x,\\
\widetilde{\mathcal{D}}_{\mathrm{d}}   
&:=\widetilde{\mathcal{D}}_{\mathrm{d},1}(\mathbf{v},\phi) 
+\widetilde{\mathcal{D}}_{\mathrm{d},1}(\partial_t\mathbf{v},\partial_t\phi) ,\\
\widetilde{I}_{\mathrm{d},1}
&:=\int  (\mathbf{M}_{\mathrm{d}} - \varrho_1 \phi   \partial_t\mathbf{v})\cdot   \mathbf{v}\,\mathrm{d}x
+ \int N_{\mathrm{d}}  \phi\,\mathrm{d}x  
+  \frac{l}{2C_{\mathrm{b}} } \int N_{\mathrm{d}}({l}^{-1}  \Psi''(\overline{\varphi})\phi-{l} \Delta  \phi )\,\mathrm{d}x ,\\
\widetilde{I}_{\mathrm{d},2}
&:=\int\Big(\partial_t \mathbf{M}_{\mathrm{d}} -\frac{\varrho_1}{2} \partial_t\phi  \partial_t\mathbf{v} \Big)\cdot  \partial_t\mathbf{v}\,\mathrm{d}x
+ \int\partial_t N_{\mathrm{d}}  \partial_t\phi\,\mathrm{d}x\\
&\quad+  \frac{l}{2C_{\mathrm{b}}} \int \partial_t N_{\mathrm{d}}({l}^{-1}  \Psi''(\overline{\varphi})\partial_t\phi-{l} \Delta  \partial_t\phi)\,\mathrm{d}x,\\
\widetilde{\mathcal{E}}_{\mathrm{d},1} (\mathbf{v},\phi)   
&:=\|  \sqrt{\overline{\rho}} \mathbf{v}\|_{0}^2+ \| \phi\|_{0}^2+ \frac{l^2}{2C_{\mathrm{b}}} \|\nabla  \phi\|^2_{0} 
+\frac{1}{2C_{\mathrm{b}}} \int {\Psi''}(\overline{\varphi})\phi^2\,\mathrm{d}x, \\
\widetilde{\mathcal{D}}_{\mathrm{d},1} (\mathbf{v},\phi) 
& :=\nu \|\nabla \mathbf{v}\|_{0}^2 
+   m\sigma{l}   \|   \Delta  \phi\|^2_{0}+   \frac{m\sigma l}{2C_{\mathrm{b}}} \|\nabla ( {l}^{-1}  \Psi''(\overline{\varphi})\phi-{l} \Delta  \phi)\|_{0}^2\\
 &\quad +  \sigma {l}\int \Big(\partial_{3}^{2}\overline{\varphi}(\partial_{1} \phi \mathbf{v}_1+ \partial_{2}\phi \mathbf{v}_2)-
\partial_{3}\overline{\varphi}\nabla \phi\cdot \nabla \mathbf{v}_3\Big)  \, \mathrm{d}x  + \int \partial_{3}\overline{\varphi}  \phi \mathbf{v}_3 \,\mathrm{d}x
 \\
 &\quad
+ m\sigma  {l}^{-1}\int {\Psi''(\overline{\varphi})}|\nabla \phi|^2 \,\mathrm{d}x
- \frac{  m\sigma  {l}^{-1}}{2} \int \partial_{3}^{2}(\Psi''(\overline{\varphi}))\phi^2  \,\mathrm{d}x\\
&\quad +  \frac{l}{2C_{\mathrm{b}}}\int ( {l}^{-1}\partial_{3}\overline{\varphi}  \Psi''(\overline{\varphi})\mathbf{v}_3
 \phi+  {l}\partial_{3}\overline{\varphi}\nabla \mathbf{v}_3 \cdot \nabla  \phi+  {l}\partial_{3}^{2}\overline{\varphi} \mathbf{v}_3 \partial_{3}\phi)\,\mathrm{d}x,
 \end{align*}
and $\widetilde{\mathcal{E}}_{\mathrm{d},1}
(\partial_t\mathbf{v},\partial_t\phi)$, resp. $\widetilde{\mathcal{D}}_{\mathrm{d},1}(\partial_t\mathbf{v},\partial_t\phi)$
are defined by $\widetilde{\mathcal{E}}_{\mathrm{d},1} (\mathbf{v},\phi) $, 
resp. $\widetilde{\mathcal{D}}_{\mathrm{d},1} (\mathbf{v},\phi)$ with $(\partial_t\mathbf{v},\partial_t\phi)$ in place of $(\mathbf{v},\phi)$. 
\end{lem}
\begin{proof}
To begin with, we derive the  lower-order energy identity.
Taking the inner products of \eqref{g=0-linear}$_3$ and $\phi$, resp. $  {l}^{-1}  \Psi''(\overline{\varphi})
 \phi -{l} \Delta  \phi$ in $L^2$, we have
\begin{align}
&\frac{ 1}{2}\frac{\mathrm{d}}{\mathrm{d}t}\| \phi\|_{0}^2
+  \int \partial_{3}\overline{\varphi}  \phi \mathbf{v}_3 \,\mathrm{d}x
+   m\sigma{l}   \|   \Delta  \phi\|_{0}^2 +   m\sigma  {l}^{-1} \int \Psi''(\overline{\varphi})|\nabla \phi|^2\,\mathrm{d}x\notag\\
&\qquad   -\frac{   m\sigma  {l}^{-1}}{2} \int \partial_{3}^{2}(\Psi''(\overline{\varphi})) \phi^2  \,\mathrm{d}x
=   \int N_{\mathrm{d}}  \phi\,\mathrm{d}x   \label{201708051714xx}
\end{align}
  and
\begin{align}
&\frac{1}{2}\frac{\mathrm{d}}{\mathrm{d}t}\Big({l} \|\nabla  \phi\|^2_{0}+ {l}^{-1}\int
 {\Psi''}(\overline{\varphi})
 \phi^2\,\mathrm{d}x\Big)+ \int ( {l}^{-1}\partial_{3}\overline{\varphi}  \Psi''(\overline{\varphi}) \mathbf{v}_3
 \phi+  {l}\partial_{3}\overline{\varphi}\nabla  \mathbf{v}_3 \cdot \nabla  \phi
 +  {l}\partial_{3}^{2}\overline{\varphi} \mathbf{v}_3 \partial_{3}\phi)\,\mathrm{d}x\nonumber \\
&\qquad+m\sigma  \|\nabla ({l}^{-1}  \Psi''(\overline{\varphi})
 \phi-{l} \Delta  \phi)\|_{0}^2   =   \int N_{\mathrm{d}}( {l}^{-1}  \Psi''(\overline{\varphi})
 \phi-{l} \Delta  \phi)\,\mathrm{d}x. \label{201708051714}
 \end{align}
Multiplying \eqref{201708051714} by ${l}/{2C_{\text{b}}}$ and adding the resultant with   \eqref{v-lower}, \eqref{201708051714xx}, we obtain
  the following lower-order identity of $( \mathbf{v},\phi)$:
\begin{align}
&\frac{1}{2}\frac{\mathrm{d}}{\mathrm{d}t}\widetilde{\mathcal{E}}_{\mathrm{d},1}
( \mathbf{v},\phi)
+\widetilde{\mathcal{D}}_{\mathrm{d},1}( \mathbf{v},\phi)=\widetilde{I}_{\mathrm{d},1}.\label{201708051714n}
 \end{align}
 
Now, taking the inner product of
\eqref{g=0-linear-t}$_1$, \eqref{g=0-linear-t}$_3$, resp. \eqref{g=0-linear-t}$_3$
and $\partial_t \mathbf{v}$,  $ \partial_t \phi$, resp. $ 
  ( \Psi''(\overline{\varphi})
 \partial_t\phi-{l}^2 \Delta  \partial_t\phi   ) /{2C_{\text{b}}}$ in $L^2$, then adding the resultants together, 
 finally we get the following identity for temporal derivatives:
\begin{align}
&\frac{1}{2}\frac{\mathrm{d}}{\mathrm{d}t}
\Big(\widetilde{\mathcal{E}}_{\mathrm{d},1}
(\partial_t\mathbf{v},\partial_t\phi)+\varrho_1\int \phi |\partial_t\mathbf{v}|^2\,\mathrm{d}x\Big)
+\widetilde{\mathcal{D}}_{\mathrm{d},1}(\partial_t\mathbf{v},\partial_t\phi)=\widetilde{I}_{\mathrm{d},2}.\label{201708051714nnn}
 \end{align}
Adding the identities \eqref{201708051714n}  and \eqref{201708051714nnn} together, we obtain the desired energy identity \eqref{d-energy}.
The proof is complete.
\end{proof}

\emph{We are ready to prove Theorem \ref{thm-1} under the assumption (b)}. 
Based on the energy identity \eqref{d-energy}, 
we can follow the idea for the proof of Theorem \ref{thm-1} under the assumption $(a)$. First of all, it is easy to check that
$\widetilde{\mathcal{E}}_{\mathrm{d}}\lesssim {\mathcal{E}}$ and $\widetilde{\mathcal{D}}_{\mathrm{d}}\lesssim\mathcal{D}$ with the assumption \eqref{aprpiosesnew}.
Applying the arguments of the second assertion in Lemma \ref{lem:2070811nn} and Lemma \ref{201708111545}, we also have
\begin{equation*}
\widetilde{\mathcal{E}}_{\mathrm{d}}(t) \lesssim {\|\mathbf{v}(t) \|_{2}^2}+\|\phi (t)\|_{5}^2  \ \mbox{ for all }t\in [0,T),
\end{equation*}
and
$$
|\widetilde{I}_{\mathrm{d},1}|+|\widetilde{I}_{\mathrm{d},2}|\lesssim \sqrt{\mathcal{E}}{\mathcal{D}}.$$
Thus, if we can further  verify that $(\mathbf{v},\phi)$ satisfies
\begin{equation}
\label{201708251047sdf}
{\mathcal{E}}\lesssim\widetilde{\mathcal{E}}_{\mathrm{d}}
\end{equation}
and for sufficiently large $\nu$ or $m$,
\begin{equation}
\label{201708251047nsdf}
  \mathcal{D}\lesssim  \widetilde{\mathcal{D}}_{\mathrm{d}},\quad \text{under the assumption } \eqref{de-assume-1},
\end{equation}
then from \eqref{d-energy}, we can obtain the exponential decay estimate of $(\mathbf{v},\phi)$ as \eqref{201708111550n},
which, together with the local well-posedness result in Proposition \ref{pro:0401n20160203}, implies the conclusion of Theorem \ref{thm-1}. 
Therefore, we shall verify \eqref{201708251047sdf} and \eqref{201708251047nsdf} to complete the proof.

For $i\in\{0,1\}$, using the Cauchy--Schwarz inequality, we get
\begin{align}
  &  \|   \Delta \partial_t^i \phi\|^2_{0}+ \frac{1}{2C_{\text{b}}} \|\nabla ( {l}^{-1}  \Psi''(\overline{\varphi})
 \partial_t^i\phi-{l} \Delta \partial_t^i \phi)\|_{0}^2\notag\\
 &\quad=
 \|   \Delta\partial_t^i  \phi\|^2_{0}
 +\frac{l^2}{2C_{\text{b}}}\|\nabla  \Delta\partial_t^i  \phi \|_{0}^2
 +   \frac{{l}^{-2}}{2C_{\text{b}}} \|\nabla( \Psi''(\overline{\varphi})
 \partial_t^i\phi)\|_{0}^2
 -\frac{1}{C_{\text{b}}}\int \nabla (\Psi''(\overline{\varphi})
 \partial_t^i\phi)\cdot \nabla \Delta\partial_t^i \phi\,\mathrm{d}x\notag\\
 &\quad\geq \|   \Delta \partial_t^i \phi\|^2_{0}+ \frac{ (\epsilon^2-1) (\epsilon l)^{-2} }{2C_{\text{b}}}\|\nabla( \Psi''(\overline{\varphi})
 \partial_t^i\phi)\|_{0}^2+\frac{{l}^2(1-\epsilon^2)}{2C_{\text{b}}}\|\nabla  \Delta \partial_t^i \phi \|_{0}^2,\label{Cb}
 \end{align}
where $\epsilon\in (0,1)$.
In addition, $\|\partial_t^i\phi\|_{2}\lesssim \|\Delta \partial_t^i\phi\|_{0}$ thanks to \eqref{201708241039}. Thus, taking $\epsilon$
sufficiently close to $1$, and using \eqref{201708131702n}, we get
\begin{equation}
 \label{201708251307}
 \| \partial_t^i\phi \|_{3}^2\lesssim     \|   \Delta \partial_t^i \phi\|^2_{0}+  \frac{1}{2C_{\mathrm{b}}}\|\nabla ( {l}^{-1}  \Psi''(\overline{\varphi})
 \partial_t^i\phi-{l} \Delta \partial_t^i \phi)\|_{0}^2.
 \end{equation}
Exploiting the Cauchy--Schwarz inequality, Poincar\'e's inequality \eqref{250508-Poincare} and \eqref{201708251307},
we can deduce that for sufficiently large viscosity $\nu$ or mobility $m$, it holds
\begin{align}
&m\sigma  {l}^{-1}\Big(\int {\Psi''(\overline{\varphi})}|\nabla \phi|^2 \,\mathrm{d}x
 - \frac{ 1}{2} \int \partial_{3}^{2}(\Psi''(\overline{\varphi})) \phi^2  \,\mathrm{d}x
 \Big)+C(\|(\mathbf{v},\partial_t\mathbf{v})\|_{1}^2+\|(\phi,\partial_t\phi)\|_{3}^2)\leq
 \widetilde{\mathcal{D}}_{\mathrm{d}}.\label{201708251837}
 \end{align}
Hence, it follows from the assumption \eqref{de-assume-1} and \eqref{201708251837} that
 \begin{equation}
 \label{201708251845}
\|(\mathbf{v},\partial_t\mathbf{v})\|_{1}^2+\|(\phi,\partial_t\phi)\|_{3}^2 \lesssim
 \widetilde{\mathcal{D}}_{\mathrm{d}} .
 \end{equation}
Finally, by the condition $ \Psi''\geq -C_\text{b}$ in \eqref{2608131333}, we easily find
\begin{equation}
\label{2017081131956safd}
\|(\mathbf{v},\partial_t\mathbf{v})\|_{0}^2+\|(\phi,\partial_t\phi)\|_{1}^2\lesssim \widetilde{\mathcal{E}}_{\mathrm{d}} .
\end{equation}
Consequently, with \eqref{201708251845} and \eqref{2017081131956safd} in hand, we immediately get
\eqref{201708251047sdf} and \eqref{201708251047nsdf} following the argument of \eqref{201708122142xx} and \eqref{201708122142xxD}. 
This completes the proof of Theorem \ref{thm-1} under the assumption $(b)$.\hfill$\square$

\subsection{Proof of Theorem \ref{thm-2}}
Under the assumption \eqref{constant-de}, we find that $\partial_{3}\overline{\varphi}$ and $\Psi'':=\Psi''(\overline{\varphi})$ are both constants. 
This property enables us to derive a refined energy identity, 
where the dissipation functional $\mathcal{D}_\mathrm{d}^\sharp$ does not include any bad terms.

\begin{lem}
\label{lemma 3.6}
   Let $(\mathbf{v},\phi)$ be a strong solution to the phase-separation problem. 
   In addition, assume that \eqref{constant-de} is satisfied. 
   Then we have the following energy identity:
\begin{align}
    \frac{1}{2}\frac{\mathrm{d}}{\mathrm{d}t}\mathcal{E}^\sharp_\mathrm{d}+\mathcal{D}^\sharp_\mathrm{d}
    =I_{\mathrm{d},1}^\sharp+I_{\mathrm{d},2}^\sharp,\label{refined-id}
\end{align}
where
\begin{align*}
   \mathcal{E}_\mathrm{d}^\sharp&:=\|  \sqrt{\overline{\rho}}\mathbf{v}\|_{0}^2+\sigma l\|\nabla\phi\|_{0}^2
   +\|  \sqrt{\overline{\rho}}\partial_t \mathbf{v}\|_{0}^2+\sigma l\|\nabla\partial_t\phi\|_{0}^2
   +\varrho_1\int \phi|\partial_t\mathbf{v}|^2\,\mathrm{d}x,\\
   \mathcal{D}^\sharp_\mathrm{d}&:=\nu( \|\nabla \mathbf{v}\|_{0}^2
   + \|\nabla \partial_t \mathbf{v}\|_{0}^2)+m\sigma ^2 l^2(\|\nabla\Delta\phi\|_{0}^2
   +\|\nabla\Delta\partial_t \phi\|_{0}^2)+m\sigma^2 \Psi''(\|\Delta\phi\|_{0}^2+\|\Delta\partial_t \phi\|_{0}^2),\\
   I_{\mathrm{d},1}^\sharp&:=\int  ({\mathbf{M}}_{\mathrm{d}} - \varrho_1 \phi  \partial_t\mathbf{v})\cdot  \mathbf{v}\,\mathrm{d}x
   -\sigma l\int {N}_\mathrm{d}\Delta\phi\,\mathrm{d}x,\\
   I_{\mathrm{d},2}^\sharp&:=\int  \Big(\partial_t{\mathbf{M}}_{\mathrm{d}} 
   - \frac{\varrho_1}{2}\partial_t \phi  \partial_t\mathbf{v}\Big)\cdot  \partial_t\mathbf{v}\,\mathrm{d}x
   -\sigma l\int \partial_t{N}_\mathrm{d}\Delta\partial_t\phi\,\mathrm{d}x.
\end{align*}
\end{lem}
\begin{proof}
    Since $\partial_{3}\overline{\varphi}$ and $\Psi''$ are constants, the equalities \eqref{rem-1} and \eqref{v-lower} reduce to
    \begin{align}
    &\frac{1}{2}\frac{\mathrm{d}}{\mathrm{d}t}\|\nabla\phi\|_{0}^2+m\sigma l\|\nabla\Delta\phi\|_{0}^2
    +m\sigma l^{-1}\Psi''\|\Delta\phi\|_{0}^2=-\int {N}_\mathrm{d}\Delta\phi\,\mathrm{d}x
    -\partial_{3}\overline{\varphi}\int\nabla\phi\cdot\nabla\mathbf{v}_3\,\mathrm{d}x,\label{rem-1'}
\end{align}
and
\begin{align}
&\frac{1}{2}\frac{\mathrm{d}}{\mathrm{d}t}\|  \sqrt{\overline{\rho}} \mathbf{v}\|_{0}^2
+ \nu \|\nabla \mathbf{v}\|_{0}^2 = \int  ({\mathbf{M}}_{\mathrm{d}} - \varrho_1 \phi  \partial_t\mathbf{v})\cdot  \mathbf{v}\,\mathrm{d}x+  \sigma {l}
\partial_{3}\overline{\varphi}\int  \nabla \phi\cdot \nabla \mathbf{v}_3 \, \mathrm{d}x.\label{v-lower'}
 \end{align}
Multiplying \eqref{rem-1'} by $\sigma l$ and adding the resultant with \eqref{v-lower'}, we get
\begin{align}
    &\frac{1}{2}\frac{\mathrm{d}}{\mathrm{d}t}\Big(\|  \sqrt{\overline{\rho}} \mathbf{v}\|_{0}^2
    +\sigma l\|\nabla\phi\|_{0}^2\Big)+\nu \|\nabla \mathbf{v}\|_{0}^2+m\sigma ^2 l^2\|\nabla\Delta\phi\|_{0}^2
    +m\sigma^2 \Psi''\|\Delta\phi\|_{0}^2\notag\\
    &\quad=\int  ({\mathbf{M}}_{\mathrm{d}} - \varrho_1 \phi  \partial_t\mathbf{v})\cdot  \mathbf{v}\,\mathrm{d}x
    -\sigma l\int {N}_\mathrm{d}\Delta\phi\,\mathrm{d}x.\label{I3d}
\end{align}
Similarly, we can derive the temporal derivative estimate of $(\mathbf{v},\phi)$, that is,
\begin{align}
    &\frac{1}{2}\frac{\mathrm{d}}{\mathrm{d}t}\Big(\|  \sqrt{\overline{\rho}}\partial_t \mathbf{v}\|_{0}^2
    +\sigma l\|\nabla\partial_t\phi\|_{0}^2+\varrho_1\int \phi|\partial_t\mathbf{v}|^2\,\mathrm{d}x\Big)\notag\\
    &\qquad+\nu \|\nabla \partial_t\mathbf{v}\|_{0}^2+m\sigma ^2 l^2\|\nabla\Delta\partial_t\phi\|_{0}^2
    +m\sigma^2 \Psi''\|\Delta\partial_t\phi\|_{0}^2\notag\\
    &\quad=\int  \Big(\partial_t{\mathbf{M}}_{\mathrm{d}} - \frac{\varrho_1}{2}\partial_t \phi  \partial_t\mathbf{v}\Big)\cdot  \partial_t\mathbf{v}\,\mathrm{d}x
    -\sigma l\int \partial_t{N}_\mathrm{d}\Delta\partial_t\phi\,\mathrm{d}x.\label{I4d}
\end{align}
Combining \eqref{I3d} and \eqref{I4d}, we obtain the energy identity \eqref{refined-id}.
\end{proof}

\emph{We are now in a position to prove Theorem \ref{thm-2}}.
We shall verify the equivalence between $\mathcal{E}$, resp. $\mathcal{D}$, and $\mathcal{E}_\mathrm{d}^\sharp$, resp. $\mathcal{D}_\mathrm{d}^\sharp$.
On the one hand, it is easy to check that $\mathcal{E}_\mathrm{d}^\sharp\lesssim\mathcal{E}$ and $\mathcal{D}_\mathrm{d}^\sharp\lesssim\mathcal{D}$.
On the other hand, under the assumption \eqref{aprpiosesnew} with sufficiently small $\delta_1$, 
using the fact $\phi=\partial_t\phi=0$ on $\Gamma$ and Poincar\'e's inequality, we get
\begin{align}
    \|(\mathbf{v},\partial_t\mathbf{v})\|_{0}^2+\|(\phi,\partial_t\phi)\|_{1}^2\lesssim\mathcal{E}^\sharp_\mathrm{d}.\label{E2d}
\end{align}
If $\Psi''\geq0$, we can easily show that
\begin{align}
\|(\mathbf{v},\partial_t\mathbf{v})\|_{1}^2+\|(\phi,\partial_t\phi)\|_{3}^2  \lesssim \mathcal{D}^\sharp_\mathrm{d}.\label{D2d}
\end{align}
If $\Psi''<0$, we can also conclude \eqref{D2d} under the condition \eqref{201708101624nnnn-de}. 
With \eqref{E2d} and \eqref{D2d} in hand, following the proof of Lemma \ref{lem:2070811nn}, we can prove that ${\mathcal{E}}_{\mathrm{d}}^\sharp$
and ${\mathcal{D}}_{\mathrm{d}} ^\sharp$ are equivalent to $ {\mathcal{E}} $
and $\mathcal{D} $, resp. Moreover, we can show that $|I_{\mathrm{d},i}^\sharp|\lesssim\sqrt{\mathcal{E}}\mathcal{D}$ for $i=1,2$, as we did in Lemma \ref{201708111545}. 
Finally, based on the refined energy identity \eqref{refined-id}, 
we can draw the conclusion of Theorem \ref{thm-2} by following the argument for the proof of Theorem \ref{thm-1}.
\hfill$\square$

\section{Phase-Separation Instability}\label{26081414582}

In this section, we present the proof of Theorem \ref{thm-3} concerning the phase-separation instability. 
In addition, we complementally derive an estimate for the largest growth rate.

\subsection{Proof of Theorem \ref{thm-3}}
We prove Theorem \ref{thm-3} in five steps as outlined in Section \ref{sec:2.4}.
First of all, we derive a Gronwall-type energy inequality of the phase-separation problem.
\begin{lem} \label{lem-1}
	Suppose that the assumptions of Theorem \ref{thm-3} are satisfied. 
    Let $(\mathbf{v},\phi)$ be a strong solution to the phase-separation problem. 
    Then, there exists a sufficiently small constant $\delta_1>0$, 
    such that if the solution $(\mathbf{v},\phi)$ further satisfies  the a priori estimate \eqref{aprpiosesnew}, 
    then the following estimate holds:
	\begin{align}
		& \|\mathbf{v}(t)\|_{2}^2 +\|\phi(t)\|_{5}^2+C_1^{\mathrm{d}} \int_0^t\mathcal{D} (\tau)\,\mathrm{d}\tau
        \leq   C_2^{\mathrm{d}}  \Big(\|\mathbf{v}_0\|_{2}^2 
        +\|\phi_0\|_{5}^2+ \int_0^t\| \phi(\tau)\|_{0}^2\,\mathrm{d}\tau \Big),\ \forall\, t\in (0,T).\label{2016121521430850}
	\end{align} 
\end{lem}
\begin{proof}
From the energy equality \eqref{d-energy'} and the definitions of ${\mathcal{D}}_{\mathrm{d}}$, $E_{\mathrm{d}}$, we have
\begin{align}
&\frac{1}{2}\frac{\mathrm{d}}{\mathrm{d}t}{\mathcal{E}}_{\mathrm{d}}+\nu\|\nabla\mathbf{v}\|_{0}^2
+\nu\|\nabla\partial_t\mathbf{v}\|_{0}^2+m\sigma  l(\| \nabla \Delta \phi\|_{0}^2+\| \nabla \Delta \partial_t\phi\|_{0}^2)\notag\\
&\quad=I_{\mathrm{d},1}+I_{\mathrm{d},2}  -m\sigma l^{-1}\Psi''(\alpha)(\|\Delta\phi\|_{0}^2+\|\Delta\partial_t\phi\|_0^2), \label{0107-1}
\end{align}
with $\mathcal{E}_{\mathrm{d}}$, $I_{\mathrm{d},1}$ and $I_{\mathrm{d},2}$ being given as in Lemma \ref{identity'}.
Next, taking the inner product of \eqref{g=0-linear}$_1$ and $\partial_t\mathbf{v}$ in $L^2$, we get
\begin{align}
	&\frac{\nu}{2}\frac{\mathrm{d}}{\mathrm{d}t}\|\nabla\mathbf{v}\|_{0}^2+\|\sqrt{\overline{\rho}}\partial_t\mathbf{v}\|_{0}^2
=\int{\mathbf{M}}_{\mathrm{d}}\cdot\partial_t\mathbf{v}\,\mathrm{d}x-\varrho_1\int\phi|\partial_t\mathbf{v}|^2\,\mathrm{d}x.\label{0107-2}
\end{align}
Similarly, taking the inner product of \eqref{g=0-linear}$_3$ and $-\Delta\partial_t\phi$ in $L^2$ yields
\begin{align}
	&\frac{m\sigma l}{2}\frac{\mathrm{d}}{\mathrm{d}t}\|\nabla\Delta\phi\|_{0}^2+\|\nabla\partial_t\phi\|_{0}^2
    =-\int {N}_{\mathrm{d}}\Delta\partial_t\phi\,\mathrm{d}x
    -m\sigma l^{-1}\Psi''(\alpha)\int \Delta \phi\Delta \partial_t\phi\,\mathrm{d}x .\label{0107-3}
\end{align}
Combining \eqref{0107-1}, \eqref{0107-2} and \eqref{0107-3}, we can deduce that
\begin{align}
 \frac{1}{2}
\frac{\mathrm{d}}{\mathrm{d}t}\widehat{\mathcal{E}}_{\mathrm{d}}
  +\widehat{\mathcal{D}}_{\mathrm{d}} &= \sum_{i=1}^3
I_{\mathrm{d},i }+ \mathcal{R}_{\mathrm{d}},\notag
\end{align}
where
\begin{align*}
\widehat{\mathcal{E}}_{\mathrm{d}}:&= \| \sqrt{\overline{\rho}}\mathbf{v}\|_{0}^2+\nu  \|
\nabla  \mathbf{v} \|_{0}^2
+\|\sqrt{\overline{\rho}}\partial_t \mathbf{v}\|_{0}^2+
 \|\nabla \phi\|_{0}^2 +m\sigma {l}\| \nabla \Delta\phi\|_{0}^2
 +\|\nabla\partial_t \phi\|_{0}^2+ \varrho_1 \int \phi |\partial_t\mathbf{v}|^2\,\mathrm{d}x,\\
  \widehat{\mathcal{D}}_{\mathrm{d}}:&= \nu  \|\nabla \mathbf{v}\|_{0}^2
  +\|\sqrt{\overline{\rho}} \partial_t\mathbf{v}\|_{0}^2 +
 \nu  \|\nabla\partial_t \mathbf{v}\|_{0}^2+ \|\nabla \partial_t\phi\|_{0}^2 +  m\sigma {l}
(\| \nabla \Delta \phi\|_{0}^2+\| \nabla \Delta \partial_t\phi\|_{0}^2 ),\\
I_{\mathrm{d},3 }:&=\int{\mathbf{M}}_{\mathrm{d}}\cdot\partial_t\mathbf{v}\,\mathrm{d}x-\varrho_1\int\phi|\partial_t\mathbf{v}|^2\,\mathrm{d}x
-\int {N}_{\mathrm{d}}\Delta\partial_t\phi\,\mathrm{d}x,\\
\mathcal{R}_{\mathrm{d}}:&=-m\sigma l^{-1}\Psi''(\alpha)\Big(\int \Delta\phi \Delta\partial_t\phi\,\mathrm{d}x
+\|\Delta\phi\|_0^2+\|\Delta\partial_t\phi\|_0^2\Big).\notag
\end{align*}
Similarly to Lemma \ref{201708111545}, we can easily verify that $|I_{\mathrm{d},3}|\lesssim\sqrt{\mathcal{E}}\mathcal{D}$.
For the term $\mathcal{R}_{\mathrm{d}}$, by H\"older's inequality, 
we find $ |\mathcal{R}_{\mathrm{d}}|\lesssim \|\phi\|_{2}^2+\|\partial_t\phi\|_2^2$.
Moreover, as $\|(\mathbf{v},\partial_t\mathbf{v})\|_1^2+\|(\phi,\partial_t\phi)\|_3^2\lesssim\widehat{\mathcal{D}}_{\mathrm{d}}$, 
we can verify the equivalence between $\widehat{\mathcal{D}}_{\mathrm{d}}$ and $\mathcal{D}$, 
following the argument as the proof of Lemma \ref{lem:2070811nn}-(1).
Hence, under the assumption \eqref{aprpiosesnew} with a sufficiently small $\delta_1>0$, 
we infer from Lemma \ref{lem:2070811nn}-(2) and the above estimates that
\begin{align}
 \frac{1}{2}
\frac{\mathrm{d}}{\mathrm{d}t}\widehat{\mathcal{E}}_{\mathrm{d}}
  +C \mathcal{D}  \lesssim \| \phi\|_{2}^2+
\| \partial_t\phi\|_{2}^2.\label{201708070905nnn2118n}
\end{align}
Furthermore, we infer from \eqref{g=0-linear}$_3$ and \eqref{201708132007} that 
\begin{align}
    \|\partial_t\phi\|_{2}\lesssim\|\phi\|_6+\|N_{\mathrm{d}}\|_2   \lesssim \|\phi\|_6+(\|\mathbf{v}\|_{2}+\|\phi\|_4)\|\phi\|_3.\label{251029}
\end{align} 
Putting \eqref{251029} into \eqref{201708070905nnn2118n}, 
by  \eqref{aprpiosesnew}, the interpolation inequality \eqref{201708190532}, and Young's inequality, we obtain
$$
 \frac{\mathrm{d}}{\mathrm{d}t}\widehat{\mathcal{E}}_{\mathrm{d}}
   +C   \mathcal{D}  \lesssim   \| \phi\|_{0}^2.
$$
Since $\mathcal{E}_{\mathrm{d}}\leq \widehat{\mathcal{E}}_{\mathrm{d}}\lesssim\mathcal{E}$,
we infer from Lemma \ref{lem:2070811nn} that $\widehat{\mathcal{E}}_\mathrm{d}$ is equivalent to  $\mathcal{E}$ and  $\widehat{\mathcal{E}}_{\mathrm{d}}|_{t=0}$ satisfies \eqref{20170811}. 
Then, integrating the above inequality with respect to $t$, we obtain the desired Gronwall-type energy inequality \eqref{2016121521430850}.
This completes the proof of Lemma \ref{lem-1}.
\end{proof}

Now, we construct unstable solutions
to the linearized perturbed Cahn--Hilliard equation via the variational method.
\begin{pro}
    \label{exist-unstable-linear}
  Suppose that $g=0$, the assumptions $(\mathbf{A0})$ and $(\mathbf{A1})$ are satisfied, $\overline{\varphi}\equiv\alpha\in(-1,1)$ is a constant equilibrium profile satisfying $\Psi''(\alpha)<0$.
  Under the instability condition \eqref{unstable}, there exists an unstable solution $\phi^{\mathrm{L}}$, 
  which is a strong solution to the linearized Cahn--Hilliard equation
    \begin{align}
       \partial_t\phi^{\mathrm{L}}
       =m\sigma(l^{-1}\Psi''(\alpha)\Delta\phi^{\mathrm{L}}
       -l\Delta^2\phi^{\mathrm{L}})\quad\text{in }\Omega\times \mathbb{R}^+,
        \label{linear-CH}
    \end{align}
    and takes the form
    \[\phi^{\mathrm{L}}=e^{\Lambda_{\mathrm{M}} t}\psi,\]  
    where $\psi\in\mathcal{A}_{\mathrm{D}}\cap H^5$ satisfies
    \begin{align}
        \|\psi\|_{0}\neq0\label{neq0}
    \end{align}
    and the largest growth rate $\Lambda_{\mathrm{M}}>0$ is given by
    \begin{align}
         \Lambda_{\mathrm{M}} =\sup_{\zeta\in \mathcal{A}_{\mathrm{D}}}(-E_\mathrm{d}(\zeta)).\label{Ed}
    \end{align} 
\end{pro}

\begin{proof}
    We consider the variational problem of the energy functional $E_\mathrm{d}$ (cf. \eqref{E-d}). 
    In view of the instability condition \eqref{unstable},  we find that $\Lambda_{\mathrm{M}}>0$. 
    Next, using the interpolation inequality, Poincare's inequality \eqref{250508-Poincare} and the estimate \eqref{201708131702n}, we obtain
    \begin{align*}
        -E_\mathrm{d}(\zeta)&=m\sigma l^{-1}|\Psi''(\alpha)|\|\Delta\zeta\|_{0}^2- m \sigma l\|\nabla\Delta\zeta\|_{0}^2\\
        &\leq  Cm \sigma l^{-1}|\Psi''(\alpha)|\|\zeta\|_{1}\|\zeta\|_{3}- m \sigma l\|\nabla\Delta\zeta\|_{0}^2\\
        &\leq Cm \sigma l^{-1} |\Psi''(\alpha)|\|\nabla\zeta\|_{0}\|\nabla\Delta\zeta\|_{0}-m \sigma l\|\nabla\Delta\zeta\|_{0}^2\\
        &\leq  Cm \sigma {l^{-3}}|\Psi''(\alpha)|^2\|\nabla\zeta\|_{0}^2- \frac{1}{2} m \sigma l\|\nabla\Delta\zeta\|_{0}^2 \\
        &\leq  Cm \sigma l^{-3}|\Psi''(\alpha)|^2,\quad \forall\,\zeta\in\mathcal{A}_{\mathrm{D}},
    \end{align*}
    which implies that $\Lambda_{\mathrm{M}}<+\infty$.
    This means that $-E_\mathrm{d}(\zeta)$ has a maximizing sequence $\{\zeta^n\}_{n\in\mathbb
    {Z}^+}\subset\mathcal{A}_{\mathrm{D}}$ and $\|\zeta^n\|_{3}\leq C$. Thus, there exists a subsequence, still labeled by $\zeta^n$, 
    and a function $\psi\in H_{\mathrm{D}}^3$, such that $\|\nabla\psi\|_{0}=1$ and
    \[\zeta^n\to\psi \text{ weakly in }H^3_\mathrm{D}\text{ and }\text{strongly in }H^2.\]
    It follows that
    \[\Lambda_{\mathrm{M}}=\limsup_{n\to+\infty}(-E_\mathrm{d}(\zeta^n))\leq -E_\mathrm{d}(\psi).\]
    In addition, $-E_{\mathrm{d}}(\psi)\leq \Lambda_{\mathrm{M}}$. Hence, $\psi$ is a maximum point of $\Lambda_{\mathrm{M}}$, that is,
    \begin{align}
        \Lambda_{\mathrm{M}} =\sup_{\zeta\in\mathcal{A}_{\mathrm{D}}}(-E_\mathrm{d}(\zeta))=-E_\mathrm{d}(\psi).\label{maximum-point}
    \end{align}
    By the variational method (i.e., the Lagrange multiplier method), we can deduce from \eqref{maximum-point} the weak formula
         \begin{align}
        \Lambda_{\mathrm{M}}\int \nabla\psi\cdot\nabla\chi\,\mathrm{d}x+m\sigma \Big(l\int \nabla\Delta\psi\cdot\nabla\Delta\chi\,\mathrm{d}x
        +l^{-1}\Psi''(\alpha)\int \Delta\psi\Delta\chi\,\mathrm{d}x\Big)=0,\quad 
       \forall\,\chi\in H^3_{\mathrm{D}}.
        \label{VM-weak}
    \end{align}

Next, we use a bootstrap argument to show that $\psi$ constructed above is indeed a strong solution to the following problem
\begin{align}
    \begin{cases}
        \Lambda_{\mathrm{M}}\psi=m\sigma(l^{-1}\Psi''(\alpha)\Delta\psi-l\Delta^2\psi)&\text{in }\Omega,\\
        \psi=\Delta\psi=0&\text{on }\Gamma.
    \end{cases}\label{BDVP}
\end{align}
To this purpose, we consider the boundary value problem
 \begin{equation}
\label{201604061413nn}
\begin{cases}
  m\sigma
 {l} \Delta  {w}  =F&\text{in }\Omega,\\
   {w}= {0}&\text{on }\Gamma, 
   \end{cases}  
 \end{equation}
with $F=m\sigma  {l}^{-1}\Psi''(\alpha)
\Delta \psi- \Lambda_{\mathrm{M}} \psi\in H_0^1$. Thanks to the existence theory of elliptic problem
with Dirichlet  boundary condition (see, e.g., \cite[Lemma 5.8]{NASII04}), 
there exists a solution $w\in {H}_{\mathrm{D}}^3$ to elliptic problem \eqref{201604061413nn}.
 Multiplying \eqref{201604061413nn}$_1$ by $-\Delta \chi$ and integrating over $\Omega$ yields that
 \begin{equation*}
  \Lambda_{\mathrm{M}}   \int \nabla \psi\cdot\nabla \chi
\,\mathrm{d}x+
 m\sigma\Big(
 {l}\int \nabla w \cdot \nabla  \Delta \chi\,\mathrm{d}x+ {l}^{-1}\Psi''(\alpha)\int
\Delta \psi\Delta \chi\,\mathrm{d}x\Big)=0,\quad\forall\,\chi\in{H}_{\mathrm{D}}^3.
\end{equation*}
Comparing the weak form above with \eqref{VM-weak}, we see that
\begin{equation}
\label{2016040614130919}
\int  \nabla(\Delta \psi- w )\cdot \nabla  \Delta \chi\,\mathrm{d}x=0,\quad\forall\,\chi\in{H}_{\mathrm{D}}^3.
\end{equation}
On the other hand,  using the existence theory of elliptic problem
again, we have a  function $\chi\in {H}_{\mathrm{D}}^3$ satisfying
 \begin{equation}
 \Delta \chi=\Delta \psi- w\in H^1_0 .\notag
  \end{equation}
Putting it into \eqref{2016040614130919} implies  that
\begin{equation}
w= \Delta \psi .\notag
\end{equation}
Recalling that $w\in H^3$, by \eqref{201708131719n}, we have $\psi\in H^5$.
Thus, inserting it into \eqref{201604061413nn}$_1$, we immediately see that
$\psi\in \mathcal{A}_{\mathrm{D}}\cap H^5$ is a strong solution to \eqref{BDVP}.
Finally, since $\|\nabla\psi\|_{0}=1$, it is easy to see that $\|\psi\|_{0}\neq0$, that is, \eqref{neq0} holds. 
Taking $\phi^{\mathrm{L}}:=e^{\Lambda_{\mathrm{M}} t}\psi$, we see that $\phi^{\mathrm{L}}$ satisfies \eqref{linear-CH}. 
The proof of Proposition \ref{exist-unstable-linear} is complete.
 \end{proof}

\begin{rem}
    \label{Neumann-unstable-mode}
    To extend Proposition \ref{exist-unstable-linear} to the Neumann counterpart, 
    we shall introduce the space $\mathcal{A}_{\mathrm{N}}=\{\widetilde{\psi}\in H_{\mathrm{N}}^3: \ \|\nabla\widetilde{\psi}\|_0=1\}$.
    Our aim is to show that, under the assumption of Proposition \ref{exist-unstable-linear} with instability condition $\varpi_{\mathrm{n}}>\varpi$ instead of \eqref{unstable}, 
    there exists an unstable solution $\widetilde{\phi}^{\mathrm{L}}$, which is a strong solution of the linearized Cahn--Hilliard equation \eqref{linear-CH} 
    and enjoys the form $\widetilde{\phi}^{\mathrm{L}}=e^{\widetilde{\Lambda}_{\mathrm{M}}t}\widetilde{\psi}$. 
    Here, $\widetilde{\psi}\in \mathcal{A}_{\mathrm{N}}\cap H_{\mathrm{N}}^5$ satisfies $\|\widetilde{\psi}\|_0\neq 0$, 
    and the largest growth rate $\widetilde{\Lambda}_{\mathrm{M}}>0$ is given by 
    $\widetilde{\Lambda}_{\mathrm{M}}=\sup_{\zeta\in \mathcal{A}_{\mathrm{N}}}(- E_{\mathrm{d}}(\zeta))$. 
    Similar to \eqref{VM-weak}, we can find a function $\widetilde{\psi}\in \mathcal{A}_{\mathrm{N}}$ satisfying
     \begin{align}
        \widetilde{\Lambda}_{\mathrm{M}}\int \nabla \widetilde{\psi}\cdot\nabla\chi\,\mathrm{d}x
        +m\sigma \Big(l\int \nabla\Delta \widetilde{\psi}\cdot\nabla\Delta\chi\,\mathrm{d}x
        +l^{-1}\Psi''(\alpha)\int \Delta \widetilde{\psi}\Delta\chi\,\mathrm{d}x\Big)=0,\quad 
       \forall\,\chi\in H^3_{\mathrm{N}}.
        \notag
    \end{align}
    Since $m\sigma l^{-1}\Psi''(\alpha)\Delta\widetilde{\psi}-\widetilde{\Lambda}_{\mathrm{M}}\widetilde{\psi}\in \underline{H}_0^1$, 
    using the classical theory of elliptic problems with Neumann boundary conditions (see, e.g., \cite[Propositions 7.4 and 7.7]{Taylor2011}) 
    and the bootstrap technique similar to \eqref{BDVP}--\eqref{2016040614130919}, 
    we can show that $\widetilde{\psi}\in \mathcal{A}_{\mathrm{N}}\cap H_{\mathrm{N}}^5$. 
    Finally, since $\|\nabla\widetilde{\psi}\|_{0}=1$, it is easy to see that $\|\widetilde{\psi}\|_{0}\neq0$. 
    This completes the proof.
\end{rem}

By Proposition \ref{exist-unstable-linear}, there exists an unstable profile
\begin{equation}\label{0501}
(\mathbf{v}^\mathrm{L},\phi^\mathrm{L} )=e^{{\Lambda_{\mathrm{M}}  t}} (\bm{0},{\psi} )\in {H}_\sigma^2\times ( H_{\mathrm{D}}^3\cap H^5)
\end{equation}
to the following linearized problem 
\begin{align*}
\begin{cases}
\overline{\rho}\partial_t \mathbf{v}^\mathrm{L}- \nu \Delta \mathbf{v}^\mathrm{L}+\nabla  q^\mathrm{L}
  =\boldsymbol{0}&\text{in }\Omega\times \mathbb{R}^+,\\
\mathrm{div}\,\mathbf{v}^\mathrm{L} =0 &\text{in }\Omega\times \mathbb{R}^+,\\
 \partial_t \phi^\mathrm{L} -m\sigma ({l}^{-1}\Psi''(\alpha)
 \Delta \phi^\mathrm{L}-{l} \Delta^2 \phi^\mathrm{L})=0&\text{in }\Omega\times \mathbb{R}^+,\\
\mathbf{v}^\mathrm{L}=\boldsymbol{0},\ \phi^\mathrm{L}=\Delta\phi^\mathrm{L}=0&\text{on }\Gamma\times \mathbb{R}^+,\\
(\mathbf{v}^\mathrm{L},\phi^\mathrm{L})|_{t=0}=(\boldsymbol{0},{\psi})&\text{in }\Omega, 
\end{cases}
 \end{align*}
with an associated constant pressure $ {q}^\mathrm{L}$. 
Here, $\overline{\varphi}\equiv \alpha\in(-1,1)$ implies that $\overline{\rho}=\alpha\varrho_1+\varrho_2>0$ is constant.
Then the pair of functions \eqref{0501} can produce a family of functions:
$$
(\mathbf{v}^\mathrm{a},\phi^\mathrm{a}):=\delta(\boldsymbol{0},\phi^{\mathrm{L}})\mbox{ for all }\delta>0,
$$
which solves the linearized problem with the initial data
$( \mathbf{v}^\mathrm{a}, \phi^\mathrm{a})|_{t=0}
=\delta (\boldsymbol{0}, {\psi}) $ and satisfies
 \begin{equation}
\label{201602071912MHnnn}
\| \phi^\mathrm{a} (t)\|_{5}\lesssim \delta e^{{\Lambda_{\mathrm{M}}  t}}.
 \end{equation}
However, the initial data of the linear solution $(\mathbf{v}^{\mathrm{a}},\phi^{\mathrm{a}})$ may fail to satisfy the compatibility condition \eqref{compatibility}. 
Indeed, since $\psi=0$ on $\Gamma$, its tangential derivatives vanish, i.e., $\nabla\psi=\partial_{\mathbf{n}}\psi \,\mathbf{n}$ on $\Gamma$. Recalling \eqref{BDVP}, for $\phi_0=\delta\psi$, we find
\[\Delta\widehat{\mu}(0)=\sigma l^{-1}\Delta(\Psi'(\delta\psi+\alpha)-\Psi'(\alpha))-\sigma l\delta \Delta^2\psi=\sigma l^{-1}\delta^2\Psi'''(\alpha)|\partial_{\mathbf{n}}\psi|^2\quad\text{on }\Gamma.\]
We therefore modify the initial data to ensure that the compatibility condition \eqref{compatibility} is satisfied.
 \begin{pro}
 \label{modified-initial-pro}
  Let $\psi$ be the function obtained in Proposition \ref{exist-unstable-linear}.
   Then there exist a constant $\delta_\ast\in (0,1)$ and a function $\phi^{\mathrm{r}}\in {H}_{\mathrm{D}}^3\cap H^5$,
  such that, for every $\delta\in(0,\delta_\ast]$, the modified initial datum
       \begin{align}
           \phi_0^\delta:=\delta\psi+\delta^2\phi^{\mathrm{r}}\in H_{\mathrm{D}}^3\cap H^5\notag
       \end{align}
    satisfies the compatibility condition \eqref{compatibility} with $\phi_0^\delta$ in place of $\phi_0$ and $\|\alpha+\phi_0^\delta\|_{L^\infty}\leq (1+|\alpha|)/2$.
Moreover, the following uniform estimate holds:
       \begin{align}
           \delta^{\frac{1}{2}}\|\phi^{\mathrm{r}}\|_5\leq \min\{\|\psi\|_0,1\}.\label{modified-uni-1}
       \end{align}
 \end{pro}
\begin{proof}
Let $d\in C^5(\overline{\Omega})$ be a regularised distance function for $\Omega$ as in \cite[Proposition III.3.7]{BF}, satisfying $d=0$ and $\nabla d=-\mathbf n$ on $\Gamma$.
Set
\[
\phi^{\mathrm{r}}=\frac{\Psi'''(\alpha)}{24l^2}d^2\psi^2.
\]
Since $H^5$ is an algebra, we have $\phi^{\mathrm{r}}\in H^5$.
Moreover, we infer from $d=\psi=0$ on $\Gamma$ that
\begin{align}
    \phi^{\mathrm{r}}=\Delta \phi^{\mathrm{r}}=0, \quad\nabla \phi^{\mathrm{r}}=\boldsymbol{0}, \quad\text{on } \Gamma,\label{phi^r-bd}
\end{align}
which implies $\phi^{\mathrm{r}}\in H_{\mathrm{D}}^3\cap H^5$.
A direct calculation yields
\[
\begin{aligned}
\Delta^2(d^2\psi^2)\big|_\Gamma
=\big[
8|\nabla d|^2|\nabla\psi|^2
+16(\nabla d\cdot\nabla\psi)^2
\big]\big|_\Gamma=24|\partial_{\mathbf n}\psi|^2,
\end{aligned}
\]
where we have used $\nabla d=-\mathbf n$ and
$\nabla\psi=\partial_{\mathbf n}\psi\,\mathbf n$
on $\Gamma$. Consequently,
\[
\Delta^2\phi^{\mathrm{r}}
=l^{-2}\Psi'''(\alpha)
|\partial_{\mathbf n}\psi|^2
\quad\text{on }\Gamma.
\]
Hence, for $\phi_0^\delta=\delta\psi+\delta^2\phi^{\mathrm{r}}$, using \eqref{phi^r-bd} and the facts $\nabla\psi=\partial_{\mathbf{n}}\psi\,\mathbf{n}$, $\Delta\psi=\Delta^2\psi=0$ on $\Gamma$, we obtain
\begin{align*}
(\Delta\widehat{\mu}(0))|_\Gamma&=\big[\sigma l^{-1}\Delta(\Psi'(\delta\psi+\delta^2\phi^{\mathrm{r}}+\alpha)-\Psi'(\alpha))-\sigma l\delta \Delta^2(\psi+\delta\phi^{\mathrm{r}})\big]|_\Gamma\\
&=\big[\sigma l^{-1}\Delta(\Psi'(\delta\psi+\delta^2\phi^{\mathrm{r}}+\alpha))-\sigma l \delta^2\Delta^2\phi^{\mathrm{r}}\big]|_\Gamma\\
&=\big[\sigma l^{-1}\Psi'''(\alpha)\delta^2|\nabla(\psi+\delta\phi^{\mathrm{r}})|^2+\sigma l^{-1}\delta\Psi''(\alpha)\Delta(\psi+\delta\phi^{\mathrm{r}})-\sigma l \delta^2\Delta^2\phi^{\mathrm{r}}\big]|_\Gamma\\
&=\sigma l^{-1}\Psi'''(\alpha)\delta^2|\partial_{\mathbf{n}}\psi|^2-\sigma l \delta^2(\Delta^2\phi^{\mathrm{r}})|_\Gamma\\
&=0,
\end{align*}
which verifies the compatibility condition \eqref{compatibility}.

Finally, let $C_*=1+\|\phi^{\mathrm{r}}\|_5$ and choose
$0<\delta_*<1$ sufficiently small such that
\[
\sqrt{\delta_*}C_*
\le\min\{\|\psi\|_0,1\},
\qquad
4\delta_*C_{\mathrm e}\|\psi\|_5
\le1-|\alpha|,
\]
where $C_\mathrm{e}>0$ is the embedding constant for $H^2\hookrightarrow L^\infty$.
Then, the estimate \eqref{modified-uni-1} follows immediately. Furthermore,
\[
\|\alpha+\phi_0^\delta\|_{L^\infty}\leq |\alpha|+\delta\|\psi\|_{L^\infty}+\delta^2\|\phi^{\mathrm{r}}\|_{L^\infty}
\leq|\alpha|+2\delta_*C_{\mathrm e}\|\psi\|_5\leq \frac{1+|\alpha|}{2}<1.
\]
This completes the proof of Proposition \ref{modified-initial-pro}.
\end{proof}

Now we use the initial data
\begin{equation}
(\mathbf{v}_0^\delta,\phi_0^\delta ):=(\boldsymbol{0},\delta\psi+\delta^2\phi^{\mathrm{r}} )\notag
\end{equation}
to construct a solution
 to the phase-separation problem for a suitable $0<\delta<\delta_\ast$. 
Since $\|\phi_0^\delta+\alpha\|_{L^\infty}\leq (1+|\alpha|)/2 <1$,
for any given $\delta\in (0,\delta_\ast]$, there exists a unique local strong solution
$$
(\mathbf{v}^\delta ,\phi^\delta  )\in C([0,T^{\max});{H}_\sigma^2\times H^5)
$$
to the phase-separation problem 
with an associated pressure $ {q}^\delta$, emanating
from the initial data $(\mathbf{v}_0^\delta ,\phi_0^\delta)$,
such that $(\mathbf{v}^\delta ,\phi^\delta) $ enjoys the regularity described in Proposition \ref{pro:0401n20160203}.

Define the difference $\phi^{\mathrm{d}} =
\phi^\delta  - \phi^\mathrm{a}$.
Then $\phi^{\mathrm{d}}$
satisfies the following error equation (cf. \eqref{g=0-linear}$_3$ and \eqref{linear-CH})
\begin{align}
 \partial_t\phi^{\mathrm{d}}-m\sigma(l^{-1}\Psi''(\alpha)\Delta\phi^{\mathrm{d}}-l\Delta^2\phi^{\mathrm{d}})=\Phi^\delta,
 \label{newfor12321057}
\end{align}
where
\begin{equation*}
 \Phi^\delta:=
m\sigma  {l}^{-1} \Delta \int_0^1\int_0^1  \Psi'''(s\tau \phi^\delta+\alpha)\tau (\phi^\delta)^2\,\mathrm{d}s\,\mathrm{d}\tau
-\mathbf{v}^\delta\cdot \nabla \phi^\delta
\end{equation*}
subject to the initial and boundary conditions
$$
\phi^{\mathrm{d}}=\Delta \phi^{\mathrm{d}}=0\quad\text{on } \Gamma_T,\qquad
\phi^{\mathrm{d}}|_{t=0}=\delta^2\phi^{\mathrm{r}}\quad\text{in }\Omega.
$$
Moreover, we can deduce the following error estimate from
the error equation \eqref{newfor12321057}.
\begin{lem}\label{lem:0401n}
For any $\delta\in (0,\delta_\ast]$, suppose that on some interval $(0,T)\subset (0,T^{\max})$, it holds
\begin{align}
&\label{201702100940}
\|\mathbf{v}^\delta(t)\|_{2}+
\| \phi^\delta(t)\|_{5} \lesssim\delta e^{\Lambda_{\mathrm{M}} t},
\end{align}
and
\begin{align}
&\label{201702100940'}
\|\phi^\delta(t)\|_5\leq 2\delta_\ast\|\psi\|_5.
\end{align}
Then
there exists a positive constant $C_3^{\mathrm{d}}$ independent of $T$, such that,
\begin{equation}\label{ereroe}
  \| \phi^{\mathrm{d}}(t) \|_{1} 
  \leq \sqrt{C_3^{\mathrm{d}}\delta ^3 e^{3\Lambda_{\mathrm{M}} t}},\quad \forall\,t\in (0,T).
 \end{equation}
\end{lem}
\begin{proof}
Taking the inner product of \eqref{newfor12321057} and $-\Delta\phi^{\mathrm{d}}$ in $L^2$ yields
\begin{align}
 &\frac{1 }{2} \frac{\mathrm{d}}{\mathrm{d}t}
  \|\nabla \phi^{\mathrm{d}}\|_{0}^2
  + {E_{\mathrm{d}}(\phi^{\mathrm{d}} )} 
   =-\int\Phi^\delta\Delta\phi^\mathrm{d}\,\mathrm{d}x.  \label{0917-1}
 \end{align} 
It follows from \eqref{201702100940'} that
$$
\|\phi^\delta(t)+\alpha\|_{L^\infty} 
\leq C_\mathrm{e}\|\phi^\delta(t)\|_{2}+|\alpha|
\leq 2\delta_\ast C_{\mathrm{e}}\|\psi\|_5+|\alpha|
\leq  \frac{1+|\alpha|}{2}<1,\  \forall\, t\in (0,T).
$$
Then we can deduce that
\begin{align}
\|\Phi^\delta\|_0
&\leq m\sigma l^{-1} \left\| \Delta \int_0^1\int_0^1  \Psi'''(s\tau \phi^\delta+\alpha)\tau (\phi^\delta)^2\,\mathrm{d}s\,\mathrm{d}\tau\right\|_0
+\|\mathbf{v}^\delta\cdot \nabla \phi^\delta\|_0 
\notag\\
&\lesssim\|\phi^\delta\|_2^2 +\|\mathbf{v}^\delta\|_1\|\phi^\delta\|_2.
\label{Phi-delta}
\end{align}
According to \eqref{201602071912MHnnn}, \eqref{201702100940}  and \eqref{Phi-delta}, from \eqref{0917-1} we infer that
\begin{align}
  \frac{1 }{2} \frac{\mathrm{d}}{\mathrm{d}t}
  \|\nabla \phi^{\mathrm{d}}\|_{0}^2
  + {E_{\mathrm{d}}(\phi^{\mathrm{d}} )}
  & \leq \|\Phi^\delta\|_{0}\|\phi^\mathrm{d}\|_{2}\leq \|\Phi^\delta\|_0(\|\phi^{\mathrm{a}}\|_2+\|\phi^{\delta}\|_2) 
  \lesssim\delta^3 e^{3\Lambda_{\mathrm{M}} t} . \notag
\end{align}
Thanks to \eqref{Ed}, we further get the following inequality
$$
\frac{\mathrm{d}}{\mathrm{d}t}
  \|\nabla \phi^{\mathrm{d}}\|_{0}^2
\leq 2 \Lambda_{\mathrm{M}} \|\nabla \phi^{\mathrm{d}}\|_{0}^2+
C \delta^3 e^{ 3\Lambda_{\mathrm{M}} t} .
$$
Noting that $\nabla \phi^{\mathrm{d}}|_{t=0}=\delta^2\nabla\phi^{\mathrm{r}}$,
applying Gronwall's lemma, we infer from \eqref{modified-uni-1} that
$$
\|\nabla \phi^{\mathrm{d}}(t)\|_{0}^2
\lesssim e^{2\Lambda_M t}\delta^4 \|\nabla\phi^{\mathrm{r}}\|_0^2+
{\delta^3} e^{3\Lambda_{\mathrm{M}} t}\lesssim {\delta^3} e^{3\Lambda_{\mathrm{M}} t}.
$$
Using Poincar\'e's inequality \eqref{250508-Poincare}, we arrive at the conclusion of \eqref{ereroe}.
 \end{proof}

\emph{Now we are in a position to prove Theorem \ref{thm-3}}. 
Denote
$$\delta_2:=\min\left\{ \frac{\delta_1} { 8\|{\psi}\|_{5} } ,\delta_\ast \right\}  .$$
 Let $\delta\in (0,\delta_2)$, we define
 \begin{align}
&T^*:=\sup\Big\{t\in [0,{T^{\max}}):{\|\mathbf{v}^\delta(\tau)\|_{2}+
\| \phi^\delta (\tau)\|_{5}}\leq 2\delta_2\|{\psi}\|_{5} \text{ for all } \tau\in [0,t]\Big\},\label{2017081815002003} \\
&  T^{**}:=\sup\Big\{t\in [0,T^{\max}):\|\phi^\delta(\tau)\|_{0}
\leq 2{\delta e^{\Lambda_{\mathrm{M}}  \tau}}\|\psi\|_0\text{ for all }\tau\in [0,t]\Big\}.\label{201708181500}
 \end{align}
An argument similar to the one in the proof of Theorem \ref{thm-1} shows that $T^*<T^{\max}$ whenever $T^{\max}<+\infty$.

 Noting that $(\mathbf{v}^\delta, \phi^\delta)\in C([0,T^{\max});{H}_\sigma^2\times H^5)$ and the facts
 \begin{align*}
&{\|\mathbf{v}^\delta(0)\|_{2}+\|\phi^\delta(0)\|_{5}}
= \|\phi_0^\delta\|_{5}\leq\delta\|{\psi}\|_{5}
+\delta^2\|\phi^{\mathrm{r}}\|_5
\leq \delta\|\psi\|_5+\delta^{\frac{3}{2}}\|\psi\|_0<2\delta_2\|{\psi}\|_{5} ,\\
&\|\phi^\delta(0)\|_{0}\leq\delta\|{\psi}\|_{0}
+\delta^2\|\phi^{\mathrm{r}}\|_0
\leq \delta\|\psi\|_0+\delta^{\frac{3}{2}}\|\psi\|_0<2{\delta}\|\psi\|_0,
 \end{align*}
 we have $T^*>0$, $T^{**}>0$.
 By the continuity of $(\mathbf{v}^\delta, \phi^\delta)$, \eqref{2017081815002003} and \eqref{201708181500}, we have
 \begin{align}
 & \|\mathbf{v}^\delta (T^*)\|_{2}+
\| \phi^\delta (T^*)\|_{5} =2\delta_2\|{\psi}\|_{5} \ \mbox{ if }T^*<+\infty,\label{0502n1}\\
&\|\phi^\delta (T^{**}) \|_{0}
=2{\delta e^{\Lambda_{\mathrm{M}}  T^{**}}} \|\psi\|_0\ \mbox{ if }T^{**}< T^{\mathrm{max}}. \label{0502n111}
\end{align}
By \eqref{2017081815002003} and the definition of $\delta_2$, we have
$$
 {\sup_{0\leq t\leq \min\{T^*,T^{**}\}}(\|\mathbf{v}^\delta(t)\|_{2} +\|\phi^\delta(t)\|_{5} )}\leq \frac{\delta_1}{4}.
$$
Thus, by Lemma \ref{lem-1},
 the solution $(\mathbf{v}^\delta,\phi^\delta)$ to the phase-separation problem satisfies the estimate \eqref{2016121521430850}, that is,
$$
\|\mathbf{v}^\delta(t)\|_{2}^2 +\|\phi^\delta(t)\|_{5}^2 
\leq   C_2^{\mathrm{d}}  \Big(\|\mathbf{v}_0^\delta\|_{2}^2 +\|\phi_0^\delta\|_{5}^2+ \int_0^t\| \phi^\delta(\tau)\|_{0}^2\,\mathrm{d}\tau \Big),\ 
\forall\, t\in (0, \min\{T^*,T^{**}\}).
$$
Recalling \eqref{201708181500}, we infer that
for all $t\in
(0,  \min\{T^*,T^{**}\})$, it holds
\begin{align}
\label{20170814}
{\|\mathbf{v}^\delta(t)\|_{2}+\|\phi^\delta(t)\|_{5}} 
\leq  C_4^{\mathrm{d}}\delta{ e^{ \Lambda_{\mathrm{M}}  t}}\mbox{ on }(0,  \min\{T^*,T^{**}\}) ,
\end{align}
where one may take $(C_4^\mathrm{d})^2=4C_2^\mathrm{d}(2\|\psi\|_5^2+\|\psi\|^2_0/\Lambda_{\mathrm{M}})$ by integrating the exponential bound in \eqref{201708181500}.

Let us define
 \begin{align}\label{definedn}
&\epsilon_0:=\min\left\{ \frac{\delta_2\|{\psi}\|_{5}}{ {C_4^{\mathrm{d}} }},
\frac{\|\psi\|_0^2}{4 C_3^{\mathrm{d}} },1 \right\}>0, \\
 \label{times}
& T^\delta:={\Lambda}^{-1}_{\mathrm{M}}\,\mathrm{ln}\Big(\frac{\epsilon_0}{\delta}\Big),\ \mbox{i.e.,}\;\
 \delta e^{\Lambda_{\mathrm{M}} T^\delta }=\epsilon_0,\\
& {T}_{\min}:= \min\{T^\delta ,T^*,T^{**}\}\  \text{and}\ \delta_0:=\min\{\delta_2,\epsilon_0\},\nonumber
 \end{align}
 where the constant $C_3^{\mathrm{d}}>0$ comes from Lemma \ref{lem:0401n} and depends on $C_4^\mathrm{d}$.  
 From now on, we assume $\delta\in (0,\delta_0)$, which ensures that $T^\delta>0$.
Recalling \eqref{2017081815002003} and \eqref{20170814}, we deduce from Lemma \ref{lem:0401n} that
\begin{equation}\label{ereroelast}
\begin{aligned}
  \|   \phi^{\mathrm{d}} (t)\|_{1} \leq \sqrt{C_3^{\mathrm{d}}\delta ^3e^{3\Lambda_{\mathrm{M}} t}} 
   ,\ \ \forall\, t\in (0,{T}_{\min}).
 \end{aligned}
 \end{equation}

We now claim
\begin{equation}
\label{201702092227}
T_{\min}=T^\delta<\min\{T^*,T^{**}\}.
\end{equation}
The conclusion \eqref{201702092227} can be achieved by a contradiction argument. 
First, assume that $T_{\min}=T^*$, then $T^*\leq T^\delta<+\infty$. Noting that
$$
\epsilon_0\leq  \delta_2\|{\psi}\|_{5}/C_4^{\mathrm{d}} ,
$$
then from \eqref{20170814} and \eqref{times}, we infer that
\begin{equation*}
\begin{aligned}
{ {\|\mathbf{v}^\delta (T^*)\|_{2}+
\| \phi^\delta (T^*)\|_{5}}}
\leq C_4^{\mathrm{d}}\delta e^{\Lambda_{\mathrm{M}} T^\ast}
\leq C_4^{\mathrm{d}}\delta e^{\Lambda_{\mathrm{M}} T^\delta}
= C_4^{\mathrm{d}}\epsilon_0<2 \delta_2\|{\psi}\|_{5},
 \end{aligned}
 \end{equation*}
which contradicts \eqref{0502n1}. Hence, $T_{\min}\neq T^*$.
Similarly, if $T_{\min}=T^{**}$, then $T^{**}\leq \min\{T^\delta,T^*\}$. 
We can deduce from \eqref{0501} and \eqref{definedn}--\eqref{ereroelast}, that
  \begin{equation*}
  \begin{aligned}
 \| \phi^\delta (T^{**})\|_{0}
&\leq \|\phi^\mathrm{a}(T^{**})\|_{0}+\|\phi^\mathrm{d} (T^{**})\|_{0}\\
&\leq  {\delta e^{{\Lambda_{\mathrm{M}}  T^{**}}}}\left(\|\psi\|_0 +\sqrt{C_3^{\mathrm{d}}\delta e^{\Lambda_{\mathrm{M}} T^{**}}} \right)\\
&\leq { \delta e^{{\Lambda_{\mathrm{M}}  T^{**}}}}\left(\|\psi\|_0 + \sqrt{C_3^{\mathrm{d}}\epsilon_0} \right)
\leq\frac{3}{2} {\delta e^{\Lambda_{\mathrm{M}}  T^{**}}}\|\psi\|_0<2{\delta e^{\Lambda_{\mathrm{M}}  T^{**}}}\|\psi\|_0,
 \end{aligned}
 \end{equation*}
which contradicts \eqref{0502n111}. 
Hence, $T_{\min}\neq T^{**}$. In summary, the claim \eqref{201702092227} is valid.

Finally,
combining the estimates \eqref{0501}, \eqref{definedn}--\eqref{201702092227}, we can deduce that
 \begin{equation*}
 	\begin{aligned}
 \|\phi^\delta(T^\delta)\|_{0}&\geq  \|\phi^{\mathrm{a}}(T^\delta )\|_{0}-\|\phi^{\mathrm{d}}(T^\delta )\|_{0}\\
& \geq \delta e^{\Lambda_{\mathrm{M}}  T^\delta }\left( \|\psi\|_0- \sqrt{ C_3^{\mathrm{d}}\delta e^{\Lambda_{\mathrm{M}} T^{\delta}}} \right)  
= \epsilon_0 \left( \|\psi\|_0-\sqrt{C_3^{\mathrm{d}}\epsilon_0}\right)\geq  \frac{ \epsilon_0}{2}\|\psi\|_0.
 \end{aligned}
\end{equation*}
This completes the proof of Theorem \ref{thm-3}.
\hfill$\square$

\subsection{Estimate of the largest growth rate}

In this subsection, we provide an estimate on the growth rate $\Lambda_{\mathrm{M}}$ constructed in
 Proposition \ref{exist-unstable-linear}. 
 Here, we call $\Lambda_{\mathrm{M}}$ the largest growth rate of the linearized problem
\begin{align}
\begin{cases}
 \overline{\rho}\partial_t \mathbf{v}- \nu \Delta \mathbf{v}+\nabla  q
  =0& \mbox{ in } \Omega\times \mathbb{R}^+,\\
\mathrm{div}\,\mathbf{v} =0& \mbox{ in } \Omega \times \mathbb{R}^+, \\
\partial_t\phi-m\sigma(l^{-1}\Psi''(\alpha)\Delta\phi-l\Delta^2\phi)=0& \mbox{ in } \Omega\times \mathbb{R}^+,\\
 \mathbf{v}=\boldsymbol{0},\ \ \phi= \Delta \phi=0   &\mbox{ on } \Gamma\times\mathbb{R}^+,\\[1mm]
{(\mathbf{v},\phi )|_{t=0}=( \mathbf{v}_0,\phi_0)}& \mbox{ in } \Omega  ,
\end{cases}\label{201708051653nn}
\end{align}
 since for any solution $(\mathbf{v},\phi)$ to \eqref{201708051653nn}, the following estimate holds
\begin{equation}
\sqrt{\|{\mathbf{v}}(t)\|_{0}^2+\|\nabla \phi(t)\|_{0}^2}
\lesssim e^{\Lambda_{\mathrm{M}} t}\sqrt{\|{\mathbf{v}}_0\|_{0}^2+\|\nabla \phi_0\|_{0 }^2},\quad\forall\,t>0.\notag
\end{equation}
Recalling the definition of $\Lambda_{\mathrm{M}}$, 
the above estimate can be easily obtained by taking the inner products of \eqref{201708051653nn}$_1$, 
resp. \eqref{201708051653nn}$_3$ and $\mathbf{v}$, resp. $-\Delta\phi$ in $L^2$, and applying  Gronwall's lemma.

In what follows, we investigate the effects of parameters $l$, and $\Psi''(\alpha)$ on the largest growth rate $\Lambda_{\mathrm{M}}$. 
Noting that $\Lambda_{\mathrm{M}}$ depends on $l$, we rewrite
 $\Lambda_{\mathrm{M}}$ as $\Lambda_l $, when other
physical parameters in the linearized problem
\eqref{201708051653nn} are given.
Thus, we can regard $\Lambda_l $ as a function of $l\in (0, l_0)$, with $$l_0:=\sqrt{|\Psi''(\alpha)|\varpi_{\mathrm{d}}},$$ 
where $\varpi_{\mathrm{d}}$ is defined in \eqref{201708101624nnnn-de}.
Similarly, we rewrite
 $\Lambda_{\mathrm{M}}$ as $\Lambda_\Upsilon$ to emphasize the dependence of $\Upsilon:=|\Psi''(\alpha)|$, 
 when other physical parameters in the linearized problem \eqref{201708051653nn} are given. 
Then $\Lambda_\Upsilon $ is also a function of $\Upsilon\in (\Upsilon_0,+\infty)$, with $$\Upsilon_0:=l^2\varpi_\mathrm{d}^{-1}.$$
Now, we state the conclusions concerning $\Lambda_l$ and $\Lambda_{\Upsilon}$.

\begin{pro} \label{growth-estimate}
	Under the assumptions of Theorem \ref{thm-3}, the following assertions hold:
     \begin{itemize}
    \item[ $(1)$]  $\Lambda_l\in C^{0,1}_{\mathrm{loc}}(0,l_0)$ and
  strictly decreases on $(0, l_0)$. Moreover,
\begin{equation}
\label{0111nnmmnxxx}
  \Lambda_l \leq   m\sigma{l}^{-5}( l_0^2
 - {l}^2
)
|\Psi''(\alpha)|^2.
\end{equation}
    \item[$(2)$]
    $\Lambda_\Upsilon\in C^{0,1}_{\mathrm{loc}}(\Upsilon_0,+\infty )$ and
   strictly increases on $(\Upsilon_0,+\infty )$. Moreover,
   \begin{align}
       \Lambda_\Upsilon\leq m\sigma l^{-3}\Upsilon(\Upsilon-\Upsilon_0).\label{upper-L}
   \end{align}
  \end{itemize}
 \end{pro}
\begin{proof}
(1) We first show the strict monotonicity  of $\Lambda_l$. 
Let us rewrite ${E_{\mathrm{d}}(\psi )}$ as  ${E_{\mathrm{d}}(\psi,l )}$ to emphasize the dependence on $l$. 
For any $0< l_1 < l_2 < l_0$, by Proposition \ref{exist-unstable-linear},  
there exists a function  $\psi^{l_2}\in \mathcal{A}_{\mathrm{D}}$ such that $\Lambda_{l_2}  = -E_{\mathrm{d}}(\psi^{l_2},l_2).$
Since $\| \nabla \Delta \psi^{l_2}\|_{0}\neq 0$, it holds
\begin{align*}
\Lambda_{l_1}&\geq -E_{\mathrm{d}}(\psi^{l_2} ,l_1)\\
&=-E_{\mathrm{d}}(\psi^{l_2},l_2)+m\sigma \Big((l_2-l_1)\| \nabla \Delta \psi^{l_2}\|_{0}^2 
+ ({l_2}^{-1}- l_1^{-1})\Psi''(\alpha)\|\Delta\psi^{l_2}\|_0^2\Big)\\
&>-E_{\mathrm{d}}(\psi^{l_2},l_2) =\Lambda_{l_2}.
\end{align*}
Hence, $\Lambda_l$ is strictly decreasing with respect to $l\in(0,l_0)$.

Next, we prove the continuity of $\Lambda_l$. Consider an arbitrary bounded interval $[b,l_0)\subset (0, l_0)$. 
For every $l\in [b,l_0)$, there is a function $\psi^l\in \mathcal{A}_{\mathrm{D}}$ satisfying $\Lambda_{l}=-E_{\mathrm{d}}(\psi^{l},l)> 0$ 
due to Proposition \ref{exist-unstable-linear}. Then it holds
\begin{equation}
  \label{201708280915}
  \|\nabla \Delta \psi^{l}\|_{0}^2\leq |\Psi''(\alpha)| {l}^{-2}  \|\Delta\psi^l\|_{0}^2.
\end{equation}
In addition, by integration-by-parts, H\"older's inequality and the fact $\|\nabla \psi^{l}\|_{0}=1$, we obtain
\begin{align}
 \|\Delta\psi^l\|_0^2=-\int \nabla\psi^l\cdot\nabla\Delta\psi^l\,\mathrm{d}x
 \leq \|\nabla\psi^l\|_0 \|\nabla\Delta\psi^l\|_0= \|\nabla\Delta\psi^l\|_0,\label{201708280915n}
\end{align}
Thus, from the estimates \eqref{201708280915}, \eqref{201708280915n} and Young's inequality, we infer that
\begin{equation}
\label{201708280702}
\| \nabla \Delta \psi^{l }\|_{0}^2 +
\|\Delta\psi^l\|_0^2\leq L,  \ \forall\,l\in [b,l_0),
\end{equation}
where the constant $L>0$ depends on $b$, $\Psi$ and $\alpha$, but is independent of $l\in [b,l_0)$.
Exploiting \eqref{201708280702}, we further deduce that, for any $l_1$, $l_2\in [b,l_0)$,
\begin{align*}
\Lambda_{l_1}-\Lambda_{l_2} 
& = E_{\mathrm{d}}(\psi^{l_2},l_2) - E_{\mathrm{d}}(\psi^{l_1},l_1)
\notag \\
&\leq E_{\mathrm{d}}(\psi^{l_1},l_2) - E_{\mathrm{d}}(\psi^{l_1},l_1)\\
&\leq    m\sigma\Big (|l_2-l_1|\| \nabla \Delta \psi^{l_1}\|_{0}^2
+ |{l_2}^{-1}- l_1^{-1}||\Psi''(\alpha)|\|\Delta\psi^{l_1}\|_0^2\Big)\leq
\widetilde{L}| l_2-l_1|,
 \end{align*}
 where $\widetilde{L}:= m\sigma L(1+ b^{-2}|\Psi''(\alpha)|)$.
 Interchanging $l_1$ and $l_2$, we also obtain $\Lambda_{l_2}-\Lambda_{l_1}\leq \widetilde{L}| l_2-l_1|$.
Consequently, we find $\Lambda_l\in C^{0,1}_{\mathrm{loc}}(0,l_0)$.

Finally, we verify the estimate \eqref{0111nnmmnxxx}. 
For any given $l$, there exists a function $\psi^l\in \mathcal{A}_{\mathrm{D}}$ satisfying $\Lambda_{l}=-E_{\mathrm{d}}(\psi^{l},l)$ thanks to Proposition \ref{exist-unstable-linear}. 
By \eqref{201708280915} and \eqref{201708280915n}, we have
$$
\|\nabla \Delta \psi^{l}\|_{0}\leq{l}^{-2}|\Psi''(\alpha)|.
$$
Recalling the definition of $l_0$, we further get
$$
-E_{\mathrm{d}}(\psi^{l},l )
 \leq
m\sigma{l}^{-1}( l_0^2
 - {l}^2
)\|\nabla \Delta \psi^{l}\|_{0}^2\leq  m\sigma{l}^{-5}( l_0^2 - {l}^2 ) |\Psi''(\alpha)|^2,
$$
which yields \eqref{0111nnmmnxxx}.   
This completes the proof of the first assertion in Proposition \ref{growth-estimate}.

(2) The monotonicity and continuity of $\Lambda_{\Upsilon}$ can be established similarly to those in (1), and we omit the details here. 
In what follows, we prove the estimate \eqref{upper-L}.
For any $\Upsilon\in (\Upsilon_0,+\infty)$, there exists a function $\psi^\Upsilon\in \mathcal{A}_{\mathrm{D}}$ 
such that $\Lambda_\Upsilon=-E_{\mathrm{d}}(\psi^\Upsilon,\Upsilon)>0$. 
Then, it holds
\begin{align}
    \|\nabla\Delta\psi^\Upsilon\|_0^2\leq \Upsilon l^{-2}\|\Delta\psi^\Upsilon\|_0^2
    \leq \Upsilon l^{-2}\|\nabla\psi^\Upsilon\|_0\|\nabla\Delta\psi^\Upsilon\|_0
    =\Upsilon l^{-2}\|\nabla\Delta\psi^\Upsilon\|_0,\notag
\end{align}
which implies that $\|\nabla\Delta\psi^\Upsilon\|_0\leq \Upsilon l^{-2}$. 
According to the definitions of $\varpi_{\mathrm{d}}$ and $\Upsilon_0$, 
it holds $\|\nabla\Delta\psi^\Upsilon\|_0^2\geq \Upsilon_0 l^{-2}\|\Delta\psi^\Upsilon\|_0^2$. 
Consequently, we conclude that
\begin{align}
    \Lambda_\Upsilon&=-m\sigma l\|\nabla\Delta\psi^\Upsilon\|_0^2+m\sigma l^{-1}\Upsilon\|\Delta\psi^\Upsilon\|_{0}^2\notag\\
    &\leq -m\sigma l^{-1}\Upsilon_0\|\Delta\psi^\Upsilon\|_0^2+m\sigma l^{-1}\Upsilon\|\Delta\psi^\Upsilon\|_{0}^2\notag\\
    &\leq m\sigma l^{-1}(\Upsilon-\Upsilon_0)\|\nabla\psi^\Upsilon\|_0\|\nabla\Delta\psi^\Upsilon\|_0\notag\\
    &\leq m\sigma l^{-3}\Upsilon(\Upsilon-\Upsilon_0),\notag
\end{align}
which implies \eqref{upper-L}. 
Hence, we complete the proof of Proposition \ref{growth-estimate}.
\end{proof}
\begin{rem}
From \eqref{0111nnmmnxxx}, \eqref{upper-L} and the definition of $\Lambda_{\mathrm{M}}$ we see that
$$
\lim_{l\to l_0^-}\Lambda_{l}=\lim_{\Upsilon\to  \Upsilon_0^+}\Lambda_{\Upsilon}=0,\quad 
\lim_{l\to 0^+}\Lambda_{l}= \lim_{\Upsilon\to  + \infty}\Lambda_{\Upsilon}=+\infty.
$$
To establish the two divergent limits, we take $u=\xi_1/\sqrt{\lambda_1}\in\mathcal A_{\mathrm D}$
in the variational characterization \eqref{Ed} of $\Lambda_{\mathrm M}$ and obtain
\[
\Lambda_{\mathrm M}
\geq -E_{\mathrm d}(\xi_1/\sqrt{\lambda_1})
=m\sigma\lambda_1(\Upsilon l^{-1}-l\lambda_1).
\]
The two divergent limits follow by letting $l\to0^+$
or $\Upsilon\to+\infty$, respectively.
The above facts indicate that a smaller interfacial thickness, or a larger absolute value of the spinodal parameter, 
leads to a stronger linear instability. 
\end{rem}

\section{Stability Results for the CHRT Problem}\label{26081414583}

In this section, we prove the stability results for the CHRT problem 
(i.e., the system \eqref{purterbation} subject to initial-boundary conditions $(\mathbf{D})$) in Theorems \ref{thm:0104n} and \ref{thm:0104nbx}. 
Now, let $(\mathbf{v},\phi)$ be a solution to the CHRT problem that satisfies the estimate \eqref{aprpiosesnew}.  
Similarly to Theorem \ref{thm-1}, the key step is to derive the \emph{a priori} exponential decay estimate of $(\mathbf{v},\phi)$ as in \eqref{201708111550n}, 
provided that $\delta_1$ is sufficiently small.
For this purpose, we rewrite \eqref{purterbation} as the following linearized equations with nonhomogeneous terms
\begin{align}
\begin{cases}
 \overline{\rho}\partial_t \mathbf{v}- \nu \Delta \mathbf{v}+\nabla  q
-\varrho_1 m\xi  \partial_{3} \mathbf{v}+\sigma {l}
  (\partial_{3}\overline{\varphi} \Delta \phi \mathbf{e}_3 
  +\partial_{3}^{2}\overline{\varphi}\nabla \phi)  +g\varrho_1 \phi \mathbf{e}_3
  =\mathbf{M}_{\mathrm{d}}- \varrho_1 \phi \partial_t \mathbf{v},
  \\
\mathrm{div}\,\mathbf{v} =0,
\\
\partial_t \phi+ \partial_{3}\overline{\varphi}\mathbf{v}_3-m\sigma \Delta({l}^{-1}  \Psi''(\overline{\varphi})
 \phi-{l} \Delta  \phi) =N_{\mathrm{d}},
\end{cases}
\label{201708051653nnx1050}
\end{align}
where $\xi$, $\mathbf{M}_\mathrm{d}$ and $N_\mathrm{d}$ are given in the proof of Lemma \ref{identity'}.
 
\subsection{Proof of Theorem \ref{thm:0104n}} 
Comparing the phase-separation problem \eqref{g=0-linear} and the CHRT problem \eqref{201708051653nnx1050}, 
the only difference is an additional gravity term in the momentum equation. 
Therefore, the proof of Theorem \ref{thm:0104n} can be carried out in the same way as the proof of Theorem \ref{thm-1}.

\emph{We first prove Theorem \ref{thm:0104n} under the assumption (a).} 
Similar to \eqref{d-energy'}, we can establish the following energy identity of the CHRT problem:
\begin{align}
&\frac{1}{2}\frac{\mathrm{d}}{\mathrm{d}t} {\mathcal{E}}_{\mathrm{R}}
+{\mathcal{D}}_{\mathrm{R}} = {I}_{\mathrm{d},1}+{I}_{\mathrm{d},2},\notag
\end{align}
where
\begin{align*}
&{\mathcal{E}}_{\mathrm{R}}   :={\mathcal{E}}_\mathrm{d},
\quad {\mathcal{D}}_{\mathrm{R}}   :={\mathcal{D}}_{\mathrm{d}}+\int g\varrho_1\phi\mathbf{v}_3\,\mathrm{d}x
+\int g\varrho_1\partial_t\phi\partial_t\mathbf{v}_3\,\mathrm{d}x
\end{align*}
and ${\mathcal{E}}_\mathrm{d}$, ${\mathcal{D}}_\mathrm{d}$,
${I}_{\mathrm{d},1}$, ${I}_{\mathrm{d},2}$ are given in Lemma \ref{identity'}.
According to Lemmas \ref{lem:2070811nn} and \ref{201708111545}, we already know the facts that
$\mathcal{E}_\mathrm{R}=\mathcal{E}_\mathrm{d}$ {is equivalent to} $\mathcal{E}$, and $I_{\mathrm{d},i}\lesssim\sqrt{\mathcal{E}}\mathcal{D}$ for $i=1,2$.
Furthermore, in analogy to \eqref{20170811319561}, under the assumption \eqref{assume1}, for sufficiently large $\nu$ or $m$, we can deduce that
\begin{equation}
\|(\mathbf{v},\partial_t\mathbf{v})\|_{1}^2+\|(\phi,\partial_t\phi)\|_{3}^2 \lesssim
 {\mathcal{D}}_{\mathrm{R}}.
 \label{DR}
\end{equation}
With the aid of \eqref{DR}, we can prove the equivalence between $\mathcal{D}_\mathrm{R}$ and $\mathcal{D}$ using the argument employed in the proof of Lemma \ref{lem:2070811nn}. 
The rest of the proof proceeds exactly as in Theorem \ref{thm-1} under the assumption (a). 

\emph{We turn to the proof of Theorem \ref{thm:0104n} under the assumption (b).} 
Similar to \eqref{d-energy}, we can establish the following energy identity of the CHRT problem
\begin{align}
&\frac{1}{2}\frac{\mathrm{d}}{\mathrm{d}t} \widetilde{\mathcal{E}}_{\mathrm{R}}
+\widetilde{\mathcal{D}}_{\mathrm{R}} = \widetilde{I}_{\mathrm{d},1}+\widetilde{I}_{\mathrm{d},2},\notag
\end{align}
where
\begin{align*}
&\widetilde{\mathcal{E}}_{\mathrm{R}}   :=\widetilde{\mathcal{E}}_\mathrm{d},
\ \widetilde{\mathcal{D}}_{\mathrm{R}}   :=\widetilde{\mathcal{D}}_{\mathrm{d}}
+\int g\varrho_1\phi\mathbf{v}_3\,\mathrm{d}x
+\int g\varrho_1\partial_t\phi\partial_t\mathbf{v}_3\,\mathrm{d}x
\end{align*}
and $\widetilde{\mathcal{E}}_\mathrm{d}$, $\widetilde{\mathcal{D}}_\mathrm{d}$,
$\widetilde{I}_{\mathrm{d},1}$, $\widetilde{I}_{\mathrm{d},2}$ are given in  Lemma \ref{identity}.
According to the proof of Theorem \ref{thm-1} under the assumption $(b)$, we know that
$\widetilde{\mathcal{E}}_\mathrm{R}=\widetilde{\mathcal{E}}_\mathrm{d}$ 
is equivalent to $\mathcal{E}$ and $\widetilde{I}_{\mathrm{d},i}\lesssim\sqrt{\mathcal{E}}\mathcal{D}$ for $i=1,2$.
Besides, under the assumption \eqref{de-assume-1}, and in analogy to \eqref{Cb}--\eqref{201708251845}, 
we obtain, for sufficiently large $\nu$ or $m$,
\begin{equation}
\|(\mathbf{v},\partial_t\mathbf{v})\|_{1}^2+\|(\phi,\partial_t\phi)\|_{3}^2 \lesssim
 \widetilde{\mathcal{D}}_{\mathrm{R}}.
\notag
\end{equation}
The rest of the proof can be carried out similarly to the proof of Theorem \ref{thm-1} under the assumption (b). 
Hence, we complete the proof of Theorem \ref{thm:0104n}.
\hfill$\square$ 
\subsection{Proof of Theorem \ref{thm:0104nbx}}
Under the assumption \eqref{constant-RT}, it holds that $\partial_{3}\overline{\varphi}>0$ and $\Psi'':=\Psi''(\overline{\varphi})$ are constants. 
This enables us to derive a refined energy identity, in which the dissipation functional $\mathcal{D}_\mathrm{R}^\sharp$ does not include any bad terms.
Then we have:
\begin{lem}
    \label{CHRT-id}
    Let $(\mathbf{v},\phi)$ be a solution to the CHRT problem  and the assumption \eqref{constant-RT} hold. 
    Then, the following energy identity holds:
    \begin{align}
 \frac{1}{2} \frac{\mathrm{d}}{\mathrm{d}t}{\mathcal{E}}^\sharp_{\mathrm{R}}
  +{  \mathcal{D}}^\sharp_{\mathrm{R}} =
I_{\mathrm{R},1}^\sharp+I_{\mathrm{R},2}^\sharp, \label{201708070905xc}
   \end{align}
   where
\begin{align*}
{\mathcal{E}}_{\mathrm{R}} ^\sharp&:= {\partial_{3}\overline{\varphi}} \| \sqrt{\overline{\rho}} \mathbf{v}\|_{0}^2+ E_{\mathrm{R}} (\nabla
\phi  ,\phi ) +
  {\partial_{3}\overline{\varphi}} \| \sqrt{\overline{\rho}} \partial_t\mathbf{v}\|_{0}^2+  E_{\mathrm{R}} (\nabla
\partial_t\phi ,\partial_t\phi) + \varrho_1 \partial_{3}\overline{\varphi} \int \phi |\partial_t\mathbf{v}|^2\,\mathrm{d}x,  \\
{\mathcal{D}}^\sharp_{\mathrm{R}}
&:=
 \partial_{3}\overline{\varphi}\nu (\|\nabla \mathbf{v}\|_{0}^2 +\|\nabla\partial_t \mathbf{v}\|_{0}^2)
 + m\sigma l (E_{\mathrm{R}} (\nabla \Delta
\phi,\Delta \phi)  + E_{\mathrm{R}} (\nabla \Delta\partial_t
\phi,\Delta\partial_t \phi))  \\
&\quad\
  + m\sigma \Psi'' l^{-1}(E_{\mathrm{R}} (\Delta
\phi,\nabla \phi )+  E_{\mathrm{R}} (\Delta
\partial_t\phi,\nabla\partial_t \phi)),\\
I_{\mathrm{R},1}^\sharp&:=\partial_{3}\overline{\varphi}\int ({\mathbf{M}}_{\mathrm{d}}- \varrho_1 \phi  \partial_t\mathbf{v})\cdot  \mathbf{v}\,\mathrm{d}x
-g\varrho_1\int N_\mathrm{d}\phi\,\mathrm{d}x 
- \sigma{l}\partial_{3}\overline{\varphi}  \int N_\mathrm{d} \Delta \phi\,\mathrm{d}x,\\
I_{\mathrm{R},2}^\sharp&:=\partial_{3}\overline{\varphi} \int\Big( \partial_t {\mathbf{M}}_{\mathrm{d}}
  -  \frac{\varrho_1}{2}\partial_t \phi  \partial_t\mathbf{v}\Big)
  \cdot \partial_t\mathbf{v}\,\mathrm{d}x -g\varrho_1\int\partial_t N_\mathrm{d}\partial_t \phi\,\mathrm{d}x 
  - \sigma{l}\partial_{3}\overline{\varphi}  \int \partial_t N_\mathrm{d} \Delta \partial_t \phi\,\mathrm{d}x.
\end{align*}
\end{lem}
\begin{proof}
Under the condition \eqref{constant-RT}, the system
\eqref{201708051653nnx1050} reduces to
\begin{align}
\begin{cases}
 \overline{\rho} \partial_t\mathbf{v}- \nu \Delta \mathbf{v}+\nabla  q
- \varrho_1 m \sigma l^{-1} \Psi''\partial_{3}\overline{\varphi}\partial_3 \mathbf{v}+ \sigma {l}
  \partial_{3}\overline{\varphi} \Delta \phi \mathbf{e}_3   +g\varrho_1 \phi \mathbf{e}_3
   =\mathbf{M}_{\mathrm{d}}- \varrho_1 \phi \partial_t \mathbf{v},\\
\mathrm{div}\,\mathbf{v} =0,\\
\partial_t \phi+ \partial_{3}\overline{\varphi}\mathbf{v}_3-m\sigma( {l}^{-1} \Psi''\Delta \phi-{l} \Delta^2  \phi ) =N_\mathrm{d}.
\end{cases}
\label{201708051653}
\end{align}
Taking the inner product of \eqref{201708051653}$_1$ and $\partial_{3}\overline{\varphi}\mathbf{v}$ in $L^2$ yields
\begin{align*}
&\frac{\partial_{3}\overline{\varphi}}{2}\frac{\mathrm{d}}{\mathrm{d}t}\|  \sqrt{\overline{\rho}}  \mathbf{v}\|_{0}^2
+\partial_{3}\overline{\varphi} \nu  \|\nabla  \mathbf{v}\|_{0}^2
+ \partial_{3}\overline{\varphi}\int( g\varrho_1 \phi  \mathbf{v}_3 -\sigma {l}\partial_{3}\overline{\varphi}
  \nabla \phi\cdot \nabla  \mathbf{v}_3  ) \mathrm{d}x\\
& = \partial_{3}\overline{\varphi}\int({\mathbf{M}}_{\mathrm{d}}- \varrho_1 \phi  \partial_t \mathbf{v})\cdot   \mathbf{v}\,\mathrm{d}x .
\end{align*}
In addition, taking the inner product of  \eqref{201708051653}$_3$ and $-g\varrho_1 \phi$ in $L^2$ yields
\begin{align*}
&-\frac{g\varrho_1 }{2 }\frac{\mathrm{d}}{\mathrm{d}t}\|\phi\|_{0}^2
-g\varrho_1\partial_{3}\overline{\varphi} \int  \phi \mathbf{v}_3 \,\mathrm{d}x-g\varrho_1 m \sigma
( {l}  \| \Delta  \phi\|^2_{0}+ {l}^{-1} \Psi'' \|\nabla \phi \|_{0}^2) = -g\varrho_1\int N_\mathrm{d}\phi\,\mathrm{d}x.
\end{align*}
Thus, adding the above two identities, we get
\begin{align}
&\frac{1}{2}\frac{\mathrm{d}}{\mathrm{d}t}
 \Big(\partial_{3}\overline{\varphi}\| \sqrt{\overline{\rho}} \mathbf{v}\|_{0}^2- {g\varrho_1 }  \|\phi\|_{0}^2 \Big)+
 \partial_{3}\overline{\varphi}\nu  \|\nabla \mathbf{v}\|_{0}^2-g\varrho_1 m\sigma
( {l}  \| \Delta  \phi\|^2_{0} + {l}^{-1} \Psi'' \|\nabla \phi \|_{0}^2)
\nonumber \\
&\qquad= \sigma {l}
  |\partial_{3}\overline{\varphi} |^2\int  \nabla \phi\cdot \nabla \mathbf{v}_3\,  \mathrm{d}x +
  \partial_{3}\overline{\varphi}\int ({\mathbf{M}}_{\mathrm{d}}- \varrho_1 \phi  \partial_t\mathbf{v})\cdot  \mathbf{v}\,\mathrm{d}x
  -g\varrho_1\int N_\mathrm{d}\phi\,\mathrm{d}x.
  \label{rt-low1}
\end{align}
To eliminate the first integral on the right-hand side of \eqref{rt-low1}, 
we take the inner product of \eqref{201708051653}$_3$ and $-\sigma{l} \partial_{3}\overline{\varphi}\Delta \phi $ in $L^2$, which gives
\begin{align*} 
 &\frac{\sigma{l}\partial_{3}\overline{\varphi} }{2}\frac{\mathrm{d}}{\mathrm{d}t}\|\nabla \phi\|_{0}^2
+  m\sigma^2 \partial_{3}\overline{\varphi} ({l}^2\|\nabla \Delta \phi\|^2_{0} + \Psi''
\|\Delta\phi\|_{0}^2)\\
&=  - \sigma{l} |\partial_{3}\overline{\varphi}|^2  \int \nabla \phi \cdot \nabla \mathbf{v}_3 \, \mathrm{d}x
- \sigma{l}\partial_{3}\overline{\varphi}  \int N_\mathrm{d} \Delta \phi\,\mathrm{d}x.
\end{align*}
Adding the above identity to \eqref{rt-low1}, we obtain the following lower-order energy equality of $(\mathbf{v},\phi)$:
\begin{align}
&\frac{1}{2}\frac{\mathrm{d}}{\mathrm{d}t}
 \Big( {\partial_{3}\overline{\varphi}} \| \sqrt{\overline{\rho}} \mathbf{v}\|_{0}^2+  E_{\mathrm{R}} (\nabla
\phi,\phi) \Big) +
 \partial_{3}\overline{\varphi} \nu \|\nabla \mathbf{v}\|_{0}^2 + m\sigma l E_{\mathrm{R}} (\nabla \Delta \phi,\Delta \phi)  
 + m\sigma \Psi'' l^{-1}E_{\mathrm{R}} (\Delta \phi,\nabla \phi)
\notag\\
&\qquad=\partial_{3}\overline{\varphi}\int ({\mathbf{M}}_{\mathrm{d}}- \varrho_1 \phi  \partial_t\mathbf{v})\cdot \mathbf{v}\,\mathrm{d}x
-g\varrho_1\int N_\mathrm{d}\phi\,\mathrm{d}x
- \sigma{l}\partial_{3}\overline{\varphi} \int N_\mathrm{d} \Delta \phi\,\mathrm{d}x,
\label{201708070905cb}
\end{align}
where $E_\mathrm{R}$ is defined as in \eqref{E-R}.

Next, we derive estimates on the temporal derivatives of $(\mathbf{v},\phi)$. Differentiating \eqref{201708051653} with respect to $t$, we obtain
\begin{align}
\begin{cases}
 \partial_t[(\overline{\rho}+\varrho_1\phi) \partial_{t}\mathbf{v}] -\nu \Delta \partial_t\mathbf{v}+\nabla\partial_t q
 -\varrho_1 m \sigma l^{-1} \Psi'' \partial_{3}\overline{\varphi} \partial_{3} \partial_t\mathbf{v}+\sigma {l}
 \partial_{3}\overline{\varphi} \Delta\partial_t \phi \mathbf{e}_3+ g\varrho_1 \partial_t\phi \mathbf{e}_3 \\
 \quad=\partial_t{\mathbf{M}}_{\mathrm{d}},\\
\mathrm{div}\,\partial_t\mathbf{v} =0,\\
\partial_{t}^{2} \phi+ \partial_{3}\overline{\varphi}\partial_t \mathbf{v}_3
-m\sigma (- {l} \Delta^2 \partial_t\phi + {l}^{-1}\Psi''
\Delta\partial_t \phi)
=\partial_tN_\mathrm{d}.
\end{cases}\label{201708252104}
\end{align}
In analogy to the derivation of \eqref{201708070905cb},  we can deduce from \eqref{201708252104} that
\begin{align}
&\frac{1}{2}\frac{\mathrm{d}}{\mathrm{d}t}
 \Big( {\partial_{3}\overline{\varphi}} \| \sqrt{\overline{\rho}}\partial_t \mathbf{v}\|_{0}^2+ E_{\mathrm{R}} (\nabla\partial_t
\phi ,\partial_t\phi) +\varrho_1 \partial_{3}\overline{\varphi} \int \phi |\partial_t\mathbf{v}|^2\,\mathrm{d}x \Big)\notag\\
&\qquad +
 \partial_{3}\overline{\varphi}\nu \|\nabla \partial_t\mathbf{v}\|_{0}^2+ m\sigma l E_{\mathrm{R}} (\nabla \Delta\partial_t
\phi,\Delta \partial_t\phi)+   m\sigma   \Psi'' l^{-1}E_{\mathrm{R}} (\Delta
\partial_t\phi,\nabla\partial_t \phi)
  \nonumber  \\
  &\quad= \partial_{3}\overline{\varphi} \int\Big( \partial_t {\mathbf{M}}_{\mathrm{d}}
  -  \frac{\varrho_1}{2}\partial_t \phi  \partial_t\mathbf{v}\Big)
  \cdot \partial_t\mathbf{v}\,\mathrm{d}x - \int
  \partial_t N_\mathrm{d}
  (g\varrho_1 \partial_t\phi+\sigma {l}\partial_{3}\overline{\varphi} \Delta \partial_t\phi )\,\mathrm{d}x.
 \label{201708070905n}
   \end{align}
Adding \eqref{201708070905cb} and \eqref{201708070905n} together, we arrive at the conclusion \eqref{201708070905xc}.
This completes the proof of Lemma \ref{CHRT-id}.
\end{proof}

\emph{We are now in a position to prove Theorem \ref{thm:0104nbx}.}
In view of the refined energy identity \eqref{201708070905xc}, 
we need to verify the equivalence between $\mathcal{E}_\mathrm{R}^\sharp$, resp. $\mathcal{D}_\mathrm{R}^\sharp$ and $\mathcal{E}$, resp. $\mathcal{D}$.
On the one hand, it is straightforward to check that $\mathcal{E}_\mathrm{R}^\sharp\lesssim\mathcal{E}$ and $\mathcal{D}_\mathrm{R}^\sharp\lesssim\mathcal{D}$.
On the other hand, noting that $\phi=\partial_t\phi=0$ on $\Gamma$,
following the argument for the inhibition estimate \eqref{201708130943n},
we can deduce from the stability condition  \eqref{201708101624} and Poincar\'e's inequality \eqref{250508-Poincare} that
$$
\|\partial_t^i\phi\|_{1}^2\lesssim E_{\mathrm{R}} (\nabla
\partial_t^i\phi  ,\partial_t^i\phi ) \mbox{ for }i=0,1,
$$
which further implies
\begin{equation}
\label{201708252216}
 \|  (\mathbf{v},\partial_t\mathbf{v})\|_{0}^2+\|(\phi,\partial_t\phi)\|_{1}^2\lesssim{\mathcal{E}}^\sharp_{\mathrm{R}},
\end{equation}
provided that $\delta_1$ in the assumption \eqref{aprpiosesnew} is sufficiently small.
Similarly, from the stability condition  \eqref{201708101624}, \eqref{201708131702n}, and \eqref{201708241039}, we can deduce that,
 for $i=0$, $1$,
\begin{equation*}
\|\partial_t^i\phi\|_{2}^2
\lesssim  E_{\mathrm{R}} (\Delta\partial_t^i
\phi ,\nabla \partial_t^i \phi )
\end{equation*}
 and
 \begin{equation}\label{201708190350dsf}
 \|\partial_t^i\phi\|_{3}^2\lesssim  E_{\mathrm{R}} (\nabla \Delta
\partial_t^i \phi,\Delta \partial_t^i \phi). \end{equation}
If $\Psi''\geq 0$, we immediately get
\begin{equation}
\label{201708252214}
\|  (\mathbf{v},\partial_t\mathbf{v})\|_{1}^2+\|(\phi,\partial_t\phi)\|_{3}^2  \lesssim {\mathcal{D}}^\sharp_{\mathrm{R}},
 \end{equation}
For $\Psi''< 0$, we use the condition \eqref{201708101624nnnn} to conclude
$$
  \|(\mathbf{v},\partial_t\mathbf{v})\|_{0}^2
 + E_{\mathrm{R}} (\nabla \Delta
\phi,\Delta \phi )  + E_{\mathrm{R}} (\nabla \Delta\partial_t \phi,\Delta \partial_t\phi)
\lesssim {\mathcal{D}}^\sharp_{\mathrm{R}},
$$
which together with \eqref{201708190350dsf} again implies \eqref{201708252214}. In view of \eqref{201708252216} and \eqref{201708252214},
we can prove the equivalence between $\mathcal{E}^\sharp_{\mathrm{R}}$, resp. $\mathcal{D}_{\mathrm{R}}^\sharp$ and $\mathcal{E}$, resp. $\mathcal{D}$, 
using the argument employed in the proof of Lemma \eqref{lem:2070811nn}.
Moreover, applying the argument of Lemma \ref{201708111545}, we also have $|I_{\mathrm{R},i}^\sharp|\lesssim \sqrt{\mathcal{E}}\mathcal{D}$ for $i=1,2$.
With the above facts in hand, we can complete the proof of Theorem \ref{thm:0104nbx} by adapting the proof of Theorem \ref{thm-1}.
\hfill$\square$


\appendix
\renewcommand{\appendixname}{}
\renewcommand{\thesection}{\Alph{section}}
\section{Admissible Equilibrium Profiles}\label{2608071516}
The results of Theorems \ref{thm-1}--\ref{thm:0104nbx} rely on the choice of suitable admissible equilibrium profiles $(\overline{\varphi},\overline{\mu})$. 
However, establishing the existence of such profiles under general assumptions on the domain, boundary conditions, 
and the nonlinear potential $\Psi$ is a nontrivial problem of independent interest. 
In what follows, we construct some examples for the admissible equilibrium profiles.

In the literature, three types of potential are often used in the Cahn--Hilliard theory for phase separation, that is,
\begin{align}
    &\Psi_{\mathrm{reg}}(r)= C_{\mathrm{R}}(1-r^2)^2,\  r\in \mathbb{R},\nonumber \\
	&\Psi_{\mathrm{log}}(r)=\frac{\Theta}{2}[(1+r)\ln(1+r)+(1-r)\ln(1-r)]-\frac{\Theta_0}{2}r^2,\ r\in[-1,1],\label{2608071527} \\
    &\Psi_{2\mathrm{obs}}(r)=
    \begin{cases}
    C_{\mathrm{OB}}(1-r^2),&\text{if}\ r\in [-1,1],\\
    +\infty,&\text{else},
    \end{cases}\nonumber 
\end{align}
where $C_{\mathrm{R}}$, $C_{\mathrm{OB}}$, $\Theta$ and $\Theta_0$ are positive constants.
For the logarithmic potential $\Psi_{\mathrm{log}}$, we adopt the convention $0\ln 0=0$, so that $\Psi_{\mathrm{log}}$ extends continuously to $[-1,1]$, with
$\Psi_{\mathrm{log}}(\pm1)=\Theta\ln 2-\Theta_0/2$.
Here are the graphs of the above functions: 
\par \quad 

\begin{tikzpicture}
\begin{groupplot}[
    group style={
        group size=3 by 1,
        horizontal sep=1.4cm,
        xlabels at=edge bottom,
        ylabels at=edge left,
        xticklabels at=edge bottom,
        yticklabels at=edge left
    },
    width=5.6cm,
    height=4.9cm,
    xmin=-1.6, xmax=1.6,
    ymin=-1.3, ymax=1.5,
    axis lines=middle,
    axis line style={->, thick},
    tick label style={font=\small},
    label style={font=\small},
    xlabel={$r$},
    ylabel={$y$},
    xtick={-1,0,1},
    ytick={0,1}
]

\nextgroupplot[
    title={$\Psi_{\mathrm{reg}}$},
    legend style={at={(0.97,0.97)}, anchor=north east},
    legend cell align={left}
]

\addplot[blue, thick, domain=-1.5:1.5, samples=200] {0.5*(x^2-1)^2};



\nextgroupplot[
    title={$\Psi_{\mathrm{log}}$},
    legend style={at={(0.97,0.97)}, anchor=north east},
    legend cell align={left}
]

\addplot[blue, thick, domain=-0.999:0.999, samples=300]
    {(1+x)*ln(1+x) + (1-x)*ln(1-x)  - 2*x^2};


\addplot[red, dashed, domain=-1.2:2.2, samples=2] coordinates {(-1,-1.2) (-1,2.2)};
\addplot[red, dashed, domain=-1.2:2.2, samples=2] coordinates {(1,-1.2) (1,2.2)};

\nextgroupplot[
    title={$\Psi_{\mathrm{2obs}}$},
    legend style={at={(0.97,0.97)}, anchor=north east},
    legend cell align={left}
]

\addplot[blue, thick, domain=-1:1, samples=100] {(1-x^2)};



\node[red, above] at (axis cs:-1.3,1.6) {$+\infty$};
\node[red, above] at (axis cs:1.3,1.6) {$+\infty$};

\addplot[red, dashed, domain=-1.2:2.2, samples=2] coordinates {(-1,-1.2) (-1,2.2)};
\addplot[red, dashed, domain=-1.2:2.2, samples=2] coordinates {(1,-1.2) (1,2.2)};

\end{groupplot}
\end{tikzpicture}
\par \quad 

 These potentials are physically relevant and closely related. 
 For instance, $\Psi_{\mathrm{reg}}$ is related to the ``shallow quenching'' limit of $\Psi_{\mathrm{log}}$ by taking $|\Theta-\Theta_0|\ll 1$ 
 and taking Taylor's expansion at $r=0$.
 Moreover, $\Psi_{2\mathrm{obs}}$ is related to the ``deep quenching'' limit of $\Psi_{\mathrm{log}}$ by formally taking $\Theta\to 0^+$ (with a proper constant shift). 
 {Finally, we mention that $\Theta>0$ in \eqref{2608071527} denotes the temperature of the mixture, 
 and $\Theta_0>0$ in \eqref{2608071527} represents the critical temperature below which phase separation processes may occur. 
 Since $\Psi_{\text{log}}'(r)\to\pm\infty$ as $r\to\pm1$, the potential $\Psi_{\text{log}}$ is classified as a singular potential (cf. \cite{CH}).}

In this work, we focus on an equilibrium profile
$\overline{\varphi}$ that is strictly separated from the pure states $\pm 1$, and study its small perturbation
$\phi$ (i.e., the solution to the perturbed AGG system). 
Hence, these solutions are strictly separated from $\pm 1$ and they only ``see'' the behavior of the potential function $\Psi$ away from the pure states $\pm 1$. 
This observation helps us with the construction of admissible equilibrium profiles.

First, we consider the constant equilibrium profile given by
$$
(\overline{\varphi},\overline{\mu})=(a,\sigma l^{-1}\Psi'(a)) \text{ for some}\ a\in (-1,1).
$$
Here, the potential $\Psi$ can be $\Psi_{\mathrm{reg}}$, $\Psi_{\mathrm{log}}$, or $\Psi_{2\mathrm{obs}}$. 
The above equilibrium profile $\overline{\varphi}$ corresponds to a homogeneous mixed state, which is interesting for the phase-separation. 
This profile is admissible if
$(\overline{\varphi},\overline{\mu})$ satisfies the homogeneous Neumann boundary condition \eqref{Neumann'}, 
or the non-homogeneous Dirichlet boundary condition \eqref{Dirichlet'} with the compatibility conditions $\varphi_\mathrm{b}=a$ and $\mu_\mathrm{b}=\sigma l^{-1}\Psi'(a)$ on $\Gamma$.
Finally, let us consider $\Psi=\Psi_{\text{log}}$ and $\overline{\varphi}$ in Theorem \ref{thm-2}, then 
$\overline{\varphi}$ must be a constant under the condition \eqref{constant-de}.
Indeed, $\Psi_{\mathrm{log}}''(r)=\Theta/(1-r^2)-\Theta
_0$ varies on every open interval of phase values, whereas \eqref{quadratic} requires it to be constant on the range of an affine equilibrium.

Next, we investigate the equilibrium profile
$(\overline{\varphi},\overline{\mu})$ that depends only on the $x_3$-variable but is not necessarily constant. 
Consider the potential function $\Psi_{\mathrm{reg}}$ with $C_{\mathrm{R}}=1/4$. 
Recalling 
$$\frac{\mathrm{d}}{\mathrm{d}r}\tanh(r)=1-(\tanh(r))^2,\quad r\in\mathbb{R},$$
we find that the following function
$$
h_1(r)= \tanh\left(\frac{r}{\sqrt{2}l}\right)\in (-1,1), \  r\in \mathbb{R}
$$
satisfies the ordinary differential equation
$$
-\frac{\mathrm{d}^2h_1}{\mathrm{d}r^2} + l^{-2}\Psi_{\mathrm{reg}}'(h_1)=0  \text{ in}\ \mathbb{R}.
$$
Let $\Omega \subset \mathbb{R}^3$ be a smooth bounded domain with boundary $\Gamma$. Define
$$
\overline{\varphi}(x_1,x_2,x_3)= h_1(x_3), \ \overline{\mu}(x_1,x_2,x_3)= 0, \  \forall\, (x_1,x_2,x_3)\in \overline{\Omega}.
$$
We find that the profile $(\overline{\varphi},\overline{\mu})$ defined above is an admissible equilibrium profile under the
Dirichlet boundary condition \eqref{Dirichlet'} with the compatibility conditions $\varphi_\mathrm{b}=\overline{\varphi}$ and $\mu_\mathrm{b}=0$ on $\Gamma$. 
In addition, $\partial_3\overline{\varphi}>0$ holds in $\Omega$ (cf. \eqref{0102}). 
However, for this example, the homogeneous Neumann boundary condition \eqref{Neumann'} cannot be met.

Finally, we look for the profile $\overline{\varphi}$ that satisfies the condition \eqref{constant-de} (or \eqref{constant-RT}). 
Consider the potential function $\Psi_{2\mathrm{obs}}$ with $C_{\mathrm{OB}}=1$. 
We find that the following function
$$
h_2(r)=c r \mbox{ in }\mathbb{R},
$$
with $c>0$ being a given constant, satisfies the ordinary differential equation
$$
-\frac{\mathrm{d}^2h_2}{\mathrm{d}r^2} + l^{-2}\Psi_{2\mathrm{obs}}'(h_2)=-2l^{-2}c r   \text{ in}\ (-c^{-1},c^{-1}).
$$
Let $\Omega \subset \mathbb{R}^3$ be a smooth bounded domain with boundary $\Gamma$. 
Moreover, $\overline{\Omega}\subset \mathbb{R}^2\times (-c^{-1},c^{-1})$. Define
$$
\overline{\varphi}(x_1,x_2,x_3)= h_2(x_3)\in (-1,1), \quad \overline{\mu}(x_1,x_2,x_3)= -2\sigma l^{-1}c x_3, \ \forall\, (x_1,x_2,x_3)\in \overline{\Omega}.
$$
We find that the profile $(\overline{\varphi},\overline{\mu})$ defined above is an admissible equilibrium profile under the
Dirichlet boundary condition \eqref{Dirichlet'} with the compatibility conditions $\varphi_\mathrm{b}=\overline{\varphi}$ and $\mu_\mathrm{b}=\overline{\mu}$ on $\Gamma$.
In addition, $\partial_3\overline{\varphi}=c$ is a positive constant, while $\Psi_{2\mathrm{obs}}''(\overline{\varphi}) = -2$ is also a constant (cf. \eqref{quadratic}). 
However, for this example, the homogeneous Neumann boundary condition \eqref{Neumann'} cannot be met.

\section*{Statements and Declarations}
\noindent \textit{Conflict of interest.} The authors have no competing interests to declare that are relevant to the content of this article.
\medskip

\noindent \textit{Funding.}
The research of Jiang was partially supported by NSFC (Grant Nos. 12525109, 12371233 and 12231016) and the NSF of Fujian Province of China (Grant  No. 2024J011011).
The research of Wu was partially supported by the Natural Science Foundation of Shanghai 25ZR1401023.
\medskip

\noindent \textit{Data availability statement.} Data sharing is not applicable to this article as no datasets were generated during the current study.
\medskip

\noindent \textit{Acknowledgments.} Wu is a member of the Key Laboratory of Mathematics for Nonlinear Sciences (Fudan University), Ministry of Education of China.

 

\end{document}